\documentclass[a4paper,12pt,reqno]{amsart}

\usepackage{fullpage}
\usepackage[dvipsnames]{xcolor}
\usepackage{amsmath,amssymb}
\usepackage{hyperref}
\usepackage{enumerate}
\usepackage{tikz}
\usetikzlibrary{math}

\newtheorem{theorem}{Theorem}[section] 

\newtheorem{algorithm}[theorem]{Algorithm}
\newtheorem{remark}[theorem]{Remark}

\def\H{\widetilde{H}}
\def\MM{\mathcal{M}}

\def\R{{\mathbb R}}
\def\Z{\mathbb {Z}}
\def\N{{\mathbb N}}

\def\KK{{\mathcal K}}
\def\NN{{\mathcal N}}
\def\OO{{\mathcal O}}
\def\PP{{\mathcal P}}
\def\QQ{{\mathcal{Q}}}
 
\def\SS{{\mathcal S}}

\def\XX{{\mathcal X}}

\def\WW{\mathcal{W}}
\def\YY{{\mathcal Y}}

\def\d{\,\mathrm{d}}
\def\V{\mathfrak{V}}
\def\K{\mathfrak{K}}
\def\W{\mathfrak{W}}
\def\diam{{\rm diam}}

\def\norm#1#2{\|#1\|_{#2}}
\def\enorm#1{|\hspace*{-.5mm}|#1|\hspace*{-.5mm}|_{\mathfrak{W}}}
\def\enormV#1{|\hspace*{-.5mm}|#1|\hspace*{-.5mm}|_{\mathfrak{V}}}
\def\seminorm#1#2{\vert #1\vert_{#2}}
\def\set#1#2{\big\{#1\,:\,#2\big\}}
\def\dual#1#2{\langle#1\,,\,#2\rangle}
\def\edualW#1#2{\langle#1\,,\,#2\rangle_{\W}}
\def\edualV#1#2{\langle#1\,,\,#2\rangle_{\V}}

\def\supp{{\rm supp}}

\def\osc{{\rm osc}}

\def\coarse{\star}

\def\fine{+}

\makeatletter
\def\@seccntformat#1{\hspace*{4mm}%
  \protect\textup{\protect\@secnumfont
    \ifnum\pdfstrcmp{subsection}{#1}=0 \bfseries\fi
    \csname the#1\endcsname
    \protect\@secnumpunct
  }%
}
\makeatother

\makeatletter
\def\section{\@startsection{section}{1}%
\z@{.7\linespacing\@plus\linespacing}{.5\linespacing}%
{\normalsize\scshape\bfseries\centering}}
\makeatother

\makeatletter
\renewcommand{\@secnumfont}{\bfseries}
\makeatother

\begin{document}

\address{TU Wien, Institute for Analysis and Scientific Computing, Wiedner Hauptstr.~8--10/E101/4, 1040 Wien, Austria}
\email{\{gregor.gantner (corresponding author), dirk.praetorius, stefan.schimanko\} @asc.tuwien.ac.at}
\author{Gregor Gantner}
\author{Dirk Praetorius}
\author{Stefan Schimanko}

\title{Stable implementation of \\ adaptive IGABEM in 2D in MATLAB}
\date{\today}

\thanks{\textbf{Acknowledgement.}\quad
The authors are supported by the Austrian Science Fund (FWF) through the research projects \textit{Optimal isogeometric boundary element method} (grant~P29096), the doctoral school \textit{Dissipation and dispersion in nonlinear PDEs} (grant~W1245), the special research program \textit{Taming complexity in partial differential systems} (grant SFB F65), and the Erwin Schr\"odinger Fellowship \textit{Optimal adaptivity for space-time methods} (grant~J4379). }

\maketitle

\begin{abstract}
We report on our {\sc Matlab} program package \texttt{IGABEM2D}, which provides an easily accessible implementation of adaptive Galerkin boundary element methods in the frame of isogeometric analysis. 

\end{abstract}


\section{Introduction}

\subsection{Isogeometric analysis}

The central idea of isogeometric analysis (IGA) is to use the same ansatz functions for the discretization of a partial differential equation (PDE), as is used for the representation of the corresponding problem geometry in computer-aided design (CAD), namely (rational) splines.
 This concept, originally invented in~\cite{hughes2005} for finite element methods (IGAFEM) has proved very fruitful in applications; see also the monograph \cite{bible}. 

\subsection{Isogeometric boundary element method}
The idea of the boundary element method (BEM) is to reformulate the considered PDE on a domain $\Omega$ as an equivalent integral equation on the boundary $\partial\Omega$. 
The solution of the latter is the missing Cauchy data, i.e., the Neumann or the Dirichlet data, which can be plugged into the so-called representation formula to obtain the PDE solution itself.
Since CAD usually provides only a parametrization of $\partial \Omega$, this makes BEM the most attractive numerical scheme, if applicable (i.e., provided that the fundamental solution of the differential operator is explicitly known); see \cite{igabem2d,igabem3d} for the first works on isogeometric BEM (IGABEM) for 2D resp.\ 3D. 
Compared to standard BEM, which uses discontinuous or continuous piecewise polynomials as ansatz functions, IGABEM typically saves degrees of freedom by using smooth splines instead. 

We refer {to~\cite{SBTR,helmholtziga,simpson,adss16,tran}} for numerical experiments, to \cite{TM,zechner,dhp16,wolf18,wolf_new,dhkmsw20} for fast IGABEM based on fast multipole, $\mathcal{H}$-matrices resp.\ $\mathcal{H}^2$-matrices, and 
to~\cite{stokesiga,keuchel,acdsss18,acdsss18,fgkss18} for some quadrature analysis.
For further references, we also mention the recent monograph~\cite{bmd20}.

\subsection{Adaptive IGABEM}

On the one hand, IGA naturally leads to high-order ansatz functions. On the other hand, however, optimal convergence behavior with higher-order discretizations is only observed in simulations, if the (given) data as well as the (unknown) solution are smooth. Therefore, {\sl a~posteriori} error estimation and related adaptive strategies are mandatory to exploit the full potential of IGA. Rate-optimal adaptive strategies for IGAFEM (using hierarchical splines) have been proposed and analyzed independently in~\cite{bg2017,ghp2017} for IGAFEM, while the earlier work~\cite{bg2016} proves only linear convergence. 
Meanwhile, these results have been extended to T-splines~\cite{gp20a}. 

As far as IGABEM is concerned, available results focus on weakly-singular integral equations with energy space $H^{-1/2}(\Gamma)$; see~\cite{igafaermann,resigabem} for {\sl a~posteriori} error estimation as well as~\cite{optigabem} for the analysis of a rate-optimal adaptive IGABEM in 2D, and~\cite{gantner17,gp20b,gp21} and \cite{bggpv21} for corresponding results for IGABEM in 3D with hierarchical splines and T-splines, respectively.
Recently, \cite{fgps19} investigated optimal preconditioning for IGABEM in 2D with locally refined meshes.
Moreover, \cite{gps20}  proved rate-optimality of an adaptive IGABEM in 2D for hyper-singular integral equations with energy space $H^{1/2}(\Gamma)$.
The recent work \cite{bggpv21} provides an overview of adaptive IGAFEM as well as IGABEM.

\subsection{Model problems}
\label{sec:model}
The present paper reports on our {\sc Matlab} implementation of adaptive IGABEM in 2D, which is available online~\cite{gps22}.
On a bounded Lipschitz domain $\Omega\subset\R^2$  with $\diam(\Omega)<1$, we consider the Laplace model problem $-\Delta P = 0$ with Dirichlet/Neumann boundary conditions.
Rewriting this PDE as boundary integral equation involves (up to) four integral operators, namely the weakly-singular operator $\mathfrak{V}$, the hyper-singular operator $\mathfrak{W}$, the double-layer operator $\mathfrak{K}$, and the adjoint double-layer operator $\mathfrak{K}'$.
Let $\Gamma\subseteq\partial\Omega$ be a connected part of the boundary. 
Given boundary densities $\phi,u:\Gamma\to \R$ and $x\in\Gamma$, the operators are formally defined by 
\begin{align}
\label{eq:V}
[\V\phi](x) &:= -\frac{1}{2\pi} \, \int_\Gamma \log|x-y| \, \phi(y) \d{y},
\\
[\W u](x) &:= \frac{1}{2\pi} \, \frac{\partial_x}{\partial \nu(x)} \, \int_\Gamma \frac{\partial_y}{\partial \nu(y)} \log|x-y| \, u(y) \d{y},
\\ \label{eq:K}
 [\K u](x) &:= - \frac{1}{2\pi} \, \int_\Gamma \frac{\partial_y}{\partial \nu(y)} \log|x-y| \, u(y) \d{y}, 
\\
 [\K'\phi](x) &:= -\frac{1}{2\pi}\int_\Gamma \phi(y)\frac{\partial_x}{\partial\nu(x)}\log|x-y|\d y.
\end{align}
For our implementation, we consider the corresponding weakly-singular integral equation
\begin{align}\label{eq:weak strong}
 [\V\phi](x) :
 = f(x)
 \quad \text{for } x \in \Gamma,
\end{align}
where $f : \Gamma \to \R$ is given and the integral density $\phi : \Gamma \to \R$ is sought. Moreover, we consider the hyper-singular integral equation
\begin{align}\label{eq:hyper strong}
 [\W u](x) 
 = g(x)
 \quad \text{for } x \in \Gamma. 
\end{align}
where $\nu$ denotes the outer unit normal vector on $\partial\Omega$, $g : \Gamma \to \R$ is given, and the integral density $u : \Gamma \to \R$ is sought. 
For $\Gamma = \partial \Omega$, 
these integral equations are equivalent formulations of the 2D Laplace problem $-\Delta P = 0$ in $\Omega$, where~\eqref{eq:weak strong} with $f = (1/2 + \K)(P|_\Gamma)$ is equivalent to the Dirichlet problem and~\eqref{eq:hyper strong} with $g = (1/2 - \K')(\partial P/\partial\nu)$ is equivalent to the Neumann problem. 
Knowing both Dirichlet as well as Neumann boundary data allows to compute $P$ via the representation formula 
$P = \widetilde{\mathfrak{V}}(\partial P/\partial\nu) - \widetilde{\mathfrak{K}}(P|_\Gamma)$, where $\widetilde{\mathfrak{V}}$ and $\widetilde{\mathfrak{K}}$ are defined as in \eqref{eq:V} and \eqref{eq:K} for $x\in\Omega$. 
While this approach is called direct, an alternative indirect approach is to make the ansatz $P=\widetilde{\mathfrak{V}}\phi$ or $P=\mathfrak{K}u$ and solve $\widetilde{\mathfrak{V}}\phi=f:=P|_\Gamma$ or $\mathfrak{W}u=g:=\partial P/\partial\nu$, respectively.
In all cases, the singularities of the involved integral kernels pose a significant challenge, which is overcome in the presented stable IGABEM implementation.

\subsection{ {\sc Matlab} library {\tt IGABEM2D} and contributions}
The library {\tt IGABEM2D} is available online~\cite{gps22}. It is distributed as a single file {\tt igabem2d.zip}. To install the library, download and unpack the zip archive, start {\sc Matlab}, and change to the root folder {\tt igabem2d/}.
You can directly run the main files  \texttt{IGABEMWeak.m} and \texttt{IGABEMHyp.m}, where the adaptive Algorithm~\ref{the algorithm} and~\ref{the algorithm2} for the problems of Section~\ref{sec:model} are implemented.
There, you may choose out of a variety of different parameters. 
These functions both automatically run \texttt{mexIGABEM.m}, which adds all required paths, compiles the {\sc C} files of \texttt{source/} and stores the resulting {\sc mex}-files in \texttt{MEX-files/}. 
The {\sc Matlab} functions provided by {\tt IGABEM2D} are contained in the folder {\tt functions/}.
Further details are also provided by {\tt README.txt}, which is found in the root folder {\tt igabem2d/}.
\subsection{Outline}
Section~\ref{sec:iga} recalls the definition of B-splines and NURBS and discusses necessary assumptions on parametrized NURBS geometries. 
Section~\ref{sec:weak} recalls the weakly-singular integral equation along with the required Sobolev spaces on the boundary and formulates an adaptive Galerkin IGABEM, which is applied in numerical experiments at the end of the section. 
Similarly, Section~\ref{sec:hyp} covers the hyper-singular integral equation. 
Details on the implementation, i.e., the stable computation of the involved (singular) boundary integrals, are given in Section~\ref{sec:implementation}.
Finally, Appendix~\ref{sec:overview} provides an overview of all functions in \texttt{IGABEM2D}.

\subsection{General notation}
Throughout and without any ambiguity, $|\cdot|$ denotes the absolute value of scalars, the Euclidean norm of vectors in $\R^2$, the measure of a set in $\R$ (e.g., the length of an interval), or the arclength of a curve in $\R^{2}$.
Throughout, mesh-related quantities have the same index, e.g., $\NN_\coarse$ is the set of nodes of the partition $\QQ_\coarse$, and $h_\coarse$ is the corresponding local mesh-width function etc. 
The analogous notation is used for partitions $\QQ_\fine$,\ $\QQ_\ell$, etc. We sometimes use  \,$\widehat{\cdot}$\, to transform notation on  the boundary to the  parameter domain or vice versa, e.g., $\widehat{\QQ}_\star$ is the partition of the parameter domain related to the partition $\QQ_\star$ of the boundary.


\section{Isogeometric analysis on 1D boundaries}\label{sec:iga} 

In this section, we recall the definition of B-splines and NURBS and discuss necessary assumptions on parametrized NURBS geometries. 
This provides the basis for the implemented IGABEM discretization discussed in Section~\ref{sec:weak}--\ref{sec:hyp}.

\subsection{B-Splines and NURBS in 1D}\label{sec:splines}
We consider an arbitrary but fixed sequence {$\widehat{\mathcal{K}}_\coarse:=(t_{\coarse, i})_{i\in\Z}$} on $\R$ with multiplicities $\#_\coarse t_{\coarse, i}:=\#\set{i'\in\Z}{t_{\coarse,i'} = t_{\coarse,i}}$ which satisfy that
$t_{\coarse, i-1}\leq t_{\coarse, i}$ for $i\in \Z$ and $\lim_{i\to \pm\infty}t_{\coarse, i}=\pm \infty$.
Let $\widehat{\mathcal{N}}_\coarse:=\set{t_{\coarse, i}}{i\in\Z}=\set{\widehat x_{\coarse, j}}{j\in \Z}$ denote the corresponding set of nodes with $\widehat x_{\coarse, j-1}<\widehat x_{\coarse, j}$ for $j\in\Z$.
For $i\in\Z$, the $i$-th \textit{B-spline}~\cite{boor86} of degree $p\in\N_0$ is defined for $t\in\R$ inductively by
\begin{align}
\begin{split}
\widehat B_{\coarse, i,0}(t)&:=\chi_{[t_{\coarse, i-1},t_{\coarse, i})}(t),\\
\widehat B_{\coarse, i,p}(t)&:=\frac{t-t_{\coarse, i-1}}{t_{\coarse, i-1+p}-t_{\coarse, i-1}}  \widehat B_{\coarse, i,p-1}(t)
+\frac{t_{\coarse, i+p}-t}{t_{\coarse, i+p}-t_{\coarse, i}} \widehat B_{\coarse, i+1,p-1}(t) 
\quad \text{for } p\in \N,
\end{split}
\end{align}
where we suppose the convention $(\cdot)/0:=0$.
It is well known that for arbitrary $I=[a,b)$ and $p\in\N_0$, the set $\set{\widehat B_{\coarse, i,p}|_I}{i\in \Z\wedge \widehat B_{\coarse, i,p}|_I\neq 0}$
is a basis for the space of all right-continuous piecewise polynomials of degree lower or equal $p$ with breakpoints $\widehat{\mathcal{N}}_\coarse$ on $I$ 
which are, at each knot $t_{\coarse, i}$, $p-\#_\coarse t_{\coarse, i}$ times continuously differentiable if $p-\#_\coarse t_{\coarse, i}\ge 0$.
Moreover, the B-splines are non-negative and have local support 
\begin{align}\label{eq:locality}
\supp(\widehat B_{\coarse, i,p})=[t_{\coarse, i-1},t_{\coarse, i+p}] \quad\text{for all }i\in\Z.
\end{align}
Indeed, $\widehat B_{\coarse, i,p}$ depends only on the knots $t_{\coarse, i-1},\dots,t_{\coarse, i+p}$. 
If $t_{\coarse, i-1}<t_{\coarse, i+p}$, then 
\begin{align}\label{eq:interpolatoric}
t_{\coarse, i}=\dots=t_{\coarse, i+p-1}\quad\Longrightarrow \quad \widehat B_{\coarse, i,p} (t_{\coarse, i}-)=1= \widehat B_{\coarse, i,p} (t_{\coarse, i}). 
\end{align}
Moreover, B-splines form a partition of unity, i.e., $\sum_{i\in\Z} B_{\coarse,i,p} = 1$ for all $p\in\N_0$.
For $p\ge 1 $ and $i\in\Z$, the right derivative of the corresponding B-spline can be computed as 
\begin{equation}\label{eq:derivative of splines}
\widehat B_{\star,i,p}'^{_r}=\frac{p}{t_{\star,i+p-1}-t_{\star,i-1}} \widehat B_{\star,i,p-1}-\frac{p}{t_{\star,i+p}-t_{\star,i}}\widehat B_{\star,i+1,p-1}.
\end{equation}
Finally, if $\sum_{i\in \Z} a_{\coarse,i} \widehat B_{\coarse,i,p}$ is a given spline and $\widehat\KK_\fine$ is obtained from $\widehat\KK_\coarse$ by adding a single knot $t'$ with $t'\in(t_{\coarse,\ell-1},t_{\coarse,\ell}]$ for some $\ell\in\Z$,
there exist coefficients $(a_{\fine,i})_{i\in\Z}$ such that
\begin{equation}\label{eq:insertion}
\sum_{k\in \Z} a_{\coarse,k} \widehat B_{\coarse,k,p}=\sum_{k\in \Z} a_{\fine,k} \widehat B_{\fine,k,p}.
\end{equation}
With the multiplicity $\#_{\fine}t'$ of $t'$ in the knots $\widehat\KK_{\fine}$, the new coefficients can be chosen as convex combinations of the old coefficients
\begin{align}
\begin{split}
a_{\fine,k}=\begin{cases}
a_{\coarse,k} \quad &\text{if } k\leq \ell-p+\#_{\fine} t'-1,\\
\frac{t_{\coarse,k-1+p}-t'}{t_{\coarse,k-1+p}-t_{\coarse,k-1}} a_{\coarse,k-1} +\frac{t'-t_{\coarse,k-1}}{t_{\coarse,k-1+p}-t_{\coarse,k-1}} a_{\coarse,k} \quad &\text{if } \ell-p+\#_{\fine} t'\leq k\leq \ell,\\
a_{\coarse,k-1}\quad &\text{if } \ell+1 \leq k.
\end{cases}
\end{split}
\end{align}
If $\#_{\coarse} t_k\le p+1$ for all $k\in\Z$, these coefficients are unique.
Proofs are found, e.g., in \cite{boor86}.

In addition to the knots $\widehat{\mathcal{K}}_\coarse=(t_{\coarse, i})_{i\in\Z}$, we consider fixed positive weights {$\mathcal{W}_\coarse:=(w_{\coarse, i})_{i\in\Z}$} with $w_{\coarse, i}>0$.
For $i\in \Z$ and $p\in \N_0$, we define the $i$-th 
NURBS (non-uniform rational B-spline) by
\begin{equation}
\widehat R_{\coarse, i,p}:=\frac{w_{\coarse, i} \widehat B_{\coarse, i,p}}{\sum_{k\in\Z}  w_{\coarse, k}\widehat B_{\coarse, k,p}}.
\end{equation}
Note that the denominator is locally finite and positive.
For any $p\in\N_0$, we define the spline space as well as the rational spline space 
\begin{equation}\label{eq:NURBS space defined} 
\widehat\SS^p(\widehat{\mathcal{K}}_\coarse):={\rm span} \,
\set{\widehat B_{\coarse, i,p}}{i\in\Z}\quad\text{and}\quad
\widehat\SS^p(\widehat{\mathcal{K}}_\coarse, \mathcal{W}_\coarse):={\rm span} \,
\set{\widehat R_{\coarse, i,p}}{i\in\Z}.
\end{equation}

\subsection{Boundary parametrization}
\label{sec:boundary parametrization}
Recall that $\Omega\subset\R^2$ is a bounded Lipschitz domain with $\diam(\Omega)<1$ and $\Gamma\subseteq\partial\Omega$ is a connected part of its boundary. 
We further assume that either $\Gamma=\partial\Omega$ is parametrized by a closed continuous and
piecewise continuously differentiable path $\gamma:[a,b]\to\Gamma$ such
that the restriction $\gamma|_{[a,b)}$ is even bijective, or that $\Gamma\subsetneqq\partial\Omega$ is parametrized by a bijective continuous and piecewise continuously differentiable path $\gamma:[a,b]\to\Gamma$.  
In the first case, we speak of \textit{closed} $\Gamma=\partial\Omega$, whereas the second case is referred to as \textit{open} $\Gamma\subsetneqq\partial\Omega$.
Throughout and by abuse of notation, we write $\gamma^{-1}$ for the inverse of $\gamma|_{[a,b)}$ and $\gamma|_{(a,b]}$.
The meaning will be clear from the context.
For the left- and right-hand derivative of $\gamma$, we assume that {$\gamma^{\prime_\ell}(t)\neq 0$ for $t\in(a,b]$ and $\gamma^{\prime_r}(t)\neq 0$  for $t\in [a,b)$.}
Moreover, we assume that $\gamma^{\prime_\ell}(t)
+c\gamma^{\prime_r}(t)\neq0$ for all $c>0$ and $t\in[a,b]$ and $t\in(a,b)$, respectively.
Note that these assumptions provide a pointwise representation of the arc-length derivative
\begin{align}\label{eq:explicit derivative}
(\partial_\Gamma u)\circ\gamma = \frac{(u\circ\gamma)'}{|\gamma'|}\quad\text{for all }u\in H^1(\Gamma). 
\end{align}
Finally and without loss of generality, we suppose that $\gamma$ is positively oriented in the sense that the outer normal vector of $\Omega$ has the form
\begin{align}
\nu(\gamma(t))=
\frac{1}{|\gamma'(t)|}\begin{pmatrix}\gamma_2(t)\\ -\gamma_1(t)\end{pmatrix}
\quad\text{for almost all }t\in[a,b].
\end{align}

While the given assumptions are sufficient for the abstract analysis of adaptive BEM, the later implementation requires that $\gamma$ is a NURBS curve of degree $p_\gamma\in\N$ in the following sense:
Let $\widehat\KK_\gamma=(t_{\gamma,i})_{i=1}^{N_\gamma}$ be a sequence of knots with 
\begin{align}
a<t_{\gamma,1}\le\dots\le  t_{\gamma,N_\gamma-p_\gamma-1}< t_{\gamma,N_\gamma-p_\gamma}=\dots= t_{\gamma,N_\gamma}=b
\end{align}
 and multiplicity $\#_\gamma t_{\gamma, i}\le p_\gamma$ for $i\in\{1,\dots,N_\gamma-p_\gamma\}$.
Let 
\begin{align}
\widehat{\mathcal{N}}_\gamma:=\set{t_{\gamma, i}}{i\in\{1,\dots N_\gamma\}}=\set{\widehat x_{\gamma, j}}{j\in \{1,\dots,n_\gamma\}}
\end{align}
 denote the corresponding set of nodes with $\widehat x_{\gamma, j-1}<\widehat x_{\gamma,j}$ for $j\in\{2,\dots,n_\gamma\}$. Note that $N_\gamma=\sum_{j=1}^{n_\gamma}\#_\gamma\widehat x_{\gamma,j}$.
With $\widehat{x}_{\gamma,0}:=a$,  we define the induced mesh on $[a,b]$,
\begin{align}
\widehat \QQ_\gamma:= \set{[\widehat x_{\gamma, j-1},\widehat x_{\gamma, j}]}{j\in\{1,\dots,n_\gamma\}}.
\end{align}
To use the definition of B-splines as in Section~\ref{sec:splines}, we extend the knot sequence arbitrarily to $(t_{\gamma,i})_{i\in \Z}$ with $t_{\gamma,-p_\gamma}=\dots=t_{\gamma,0}=a$, $t_{\gamma,i}\le t_{\gamma,i+1}$, and $\lim_{i\to \pm\infty}t_{\gamma,i}=\pm \infty$.
For the extended sequence, we also write  $\widehat{\mathcal{K}}_\gamma$.
Moreover, let $\mathcal{W}_\gamma=(w_{\gamma,i})_{i=1-p_\gamma}^{N_\gamma-p_\gamma}$ and $\mathcal{C}_\gamma=(C_i)_{i=1-p_\gamma}^{N_\gamma-p_\gamma}$ be given positive weights and control points in $\R^2$, respectively, 
which satisfy that  $w_{\gamma,1-p_\gamma}=w_{\gamma,N_\gamma-p_\gamma}$ and $C_{\gamma,1-p_\gamma}=C_{\gamma,N_\gamma-p_\gamma}$  in the case of $\Gamma=\partial\Omega$.
We extend  $\mathcal{W}_\gamma$ and $\mathcal{C}_\gamma$ arbitrarily to $(w_{\gamma,i})_{i\in\Z}$ and $(C_{\gamma,i})_{i\in\Z}$  with positive weights $w_{\gamma,i}$ and control points $C_{\gamma,i}$ in $\R^2$, respectively.
With the standard NURBS $ \widehat R_{\gamma,i,p_\gamma}$ defined in Section~\ref{sec:splines}, we suppose that $\gamma$ has the form
\begin{align}
\gamma|_{[a,b)}=\sum_{i\in\Z} C_{\gamma,i} \widehat R_{\gamma,i,p_\gamma}|_{[a,b)}=\sum_{i=1-p_\gamma}^{N_\gamma-p_\gamma} C_{\gamma,i} \widehat R_{\gamma,i,p_\gamma}|_{[a,b)},
\end{align}
where the second equality follows from the locality~\eqref{eq:locality} of the B-splines resp.\ NURBS. 
The locality of NURBS even shows that this definition does not depend on how the knots, the weights, and the control points are precisely extended. 
We note that the assumptions $w_{\gamma,1-p_\gamma}=w_{\gamma,N_\gamma-p_\gamma}$ and $C_{\gamma,1-p_\gamma}=C_{\gamma,N_\gamma-p_\gamma}$ together with~\eqref{eq:interpolatoric} below ensure that $\gamma(a)=\gamma(b-)$ in the case of closed $\Gamma=\partial\Omega$.

\subsection{Rational splines on $\Gamma$}\label{sec:splines on Gamma}
Let $p\in\N_0$ be an arbitrary but fixed polynomial degree.
Let $\widehat\KK_0=(t_{0,i})_{i=1}^{N_0}$ be a sequence of initial knots with 
\begin{align}
a<t_{0,1}\le\dots \le t_{0,N_0-p-1}< t_{0,N_0-p}=\dots= t_{0,N_0}=b
\end{align} 
and multiplicity $\#_0 t_{0, i}\le p+1$ for $i\in\{1,\dots,N_0-p\}$.
Let 
\begin{align}
\widehat{\mathcal{N}}_0:=\set{t_{0, i}}{i\in\{1,\dots N_0\}}=\set{\widehat x_{0, j}}{j\in \{1,\dots,n_0\}}
\end{align}
denote the corresponding set of nodes with $\widehat x_{0, j-1}<\widehat x_{0,j}$ for $j\in\{2,\dots,n_0\}$.
In order to apply standard quadrature rules, we assume that $\widehat \NN_\gamma\subseteq\widehat \NN_0$.
Moreover, let $\mathcal{W}_0=(w_{0,i})_{i=1-p}^{N_0-p}$ be given positive weights. 
We extend the knots and weights as in Section~\ref{sec:boundary parametrization} and define the weight function 
\begin{align}
\widehat W_0|_{[a,b)}:=\sum_{k\in\Z} w_{0,k}\widehat B_{0,k,p}|_{[a,b)}=\sum_{k=1-p}^{N_0-p} w_{0,k}\widehat B_{0,k,p}|_{[a,b)},
\end{align}
where $\widehat B_{0,k,p}$ are the standard B-splines defined in Section~\ref{sec:splines}.
We define $\widehat W_0:[a,b]\to\R$ by continuously extending this function at $b$.

Now, let $\widehat \KK_\coarse=(t_{0,i})_{i=1}^{N_\coarse}$ be a finer knot vector, i.e., $\widehat\KK_0$ is a subsequence of $\widehat \KK_\coarse$ which satisfies the same properties as $\widehat\KK_0$.
Outside the interval $(a,b]$, we extend  $\widehat\KK_\coarse$ as $\widehat\KK_0$.
Let 
\begin{align}
\widehat{\mathcal{N}}_\coarse:=\set{t_{\coarse, i}}{i\in\{1,\dots N_0\}}=\set{\widehat x_{\coarse, j}}{j\in \{1,\dots,n_\coarse\}}
\end{align}
denote the corresponding set of nodes on $[a,b]$ with $\widehat x_{\coarse, j-1}<\widehat x_{\coarse,j}$ for $j\in\{1,\dots,n_\coarse\}$, and let $\NN_\coarse:=\set{\gamma(\widehat{x}_{\coarse, j})}{j\in\{0,\dots,n_\coarse\}}$ denote the corresponding nodes on $\Gamma$.
Note that $\widehat\NN_\coarse$ does not contain the node $\widehat{x}_{\coarse, 0}=a$, which is natural 
for the case $\Gamma=\partial\Omega$,   since then $\gamma(\widehat x_{\coarse,0})=\gamma(\widehat x_{\coarse,n_\coarse})$.
We define the induced mesh on $[a,b]$, 
\begin{align}
\widehat \QQ_\coarse:= \set{[\widehat x_{\coarse, j-1},\widehat x_{\coarse, j}]}{j\in\{1,\dots,n_\coarse\}}
\end{align}
where we set $\widehat{x}_{\coarse,0}:=a$, and the induced mesh on $\Gamma$, $\QQ_\coarse:=\set{\gamma(\widehat Q)}{\widehat Q\in\widehat\QQ_\coarse}$. 
Recall that the non-vanishing B-splines on $[a,b)$ form a basis; see Section~\ref{sec:splines}.
This proves the existence and uniqueness of   weights $\WW_\coarse=(w_{\coarse, i})_{i=1-p}^{N_\coarse-p}$ such that 
\begin{align}\label{eq:w}
\widehat W_0|_{[a,b)}=\sum_{k=1-p}^{N_0-p} w_{0,k}\widehat B_{0,k,p}|_{[a,b)}=\sum_{k=1-p}^{N_\coarse-p} w_{\coarse, k} \widehat B_{\coarse, k,p}|_{[a,b)}.
\end{align}
Choosing these weights, we ensure that the denominator of the rational splines  does not change.
These weights are just convex combinations of the initial weights $\WW_0$ and can be computed via the knot insertion procedure \eqref{eq:insertion}. 
Finally, we extend $\mathcal{W}_\coarse$ arbitrarily to $(w_{\coarse, i})_{i\in\Z}$ with $w_{\coarse, i}>0$ and define the space of (transformed) rational splines on $\Gamma$ as  
\begin{equation}
\SS^p({\mathcal{K}}_\coarse,\mathcal{W}_\coarse):=\Big\{\frac{\widehat S_\coarse}{\widehat W_\coarse}\circ\gamma^{-1}:\,\widehat S_\coarse \in \widehat\SS^p(\widehat{\mathcal{K}}_\coarse)\Big\},
\end{equation}
where $\widehat\SS^p(\widehat{\mathcal{K}}_\coarse)$ denotes the space of all right-continuous piecewise polynomials of degree lower or equal $p$ with breakpoints $\widehat{\mathcal{N}}_\coarse$ on $[a,b)$ 
which are, at each knot $t_{\coarse, i}$, $p-\#_\coarse t_{\coarse, i}$ times continuously differentiable if $p-\#_\coarse t_{\coarse, i}\ge 0$; see also Section~\ref{sec:splines}.
Here, we extend each quotient $\widehat S_\coarse/\widehat W_\coarse$ left-continuously at $b$. 
The locality~\eqref{eq:locality} of  B-splines shows that this definition does not depend on how the knots and the weights are precisely extended. 
With the standard NURBS $\widehat R_{\coarse, i,p}$ from Section~\ref{sec:splines}, a basis of $\SS^p({\mathcal{K}}_\coarse,\mathcal{W}_\coarse)$ is given by 
\begin{align}
\SS^p({\mathcal{K}}_\coarse,\mathcal{W}_\coarse)= {\rm span} \set{{R}_{\coarse, i,p}}{i\in\{1-p,\dots,N_\coarse-p\}}
\quad\text{with}\quad{R}_{\coarse, i,p}:=
\widehat R_{\coarse, i,p}\circ\gamma^{-1},
\end{align}
where the functions $\widehat R_{\coarse, i,p}$  are again left-continuously extended at $b$.



\section{\texttt{IGABEM2D} for weakly-singular integral equation}\label{sec:weak}

\subsection{Sobolev spaces}\label{sec:Sobolev spaces}
The usual Lebesgue and Sobolev spaces on $\Gamma$ are denoted by $L^2(\Gamma)=H^0(\Gamma)$ and  $H^1(\Gamma)$. 
With the weak arclength derivative $\partial_\Gamma$, the $H^1$-norm reads
\begin{align}
\norm{u}{H^1(\Gamma)}^2= \norm{u}{L^2(\Gamma)}^2+ \norm{\partial_\Gamma u}{L^2(\Gamma)}^2\quad\text{for all }u\in H^1(\Gamma),
\end{align}
We stress that $H^1(\Gamma)$ is continuously embedded in $C^0(\Gamma)$. 
Moreover, we equip the space $\H^1(\Gamma):=\set{v\in H^1(\partial\Omega)}{\supp(v)\subseteq\Gamma}$ with the same norm. 
We define the Sobolev space $H^{1/2}(\Gamma)$ as the space of all functions $u\in L^2(\Gamma)$ with finite Sobolev-Slobodeckij norm
\begin{align}
\norm{u}{H^{1/2}(\Gamma)}^2:=\norm{u}{L^2(\Gamma)}^2+\seminorm{u}{H^{1/2}(\Gamma)}^2\quad\text{with}\quad \seminorm{u}{H^{1/2}(\Gamma)}^2:=\int_\Gamma\int_\Gamma \frac{|u(x)-u(y)|^2}{|x-y|^2} \d y \d x.
\end{align}
We will also need  the seminorm $\seminorm{u}{H^{1/2}(\omega)}$ for subsets $\omega\subseteq\Gamma$, which is defined analogously.
Moreover, $\H^{1/2}(\Gamma):=\set{v\in H^{1/2}(\partial\Omega)}{\supp(v)\subseteq\Gamma}$ equipped with the same norm.

Sobolev spaces of negative order are defined by duality
$H^{-1/2}(\Gamma) := \H^{1/2}(\Gamma)^*$ and 
$\H^{-1/2}(\Gamma) := H^{1/2}(\Gamma)^*$, where duality is understood with
respect to the extended $L^2(\Gamma)$-scalar
product $\dual\cdot\cdot_\Gamma$. 
Note that  $\H^{\pm \sigma}(\partial\Omega) = H^{\pm \sigma}(\partial\Omega)$ in case of $\Gamma=\partial\Omega$
even with equal norms for $\sigma\in\{0,1/2,1\}$.

All details and equivalent definitions of the Sobolev spaces 
are, for instance, found in the monographs~\cite{mclean00,hw08,steinbach08,ss11}.

\subsection{Weakly-singular integral equation}
\label{sec:weaksing}
For $\sigma\in\{0,1/2,1\}$, the single-layer operator 
$\V:\H^{\sigma-1}(\Gamma)\to H^{\sigma}(\Gamma)$ and the double-layer operator $\K:\H^{\sigma}(\Gamma)\to H^{\sigma}(\Gamma)$ are well-defined, linear, and continuous. 
As $H^1(\Gamma)\subset C^0(\Gamma)$, the case $\sigma=1$ yields continuous functions, which is essential for the implementation described in Section~\ref{sec:implementation}.

For $\sigma=1/2$, 
$\V:\H^{-1/2}(\Gamma)\to H^{1/2}(\Gamma)$ is symmetric and 
elliptic under the assumption that $\diam(\Omega)<1$.
In particular,
\begin{align}
\edualV{\phi}{\psi}:=\dual{\V \phi}{\psi}_\Gamma
\end{align}
 defines an equivalent scalar product on  
$\H^{-1/2}(\Gamma)$ with corresponding norm $\enormV{\cdot}$.
With this notation, the strong form~\eqref{eq:weak strong} with data $f=u$ for some $u\in H^{1/2}(\Gamma)$ or $f=(1/2+\K)u$ for some $u\in\H^{1/2}(\Gamma)$ is equivalently stated as
\begin{align}
\label{eq:weak weak}
\edualV{\phi}{\psi} = \dual{f}{\psi}_\Gamma
 \quad\text{for all }\psi\in \H^{-1/2}(\Gamma).
\end{align}
The Lax--Milgram lemma applies and hence  \eqref{eq:weak strong} admits a unique solution $\phi\in \H^{-1/2}(\Gamma)$.
If $f=u$, the approach is called indirect, otherwise if $f=(1/2+\K)u$, it is called direct.
More details and proofs are found, e.g.,  in  \cite{mclean00,hw08,steinbach08,ss11}.

\subsection{Galerkin IGABEM}
\label{sec:spaces} 
Let $p\in\N_0$ be a fixed polynomial degree.
Moreover, let $\widehat\KK_\coarse$ and $\WW_\coarse$ be knots and weights as in Section~\ref{sec:splines on Gamma}.
We introduce the  ansatz space
\begin{equation}\label{eq:weaksing X0}
\XX_\coarse:=\SS^p({\mathcal{K}}_\coarse, \mathcal{W}_\coarse)\subset L^2(\Gamma). 
\end{equation}
Recall that
\begin{align}\label{eq:weaksing basis}
\XX_\coarse={\rm span}{\set{R_{\coarse, i,p}}{i=1-p,\dots,N_\coarse-p}},
\end{align}
where the set even forms a basis.
We define the Galerkin approximation $\Phi_\coarse\in\XX_\coarse$ of $\phi$  by
\begin{align}\label{eq:weaksing Galerkin}
 \edualV{\Phi_\coarse}{\Psi_\coarse} = \dual{f}{\Psi_\coarse}_\Gamma  
\,\text{ for all }\Psi_\coarse\in\XX_\coarse. 
\end{align}

Note that \eqref{eq:weaksing Galerkin} is equivalent to solving the finite-dimensional linear system 
\begin{align}\label{eq:Galerkin_weak}
\big(\edualV{R_{\coarse,i',p}}{R_{\coarse,i,p}}\big)_{i,i'=1-p}^{N_\coarse-p} \cdot (c_{\coarse,i'})_{i'=1-p}^{N_\coarse-p}= \big(\dual{f}{R_{\coarse,i,p}}_\Gamma\big)_{i=1-p}^{N_\coarse-p},
\end{align}
where $\Phi_\coarse = \sum_{i'=1-p}^{N_\coarse-p} c_{\coarse,i'} R_{\coarse,i',p}$.
The computation of the matrix and the right-hand side vector in~\eqref{eq:Galerkin_weak} is realized in {\tt VMatrix.c} and {\tt RHSVectorWeak.c} and can be called in {\sc Matlab} via {\tt buildVMatrix} and {\tt buildRHSVectorWeak}; see Appendix~\ref{sec:overview}.  
For details on the implementation, we refer to Section~\ref{sec:matrices} and \ref{sec:right vectors}.

\subsection{\textsl{A~posteriori} error estimation}
\label{sec:a_posteriori_weak}

Let $\widehat \KK_\coarse$ be a knot vector as in Section~\ref{sec:spaces}  with corresponding ansatz space $\XX_\coarse$. 
We define the mesh-size function $h_\coarse\in L^\infty(\Gamma)$ by $h_\coarse|_Q:=|Q|$ for all $Q\in\QQ_\coarse$. 
Moreover, we abbreviate the node patch
\begin{align}
\omega_\coarse(x):=\bigcup \set{Q\in\QQ_\coarse}{x\in Q}\quad\text{for all }x\in\NN_\coarse.
\end{align}
We consider the following three different node-based estimators:
the $(h-h/2)$-estimator
\begin{align}\label{eq:hhV}
\eta_{\V,\rm hh2,\coarse}^2:=\sum_{x\in\NN_\coarse} \eta_{\V,\rm hh2,\coarse}(x)^2\quad\text{with}\quad  \eta_{\V,\rm hh2,\coarse}(x)^2:= \norm{h_\coarse^{1/2}(\Phi_\fine-\Phi_\coarse)}{L^2(\omega_\coarse(x))}^2,
\end{align}
where  $\widehat\KK_\fine$ are the uniformly refined knots obtained from $\widehat \KK_\coarse$ by adding the midpoint of each element $\widehat Q\in\widehat\QQ_\coarse$ with multiplicity one (see Algorithm~\ref{alg:refinement});
the Faermann estimator
\begin{align}\label{eq:fae}
\eta_{\V,\rm fae,\coarse}^2:=\sum_{x\in\NN_\coarse} \eta_{\V,\rm fae,\coarse}(x)^2\quad\text{with}\quad  \eta_{\V,\rm fae,\coarse}(x)^2:= \seminorm{f-\V\Phi_\coarse}{H^{1/2}(\omega_\coarse(x))}^2;
\end{align}
and the weighted-residual estimator
\begin{align}\label{eq:resV}
\eta_{\V,\rm res,\coarse}^2:=\sum_{x\in\NN_\coarse} \eta_{\V,\rm res,\coarse}(x)^2\text{ with }  \eta_{\V,\rm res,\coarse}(x)^2:= \norm{h_\coarse^{1/2}\partial_\Gamma(f-\V\Phi_\coarse)}{L^{2}(\omega_\coarse(x))}^2,
\end{align}
which requires the additional regularity $f\in H^1(\Gamma)$ to ensure that $f-\mathfrak{V}\Phi_\coarse\in H^1(\Gamma)$ and hence \eqref{eq:resV} is well-defined.
If $f=(\K+1/2)u$ stems from a Dirichlet problem as in Section~\ref{sec:model}, Section~\ref{sec:weaksing} shows that this is satisfied for $u\in \H^1(\Gamma)$.

The computation of the estimators is implemented in \texttt{HHEstWeak.c}, \texttt{FaerEstWeak.c}, and \texttt{ResEstWeak.c} 
and can be called in {\sc Matlab} via \texttt{buildHHEstWeak}, \texttt{buildFaerEstWeak}, and \texttt{buildResEstWeak}; see Appendix~\ref{sec:overview}.
The residual estimators require the evaluation of the boundary integral operators $\mathfrak{V}$ and $\mathfrak{K}$ applied to some function, which is implemented in \texttt{ResidualWeak.c}. 
The evaluations can be used in {\sc Matlab} via \texttt{evalV} and \texttt{evalRHSWeak}; see Appendix~\ref{sec:overview}.
The implementation of the estimators and the evaluations are discussed in Section~\ref{sec:implementation}.

\subsection{Adaptive IGABEM algorithm}

To refine a given ansatz space $\XX_\coarse$, we compute one of the corresponding error estimators $\eta_\coarse$ and determine a set of marked nodes with large error indicator $\eta_\coarse(x)$.
By default, we then apply the following refinement algorithm, which uses $h$-refinement controlling the maximal mesh-ratio in the parameter domain  
\begin{align}\label{eq:meshratio}
\widehat{\kappa}_\coarse&:=\max\Big\{\frac{|\gamma^{-1}(Q)|}{|\gamma^{-1}(Q')|}:{Q},{Q}'\in\QQ_\coarse \text{ with }  Q\cap Q'\neq \emptyset\Big\}
\end{align}
as well as multiplicity increase.

 \begin{algorithm}\label{alg:refinement}
\textbf{Input:} Polynomial degree $p\in\N_0$, initial maximal mesh-ratio $\widehat\kappa_0$, knot vector  $\widehat\KK_\coarse$ with $\widehat\kappa_\coarse\le 2\widehat\kappa_0$, marked nodes $\MM_\coarse\subseteq\NN_\coarse$.
\begin{enumerate}[\rm(i)]
\item Define the set of marked elements $\MM_{\coarse}':=\emptyset$.
\item If both nodes of an element $Q\in \QQ_\coarse$ belong to $\MM_\coarse$, mark  $Q$ by adding it to $\MM_{\coarse}'$.
\item For all other nodes in $\MM_\coarse$, increase the multiplicity if it is less or equal to $p$. 
Otherwise mark the elements which contain one of these nodes,  by adding them to $\MM_{\coarse}'$.
\item
Bisect all $Q\in \MM_\coarse'$ in the parameter domain by inserting the midpoint of  $\gamma^{-1}(Q)$ with multiplicity one to the current knot vector.  
Use a minimal number of further bisections to guarantee that the new knot vector $\widehat\KK_\fine$ satisfies that $\widehat\kappa_\fine\le 2\widehat\kappa_0$.
\end{enumerate}
\textbf{Output:} Refined knot vector $\widehat\KK_\fine$.
\end{algorithm}

The marked elements $\MM_\ell'$ and the nodes whose multiplicity should be increased are determined in {\tt markNodesElements.m}.
The multiplicity increase and the computation of the new weights are realized in {\tt increaseMult.m}. 
An optimal 1D bisection algorithm as in {\rm (iv)} is discussed and analyzed in \cite{affkp13}. 
Together with the computation of the new weights, it is realized in  {\tt refineBoundaryMesh.m}.  
We stress that we have also implemented two further relevant strategies that rely on $h$-refinement only:
Replace {\rm (iii)} by adding all elements $Q\in\QQ_\coarse$ containing one of the other nodes in $\MM_\coarse$ to $\MM_{\coarse}'$;
or replace {\rm (iii)} by adding all elements $Q\in\QQ_\coarse$ containing one of the other nodes in $\MM_\coarse$ to $\MM_{\coarse}'$
and  insert the midpoints in {\rm (iv)} with multiplicity $p+1$.  
The first strategy leads to refined splines of full regularity, whereas the second one leads to lowest regularity.

We fix the considered error estimator $\eta_\coarse\in \{\eta_{\V,\rm hh2,\coarse},$ $\eta_{\V,\rm fae,\coarse},$ $\eta_{\V,\rm res,\coarse}\}$.
The corresponding error indicators $\eta_\coarse(x)$ are defined accordingly.

\begin{algorithm}\label{the algorithm}
\textbf{Input:} Polynomial degree $p\in \N_0$, initial knot vector $\widehat\KK_0$, initial weights~$\mathcal{W}_0$, marking parameter $0< \theta\le1$.\\
\textbf{Adaptive loop:} For each $\ell=0,1,2,\dots$ iterate the following steps {\rm(i)--(iv)}:
\begin{itemize}
\item[\rm(i)] Compute  approximation $\Phi_\ell\in\XX_\ell$ by solving~\eqref{eq:weaksing Galerkin}.
\item[\rm(ii)] Compute refinement indicators $\eta_\ell({x})$
for all nodes ${x}\in\NN_\ell$.
\item[\rm(iii)] Determine a minimal set of nodes $\MM_\ell\subseteq\NN_\ell$ such that
\begin{align}\label{eq:Doerfler}
 \theta\,\eta_\ell^2 \le \sum_{{x}\in\MM_\ell}\eta_\ell({x})^2.
\end{align}
\item[\rm(iv)] Generate refined knot vector $\widehat\KK_{\ell+1}$ via Algorithm~\ref{alg:refinement}.
\end{itemize}
\textbf{Output:} Approximate solutions $\Phi_\ell$ and estimators $\eta_\ell$ for all $\ell \in \N_0$.
\end{algorithm}

The adaptive algorithm is realized in the main function {\tt IGABEMWeak.m}. 
The considered problem as well as the used parameters can be changed there by the user.
In Section~\ref{sec:implementation}, we discuss how the arising (singular) boundary integrals are computed.
We also mention that Algorithm~\ref{the algorithm}~(iii) is realized in {\tt markNodesElements.m}, which also determines the corresponding set of marked elements $\MM_\ell'$ and the nodes whose multiplicity should be increased from Algorithm~\ref{alg:refinement}~(ii)--(iii).



\subsection{Numerical experiments}
\label{section:numerics_weak}

\begin{figure}
\centering
	\includegraphics[width=.475\textwidth,clip=true]{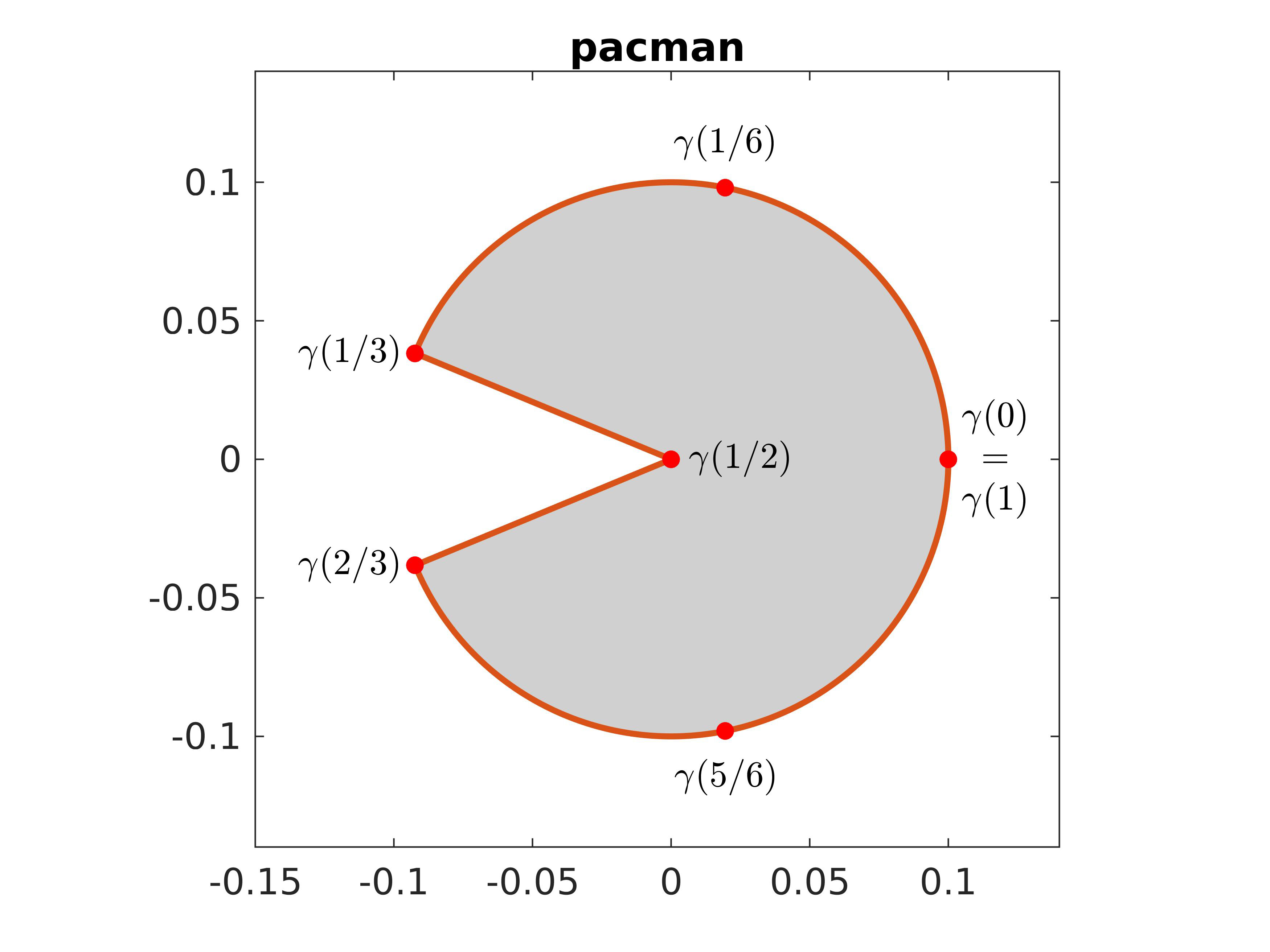}
	\includegraphics[width=.475\textwidth,clip=true]{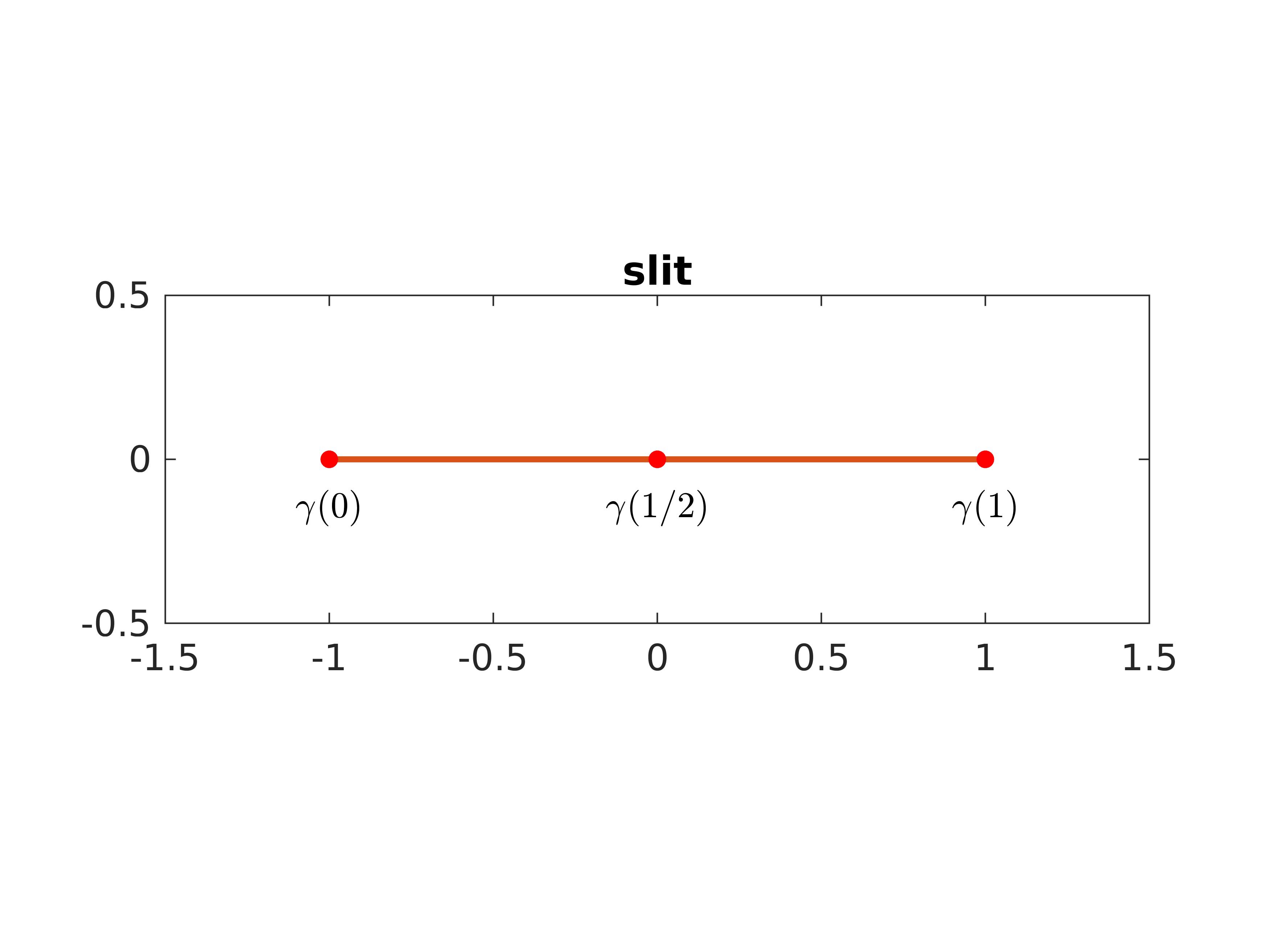}
	\caption{Geometries and initial nodes for examples from Section~\ref{section:numerics_weak}.}
	\label{fig:geometries_weak}
\end{figure}

In this section, we empirically investigate the performance of Algorithm~\ref{the algorithm} for different geometries, ansatz spaces, and error estimators. 
Figure~\ref{fig:geometries_weak} shows the different geometries, namely a pacman geometry and a slit. 
The boundary of the pacman geometry can be parametrized via rational splines of degree~$2$, while the slit can be parametrized by splines of degree~$1$.
As initial ansatz space, we either consider the same (rational) splines, i.e., $\widehat\KK_0=\widehat\KK_\gamma$ and $\WW_0=\WW_\gamma$, splines of degree $p$ and smoothness $C^{p-1}$, 
i.e., 
\begin{align*}
\widehat\KK_0 =(\widehat x_{\gamma,1}, \widehat x_{\gamma,2},\dots,\widehat x_{\gamma,n_\gamma-1},\underbrace{\widehat x_{\gamma,n_\gamma},\dots,\widehat x_{\gamma,n_\gamma}}_{\# = p+1})
\end{align*} 
and $\WW_0=(1,\dots,1)$, 
or piecewise polynomials of degree $p$, i.e., 
\begin{align*}
\widehat\KK_0 =(\underbrace{\widehat x_{\gamma,1},\dots, \widehat x_{\gamma,1}}_{\#=p+1},\underbrace{\widehat x_{\gamma,2},\dots, \widehat x_{\gamma,2}}_{\#=p+1},\dots,\underbrace{\widehat x_{\gamma,n_\gamma},\dots,\widehat x_{\gamma,n_\gamma}}_{\#=p+1})
\end{align*} 
and $\WW_0=(1,\dots,1)$.
In the latter case, we always consider $h$-refinement with new knots having multiplicity $p+1$ as explained after Algorithm~\ref{alg:refinement}.

\subsubsection{Singular problem on pacman geometry}
\label{section:weak_singular_pacman}
We prescribe an exact solution $P$  of the Laplace problem in polar coordinates $(x_1,x_2)=r(\cos\beta,\sin\beta)$ with $\beta\in(-\frac{\pi}{2\tau},\frac{\pi}{2\tau})$ and $\tau:=4/7$ as
\begin{align*}
P(x,y):= r^{\tau}\cos(\tau\beta)
\end{align*}
on the pacman domain, cf. Figure~\ref{fig:geometries_weak}. The solution $\phi$ of the weakly-singular integral equation \eqref{eq:weak strong} with $f=(1/2+\mathfrak{K})(P|_\Gamma)$ is just the normal derivative  of $P$, which reads 
\begin{align*}
\phi(x,y)=\left(
\begin{array}{c}
\cos(\beta)\cos(\tau\beta)+\sin(\beta)\sin(\tau\beta) \\
\sin(\beta)\cos(\tau\beta)-\cos(\beta)\sin(\tau\beta)
\end{array}
\right)\nu(x,y)\,\tau\, r^{\tau-1}
\end{align*}
and has a generic singularity at the origin.

\begin{figure}
\centering
	\small{\quad\quad Pacman geometry}\\
	\vspace{2mm}
	\tiny{\hspace*{12mm} Without multiplicity increase}\\
	\vspace{2mm}
	\includegraphics[width=.475\textwidth,clip=true]{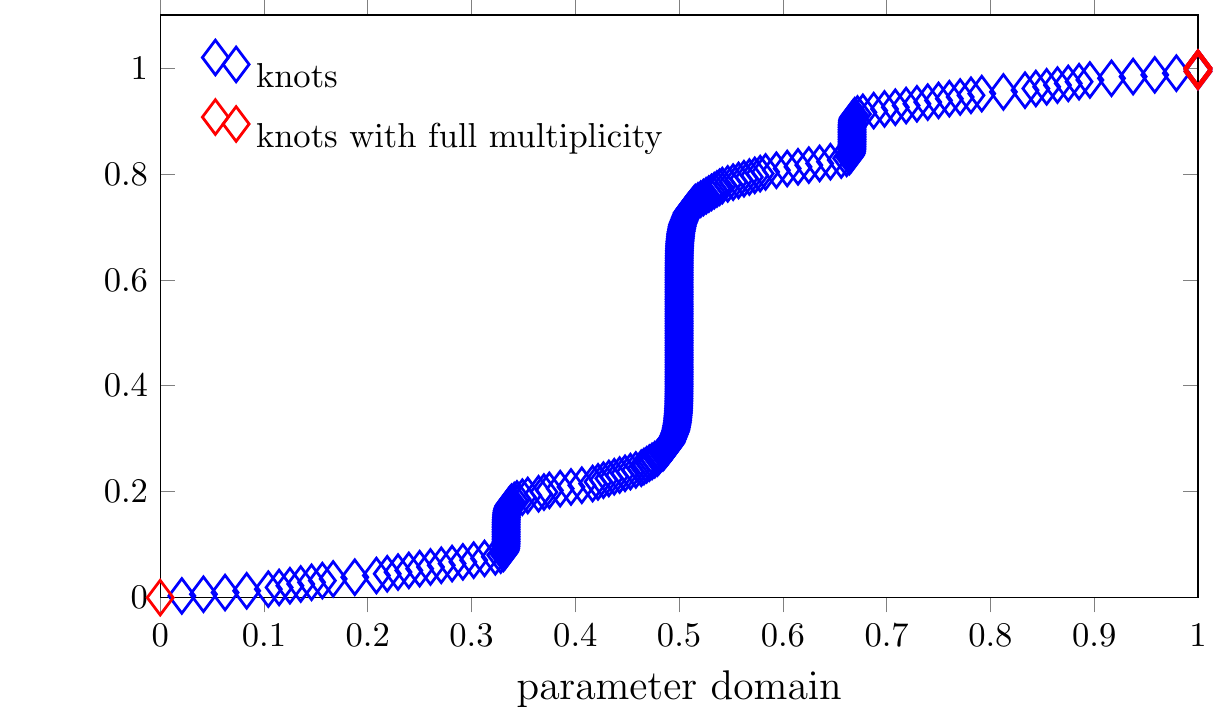}\quad
	\includegraphics[width=.475\textwidth,clip=true]{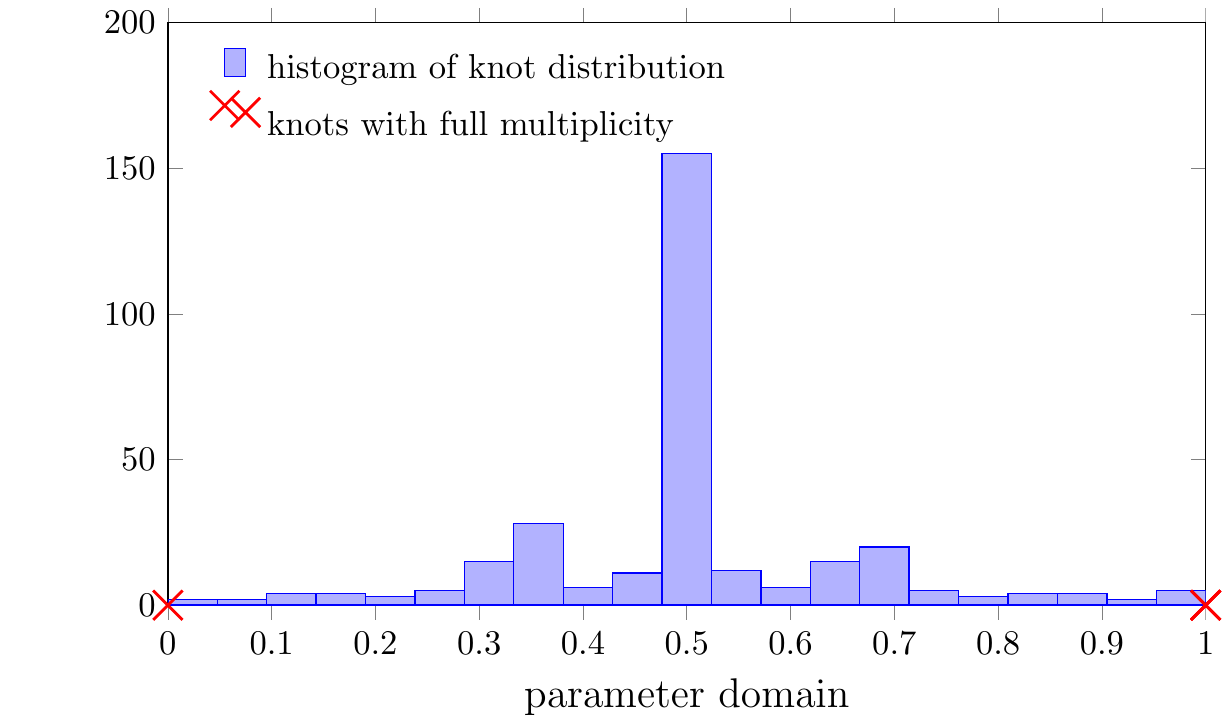}
	\tiny{\hspace*{12mm} With multiplicity increase}\\
	\vspace{2mm}
	\includegraphics[width=.475\textwidth,clip=true]{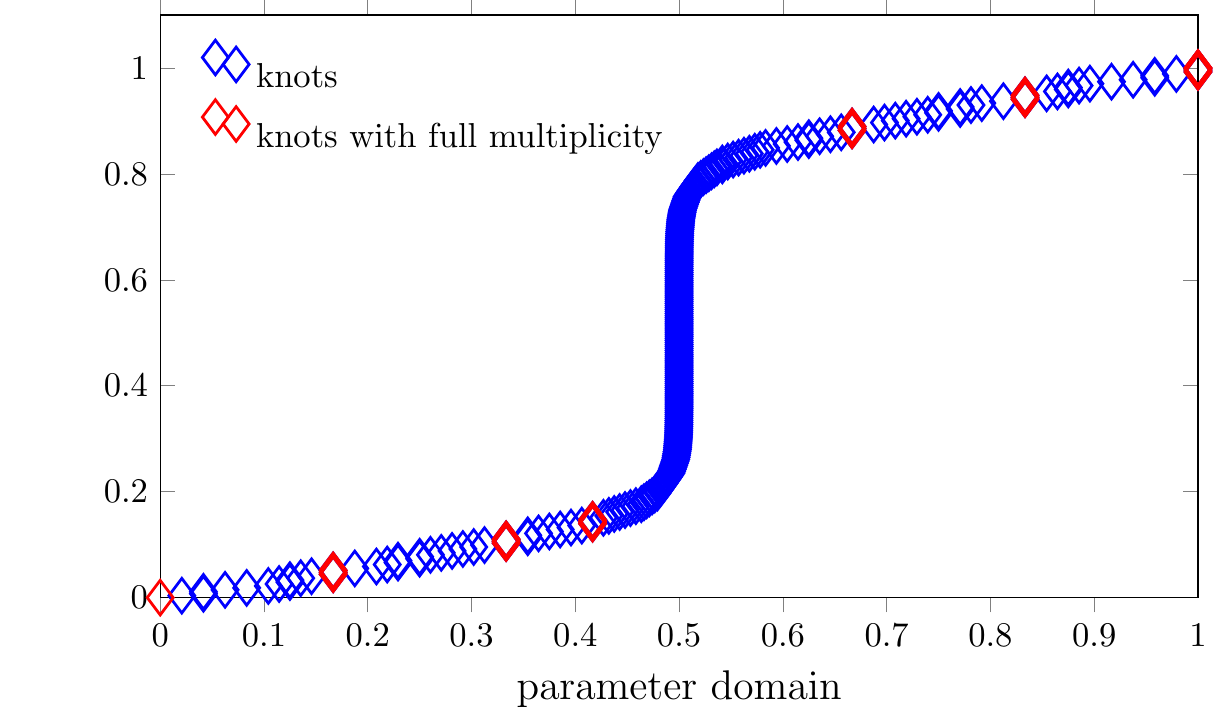}\quad
	\includegraphics[width=.475\textwidth,clip=true]{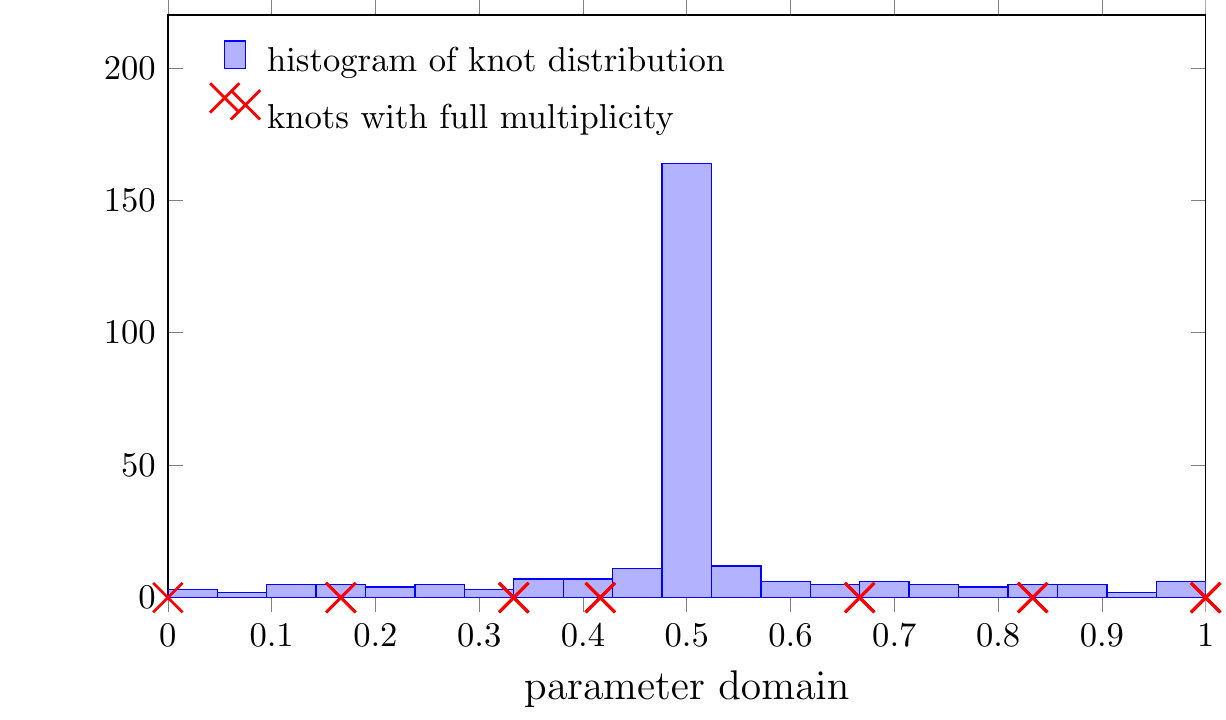}
	\caption{Example from Section~\ref{section:weak_singular_pacman}: Distribution and histogram of the knots of the first step of Algorithm~\ref{the algorithm} where the error estimator is below $10^{-6}$ without multiplicity increase (top) as well as with mulitplicity increase (bottom) for splines of order $p=2$ and the weighted-residual estimator with $\theta=0.75$.}
	\label{fig:weak_pacman_dist}
\end{figure}

Figure~\ref{fig:weak_pacman_dist} shows the knot distribution (i.e., the relative number of knots lower or equal than the parameter points on the $x$-axis) as well as a histogram of the knots in the parameter domain $[a,b]$ of the first step of Algorithm~\ref{the algorithm} where the residual error estimator is below $10^{-6}$. We compare pure $h$-refinement (top) to the proposed refinement strategy including multiplicity increase (bottom), i.e., Algorithm~\ref{alg:refinement}, where the knots with full multiplicity $p+1$ are plotted in red. Splines of order $p=2$ are used as ansatz spaces, and the weighted-residual estimator with $\theta=0.75$ is used to steer the adaptive refinement. It can be seen that Algorithm~\ref{the algorithm} heavily refines the mesh towards the reentrant corner $(0,0)\in\Gamma$, which corresponds to $\gamma(1/2)$ and where $\phi$ has a singularity. The solution $\phi$ additionally has jumps at $\gamma(1/3)$ and $\gamma(2/3)$. Since we use splines of  polynomial degree $p=2$ and there are no knots with multiplicity $p+1$ in $(0,1)$ for Algorithm~\ref{the algorithm} without multiplicity increase (top), we know that the resulting approximation is at least continuous at each knot, cf.~Section~\ref{sec:splines}, and that it can thus not resolve these jumps appropriately. As a result, the algorithm uses $h$-refinement at these jump points to reduce the error. Since Algorithm~\ref{the algorithm} with multiplicity increase (bottom) allows for knots with full multiplicity $p+1$, i.e., it allows for jumps in the resulting approximation, the jumps at $\gamma(1/3)$ and $\gamma(2/3)$ are simply resolved by knots with full multiplicity.

\begin{figure}
\centering
	\small{\hspace*{10mm} Pacman geometry}\\
	\vspace{2mm}
	\begin{minipage}{.475\textwidth}
	\centering
	\tiny{\hspace{15mm} Faermann estimator}
	\vspace{2mm}
	\end{minipage}\quad
	\begin{minipage}{.475\textwidth}
	\centering
	\tiny{\hspace{15mm} Weighted-residual estimator}
	\vspace{2mm}
	\end{minipage}
	\tiny{\hspace*{12mm} Error}\\
	\vspace{2mm}
	\includegraphics[width=.475\textwidth,clip=true]{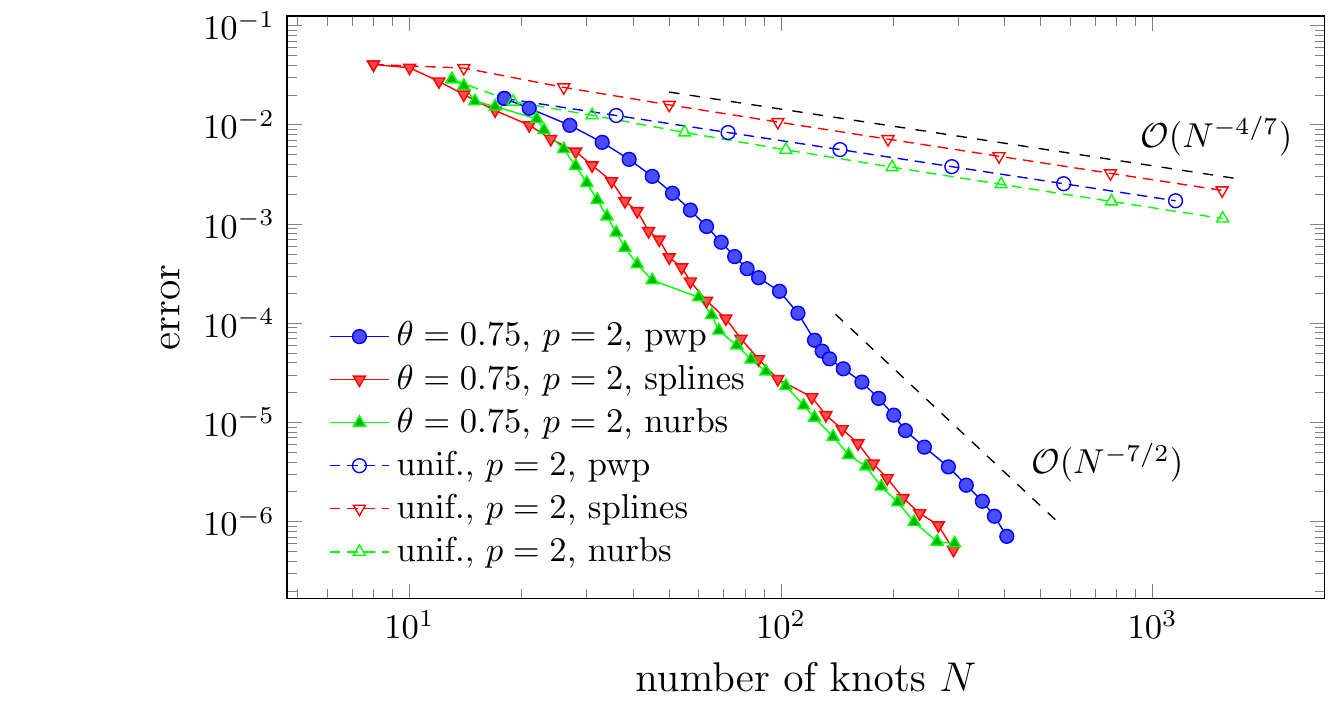}\quad
	\includegraphics[width=.475\textwidth,clip=true]{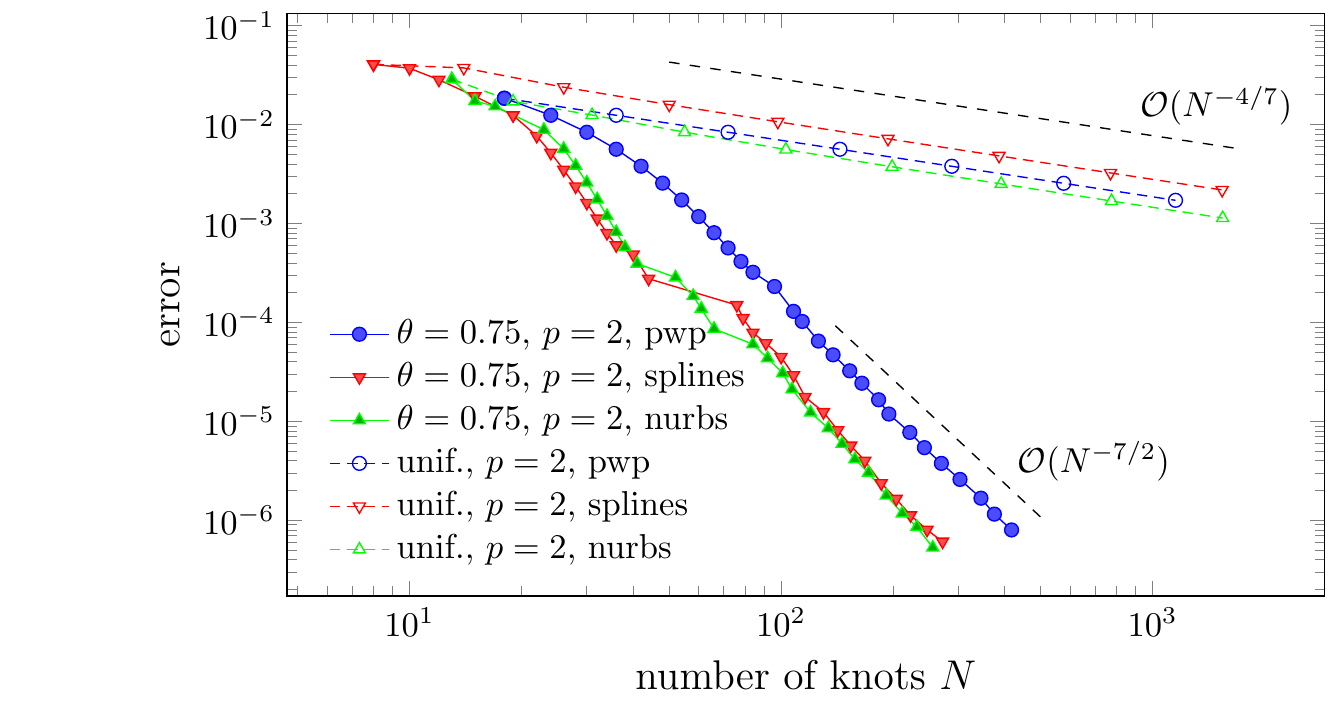}	
	\vspace{2mm}
	\tiny{\hspace*{12mm} Piecewise polynomials}\\
	\vspace{2mm}
	\includegraphics[width=.475\textwidth,clip=true]{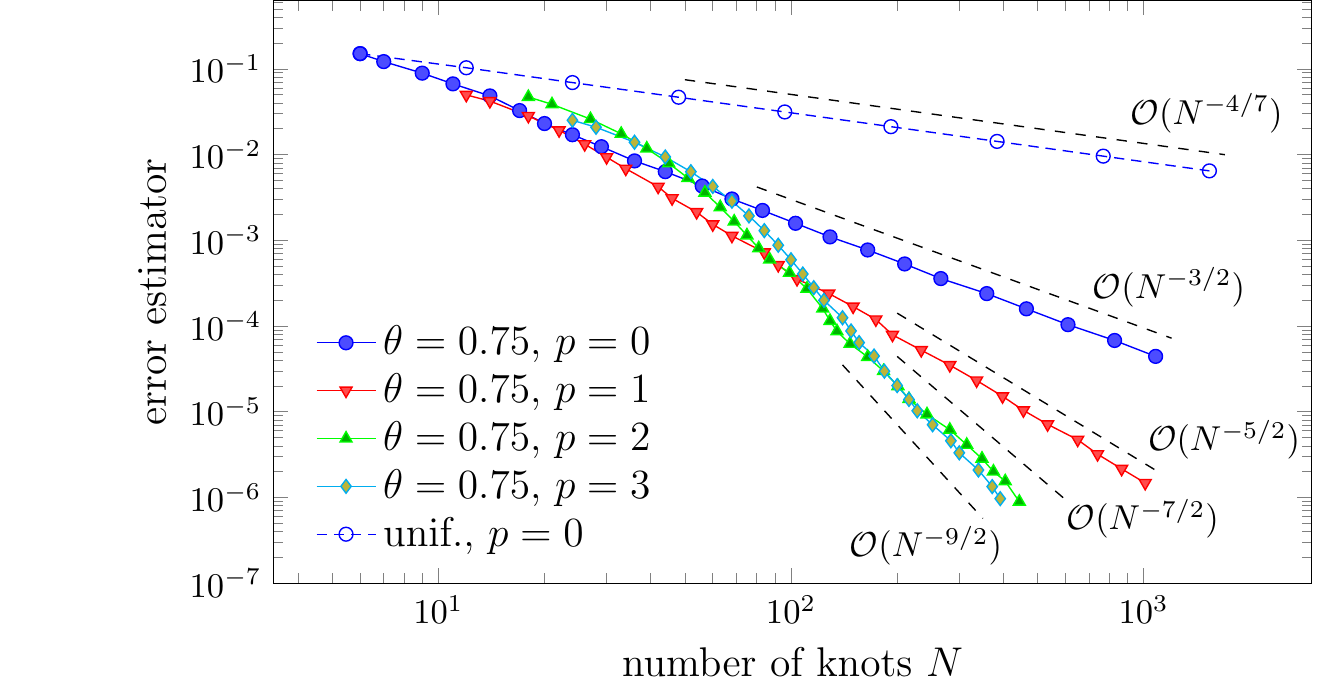}\quad
	\includegraphics[width=.475\textwidth,clip=true]{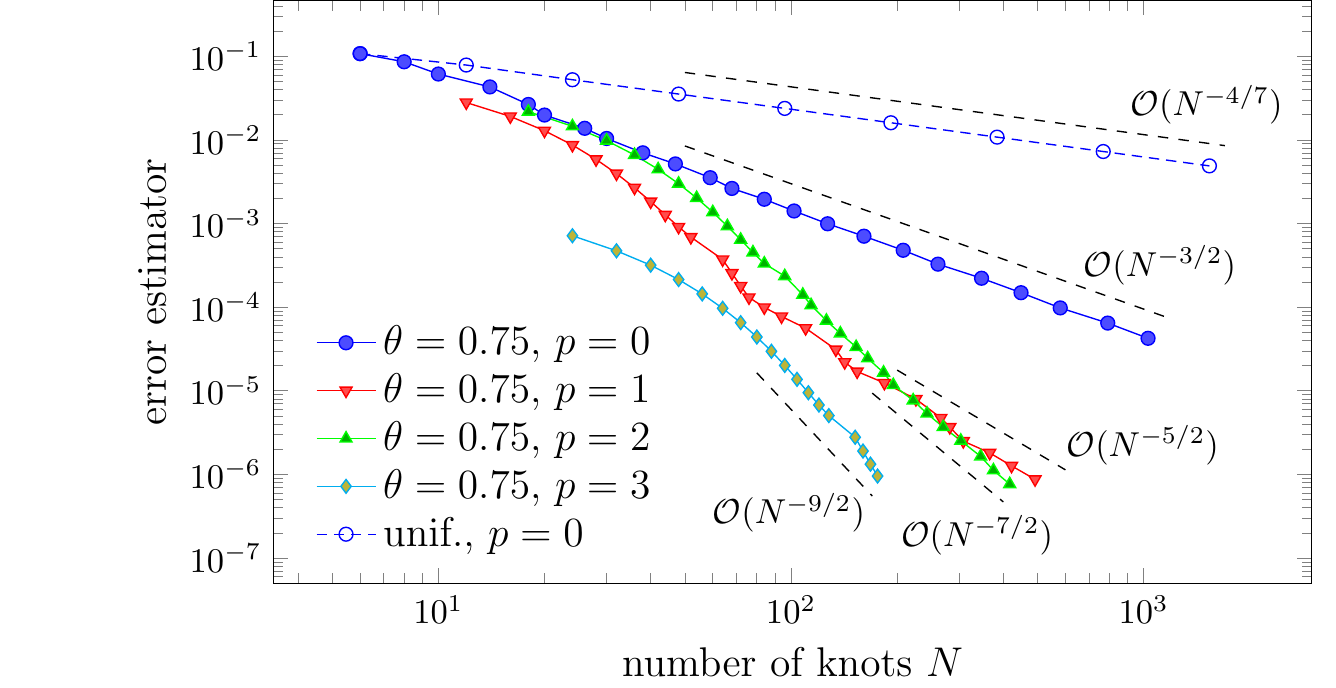}	
	\vspace{2mm}	
	\tiny{\hspace*{12mm} Splines}\\
	\vspace{2mm}	
	\includegraphics[width=.475\textwidth,clip=true]{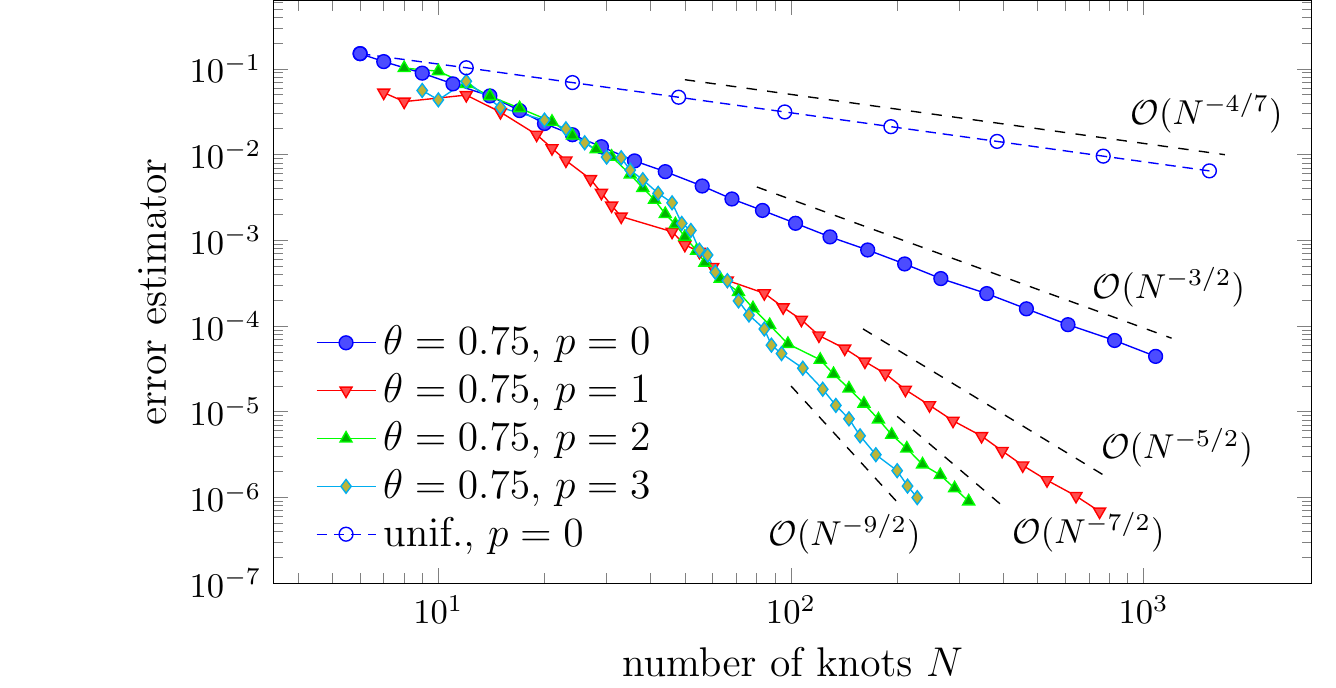}\quad
	\includegraphics[width=.475\textwidth,clip=true]{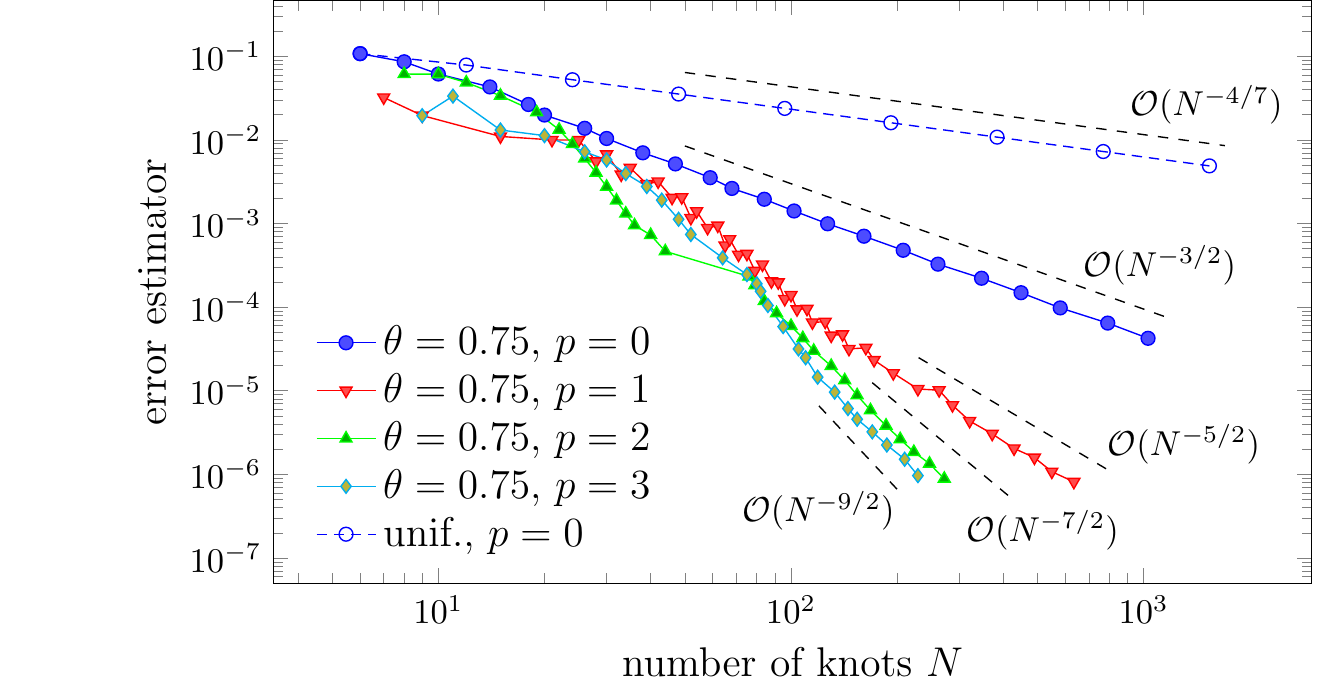}

	\caption{Example from Section~\ref{section:weak_singular_pacman}: (Approximate) error $\norm{\phi-\Phi_\ell}{\mathfrak V}$ and error estimator $\eta_\ell$ with respect to the number of knots $N$ for $\theta\in\{0.75,1\}$, $p\in\{0,1,2,3\}$, different ansatz spaces (piecewise polynomials, splines), and different estimators (Faermann, weighted-residual). Left: Faermann estimator. Right: Weighted-residual estimator. Top: Comparison of $p\in\{0,1,2,3\}$ for piecewise polynomials. Middle: Comparison of $p\in\{0,1,2,3\}$ for splines. Bottom: Comparison of $p\in\{1,2,3\}$ for different ansatz spaces.}
	\label{fig:weak_pacman_comparison}
\end{figure}

In Figure~\ref{fig:weak_pacman_comparison}, we compare uniform refinement with $\theta=1$ and adaptive refinement with $\theta=0.75$, different ansatz spaces, and different estimators for Algorithm~\ref{the algorithm}. 
First (top), for $p=2$ fixed, we compare (an approximation of the) error for the different ansatz spaces and the algorithm steered either by the Faermann estimator (top, left) or the weighted-residual estimator (top, right). 
Instead of the exact error $\norm{\phi-\Phi_\ell}{\mathfrak V}$, we compute $\norm{\Phi_{+}-\Phi_\ell}{\mathfrak V}$,  where  $\widehat\KK_\fine$ are the uniformly refined knots, which result from $\widehat \KK_\ell$ by adding the midpoint of each element $\widehat Q\in\widehat\QQ_\coarse$ with multiplicity one (see also Algorithm~\ref{alg:refinement}).
We get similar plots, where uniform mesh refinement leads to the suboptimal rate of convergence $\OO(N^{-4/7})$ because the solution lacks regularity. However, for all adaptive cases, Algorithm~\ref{the algorithm} regains the optimal rate of convergence $\OO(N^{-3/2-p})$, where splines and NURBS exhibit a significantly better multiplicative constant than piecewise polynomials.
Next (middle), we consider the estimators for piecewise polynomials. For both the Faermann estimator (middle, left) and the weighted-residual estimator (middle, right), we compare different orders $p\in\{0,1,2,3\}$ for the ansatz spaces. Again, uniform mesh refinement leads to the suboptimal rate $\OO(N^{-4/7})$ (only displayed for $p=0$), whereas the adaptive strategy regains the optimal order of convergence $\OO(N^{-3/2-p})$. 
Lastly, we get similar results for splines (bottom).

\subsubsection{Crack problem on slit geometry}
\label{section:weak_crack_slit}
We consider a crack problem on the slit $\Gamma=[-1,1]\times\{0\}$, cf. Figure~\ref{fig:geometries_weak}. Let $f$ in \eqref{eq:weak weak} be defined as $f(x,0):= -x/2$. Then, the exact solution of \eqref{eq:weak weak} reads
\begin{align}
\phi(x,0)=\frac{-x}{\sqrt{1-x^2}}
\end{align}
and has singularities at the tips $x=\pm 1$.

\begin{figure}
\centering
	\small{\quad\quad Slit geometry}\\
	\vspace{2mm}
	\tiny{\hspace*{12mm} Without multiplicity increase}\\
	\vspace{2mm}
	\includegraphics[width=.475\textwidth,clip=true]{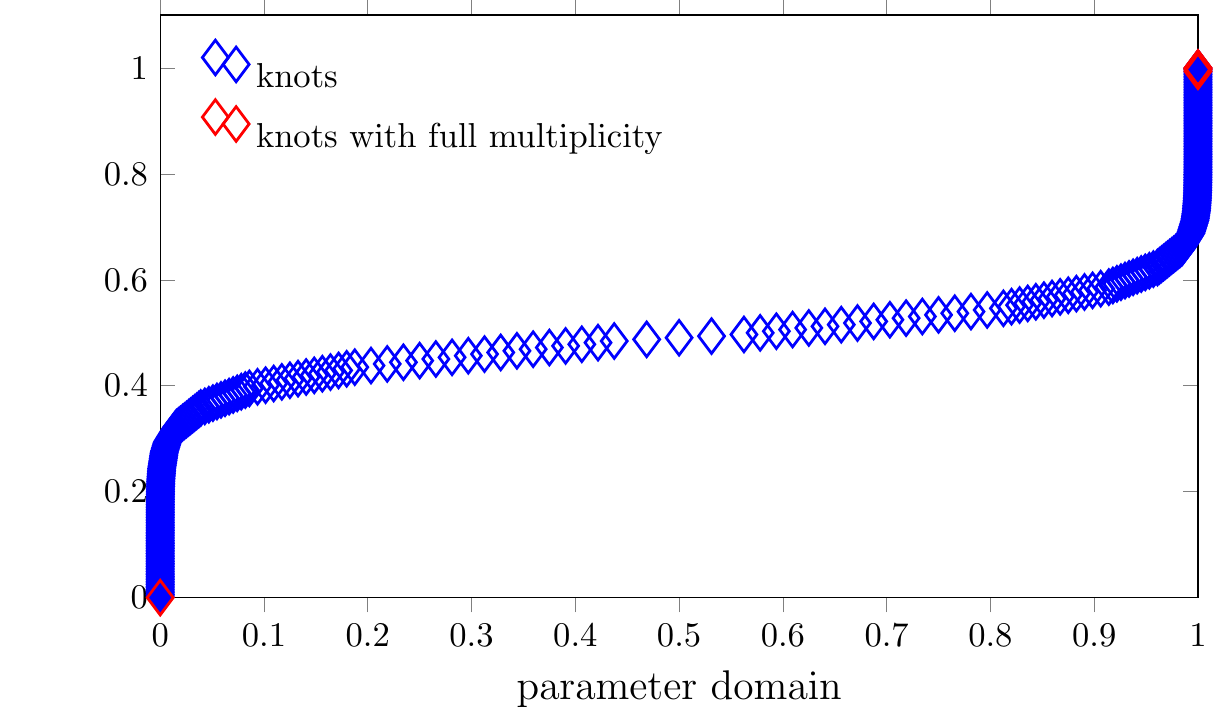}\quad
	\includegraphics[width=.475\textwidth,clip=true]{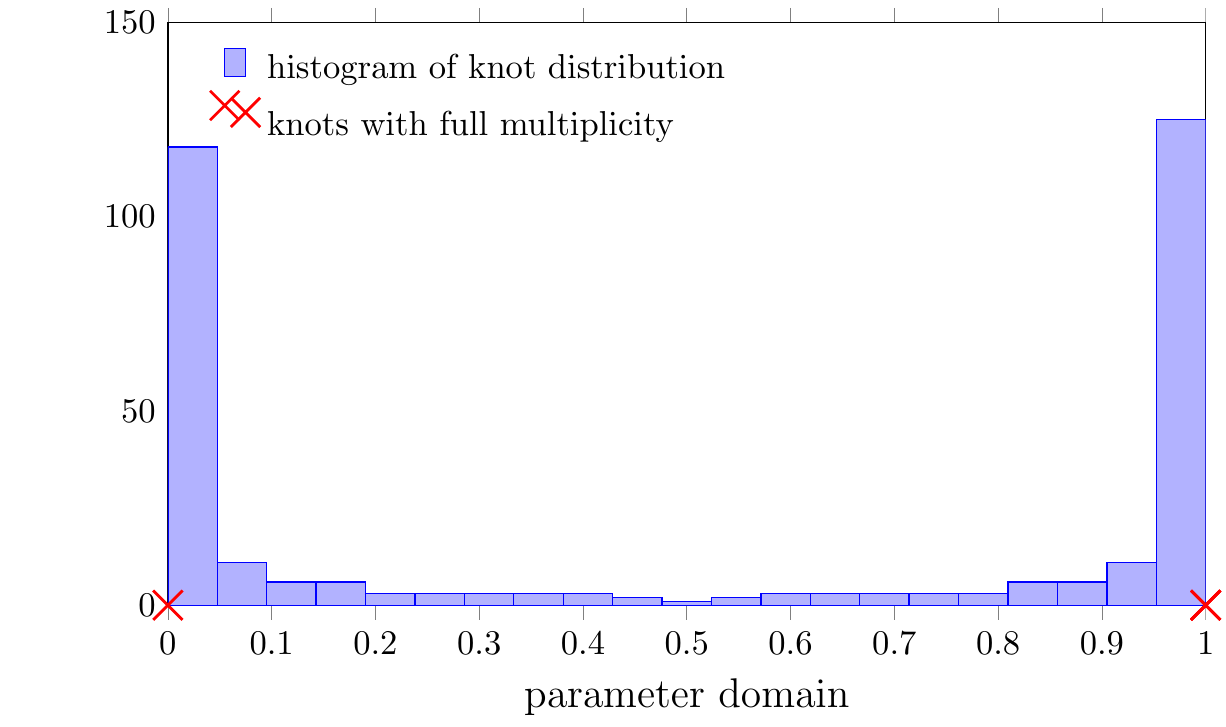}
	\tiny{\hspace*{12mm} With multiplicity increase}\\
	\vspace{2mm}
	\includegraphics[width=.475\textwidth,clip=true]{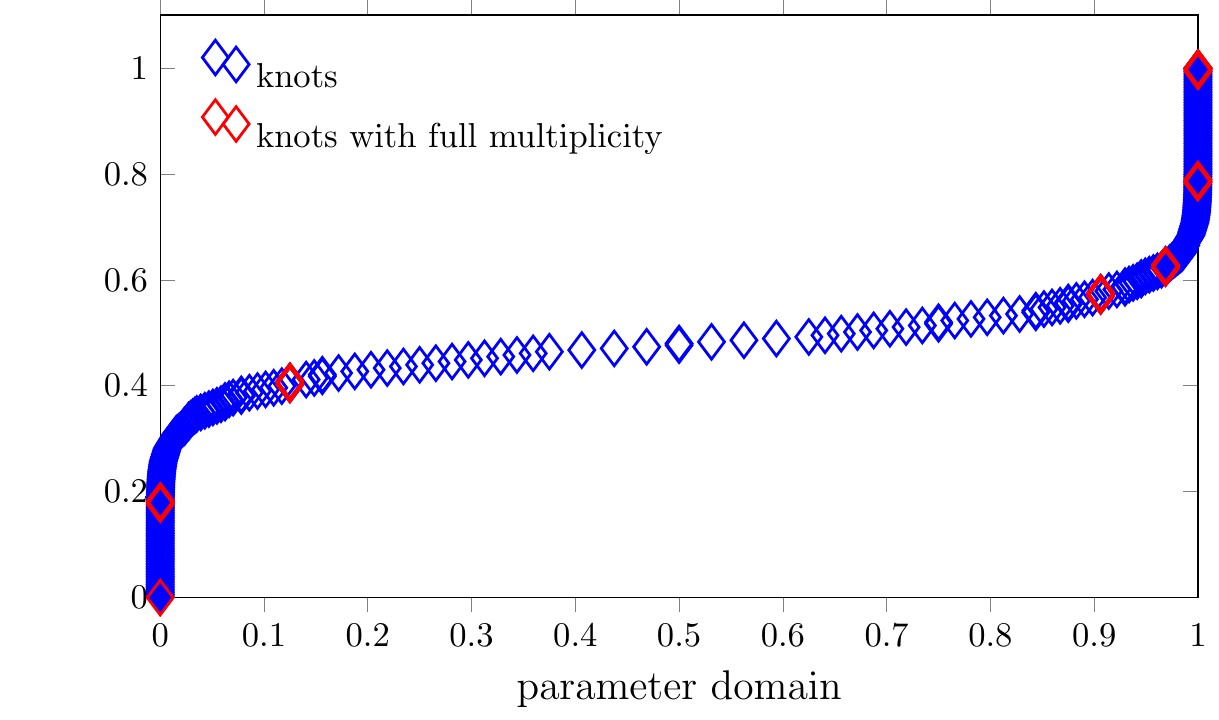}\quad
	\includegraphics[width=.475\textwidth,clip=true]{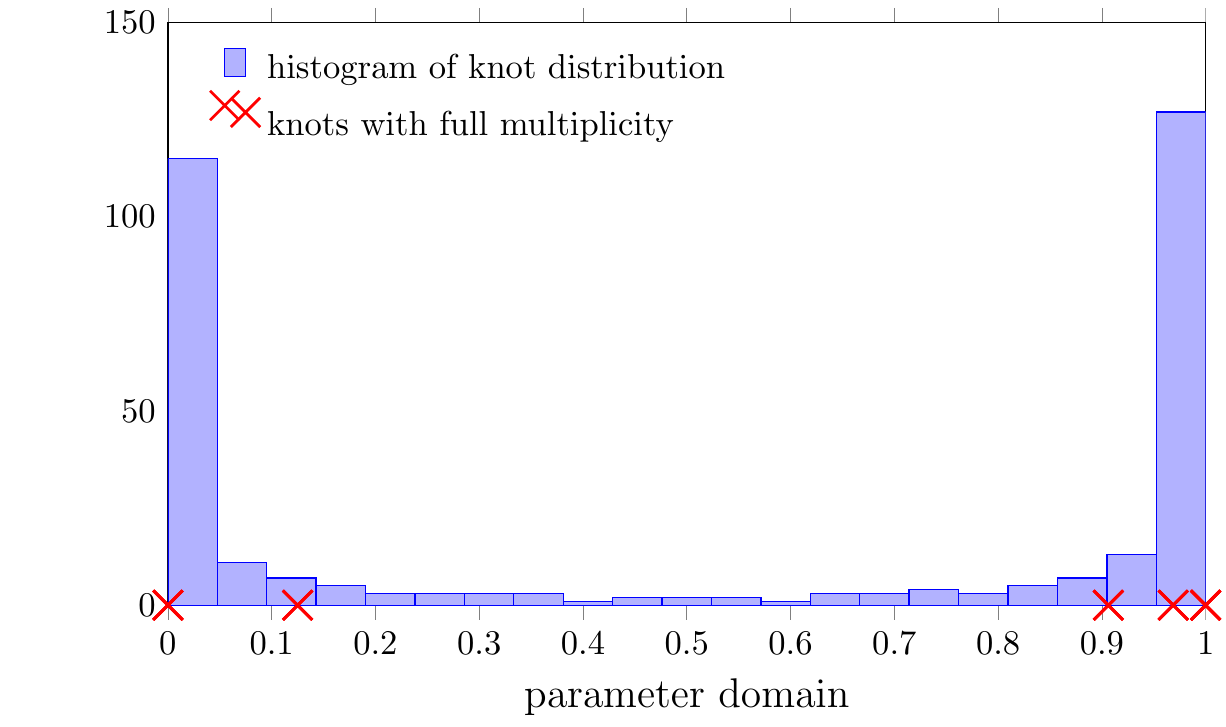}
	\caption{Example from Section~\ref{section:weak_crack_slit}: Distribution and histogram of the knots of the first  step of Algorithm~\ref{the algorithm}, where the error estimator is below $10^{-6}$ without multiplicity increase (top) as well as with mulitplicity increase (bottom) for  splines of order $p=2$ and the weighted-residual estimator with $\theta=0.75$.}
	\label{fig:weak_slit_dist}
\end{figure}

In Figure~\ref{fig:weak_slit_dist}, the knot distribution (i.e., the relative number of knots lower or equal than the parameter points on the $x$-axis) as well as a histogram of the knots in the parameter domain $[a,b]$ of the first step of Algorithm~\ref{the algorithm}, where the residual error estimator is below $10^{-6}$ are plotted over the parameter domain. 
For splines of degree $p=2$ and the weighted-residual estimator with $\theta=0.75$, we  compare Algorithm~\ref{the algorithm} without (top) and with (bottom) multiplicity increase. 
We see that for this example the difference between these two approaches is smaller than for the example of Section~\ref{section:weak_singular_pacman}. Since the solution $\phi$ is continuous in $(a,b)$, the knot distribution looks quite similar for both cases. However, we can see that multiplicity increase takes place nonetheless (bottom) and Algorithm~\ref{the algorithm} refines towards the two tips $x=\pm 1$, where the solution has singularities.

\begin{figure}
\centering
	\small{\hspace*{10mm} Slit geometry}\\
	\vspace{2mm}
	\begin{minipage}{.475\textwidth}
	\centering
	\tiny{\hspace{15mm} Faermann estimator}
	\vspace{2mm}
	\end{minipage}\quad
	\begin{minipage}{.475\textwidth}
	\centering
	\tiny{\hspace{15mm} Weighted-residual estimator}
	\vspace{2mm}
	\end{minipage}	
	\tiny{\hspace*{12mm} Error}\\
	\vspace{2mm}
	\includegraphics[width=.475\textwidth,clip=true]{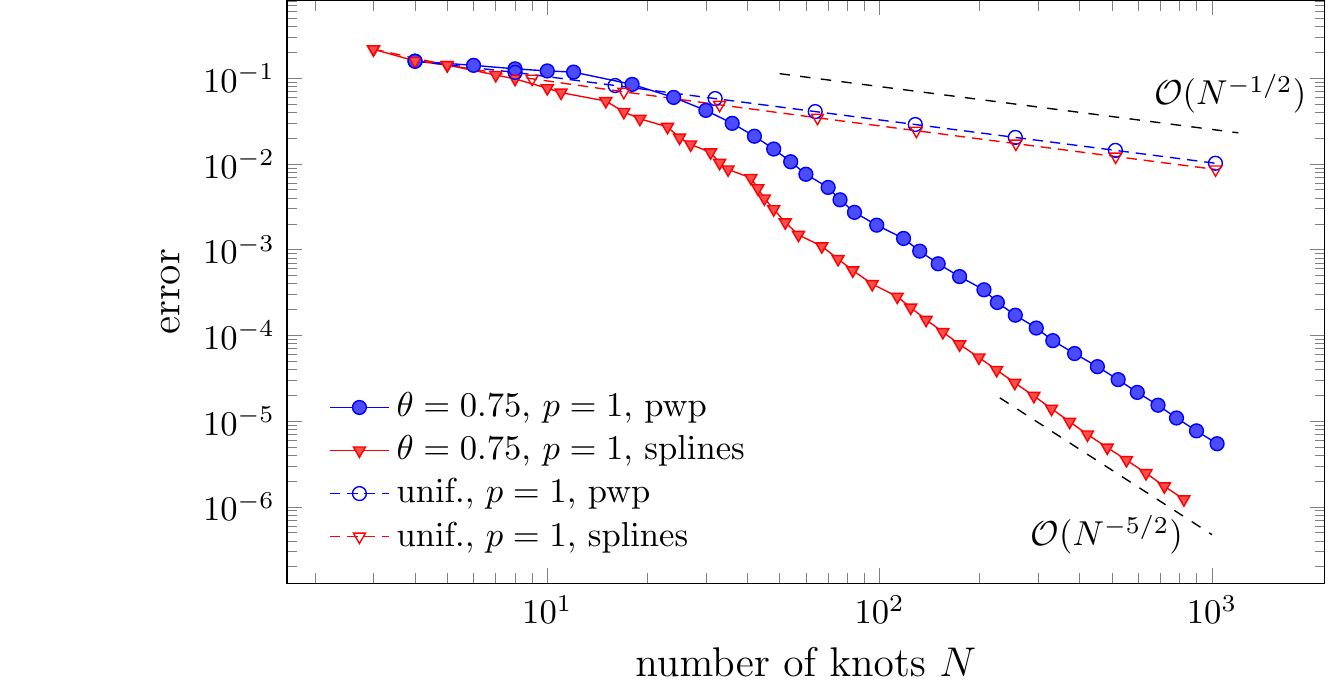}\quad
	\includegraphics[width=.475\textwidth,clip=true]{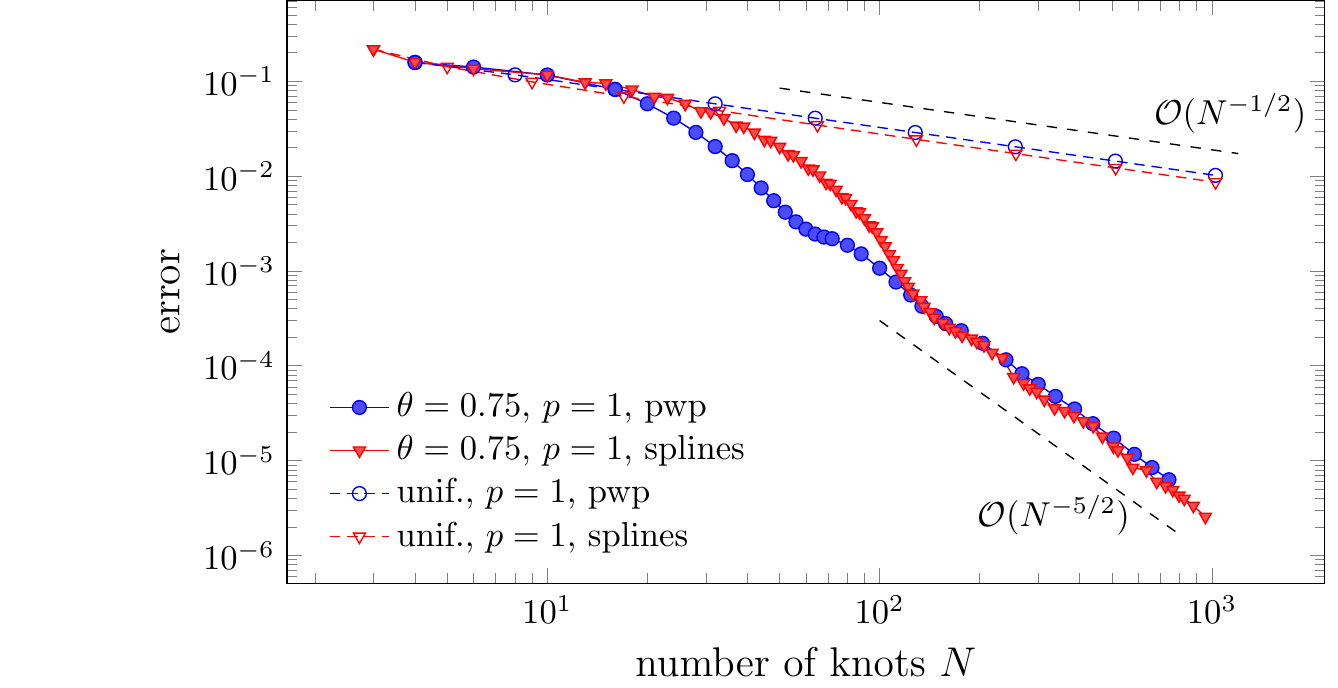}
	\vspace{2mm}
	\tiny{\hspace*{12mm} Piecewise polynomials}\\
	\vspace{2mm}
	\includegraphics[width=.475\textwidth,clip=true]{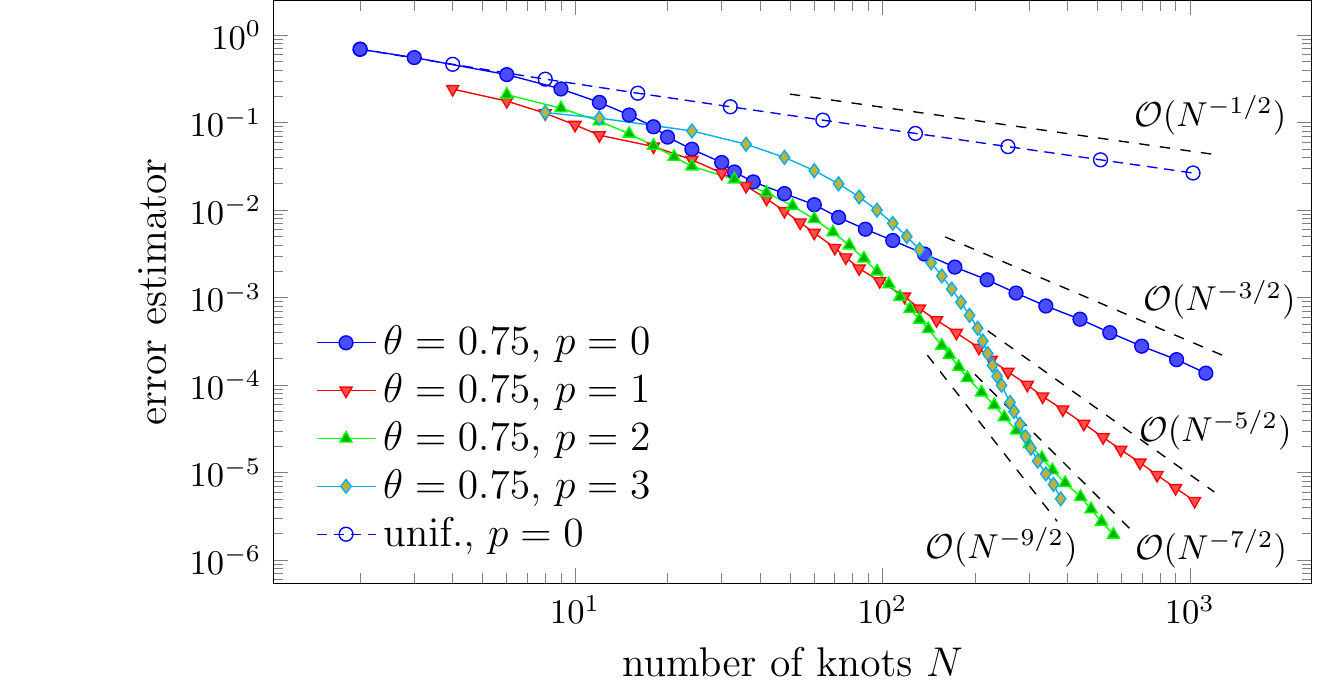}\quad
	\includegraphics[width=.475\textwidth,clip=true]{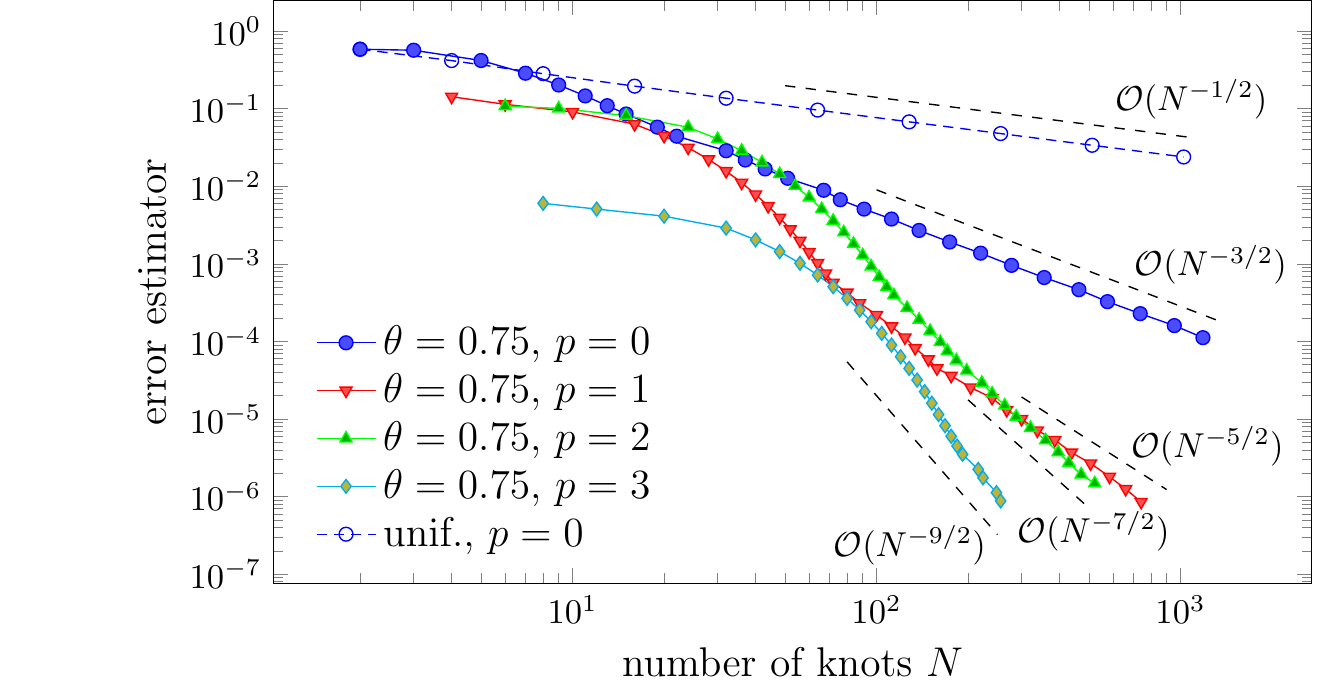}
	\vspace{2mm}
	\tiny{\hspace*{12mm} Splines}\\
	\vspace{2mm}
	\includegraphics[width=.475\textwidth,clip=true]{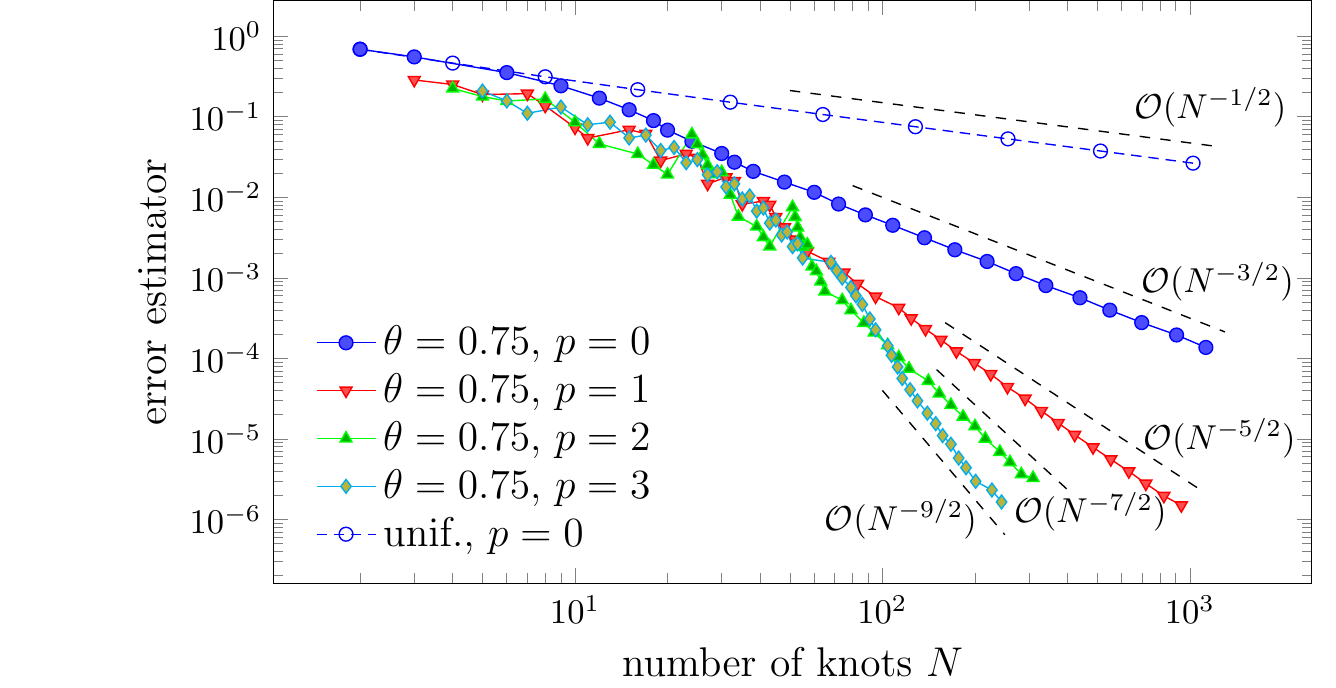}\quad
	\includegraphics[width=.475\textwidth,clip=true]{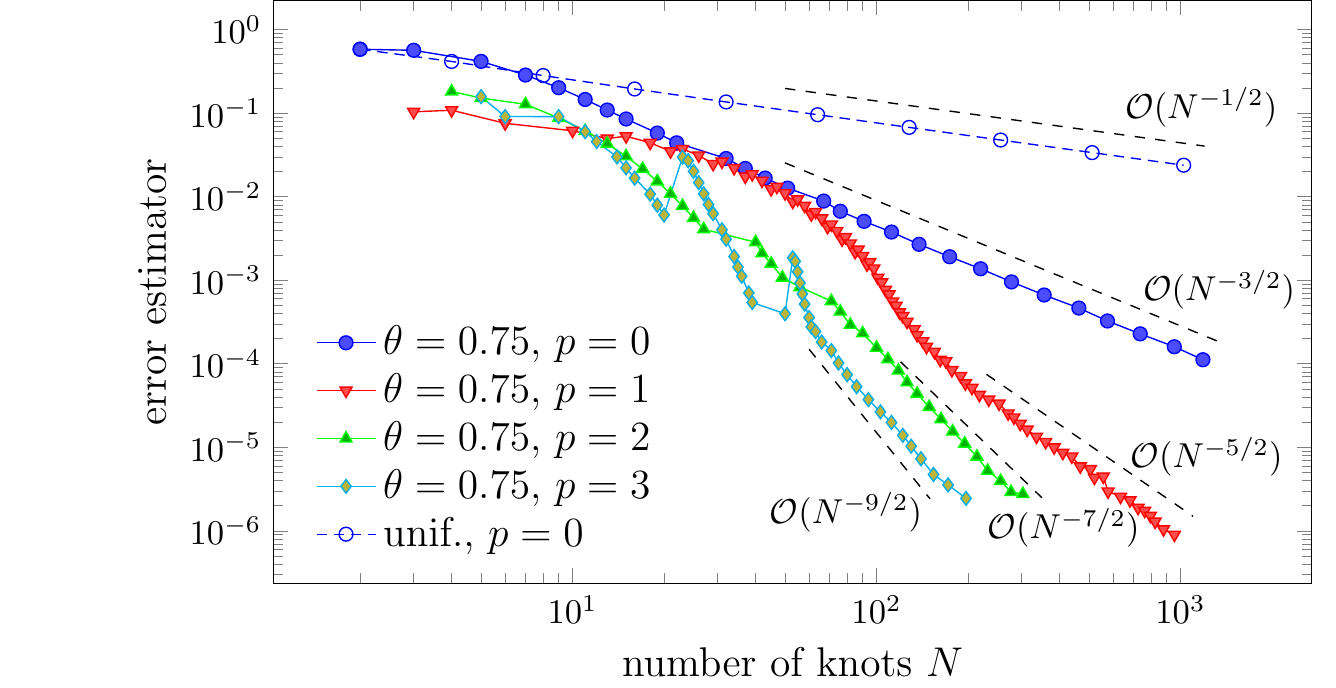}

	\caption{Example from Section~\ref{section:weak_crack_slit}: (Approximate) error $\norm{\phi-\Phi_\ell}{\mathfrak V}$ and error estimator $\eta_\ell$ with respect to the number of knots $N$ for $\theta\in\{0.75,1\}$, $p\in\{0,1,2,3\}$, different ansatz spaces (piecewise polynomials, splines), and different estimators (Faermann, weighted-residual). Left: Faermann estimator. Right: Weighted-residual estimator. Top: Comparison of $p\in\{0,1,2,3\}$ for piecewise polynomials. Middle: Comparison of $p\in\{0,1,2,3\}$ for splines. Bottom: Comparison of $p\in\{1,2,3\}$ for different ansatz spaces.}
	\label{fig:weak_slit_comparison}
\end{figure}

In Figure~\ref{fig:weak_slit_comparison}, we compare uniform refinement with $\theta=1$ and adaptive refinement with $\theta=0.75$, different ansatz spaces, and different estimators for Algorithm~\ref{the algorithm}. First (top), for $p=1$ fixed, we compare  (an approximation of the) error for the different ansatz spaces  and the algorithm steered either by  the Faermann estimator (top, left) or the weighted-residual estimator (top, right). 
Instead of the exact error $\norm{\phi-\Phi_\ell}{\mathfrak V}$, we compute $\norm{\Phi_{+}-\Phi_\ell}{\mathfrak V}$,  where  $\widehat\KK_\fine$ are the uniformly refined knots, which result from $\widehat \KK_\ell$ by adding the midpoint of each element $\widehat Q\in\widehat\QQ_\coarse$ with multiplicity one (see also Algorithm~\ref{alg:refinement}).
We get similar plots, where uniform mesh refinement leads to the suboptimal rate of convergence $\OO(N^{-1/2})$ because the solution lacks regularity at the tips $x=\pm 1$. However, for all adaptive cases, Algorithm~\ref{the algorithm} regains the optimal rate of convergence $\OO(N^{-3/2-p})$, where splines exhibit a better multiplicative constant than piecewise polynomials. Next (middle), we consider  the estimators for piecewise polynomials. For both the Faermann estimator (middle, left) and the weighted-residual estimator (middle, right), we compare different orders $p\in\{0,1,2,3\}$ for the ansatz spaces. Again, uniform mesh refinement leads to the suboptimal rate $\OO(N^{-1/2})$ (only displayed for $p=0$), whereas the adaptive strategy regains the optimal order of convergence $\OO(N^{-3/2-p})$. Lastly, we get similar results for splines (bottom).


\section{\texttt{IGABEM2D} for hyper-singular integral equation}\label{sec:hyp}

\subsection{Sobolev spaces}

Besides the Sobolev spaces from Section~\ref{sec:Sobolev spaces}, the treatment of the closed boundary $\Gamma=\partial\Omega$ in the hyper-singular case
requires the definition of $H^{\pm \sigma}_0(\partial\Omega) := \set{u\in H^{\pm \sigma}(\partial\Omega)}{\dual{u}{1}_{\partial\Omega} = 0}$
for $\sigma\in\{0,1/2,1\}$.

\subsection{Hyper-singular integral equation}
\label{sec:hypsing}
For $\sigma\in\{0,1/2,1\}$, the hyper-singular integral operator 
$\W:\H^{\sigma}(\Gamma)\to H^{\sigma-1}(\Gamma)$  and the adjoint double-layer operator $\K':\H^{\sigma-1}(\Gamma)\to H^{\sigma-1}(\Gamma)$ are well-defined, linear, and continuous. 
Moreover, $\K'$ is indeed the adjoint operator of $\K$.

For $\Gamma\subsetneqq\partial\Omega$ and $\sigma=1/2$, 
$\W:\H^{1/2}(\Gamma)\to H^{-1/2}(\Gamma)$ is symmetric and 
elliptic. 
Hence,
\begin{align}
\edualW{u}{v}:=\dual{\W u}{v}_\Gamma
\end{align}
 defines an equivalent scalar product on  
$\H^{1/2}(\Gamma)$ with corresponding norm $\enorm{\cdot}$.

For  $\Gamma=\partial\Omega$, the 
operator $\W$ is symmetric and elliptic up to the constant functions, i.e., 
$\W:H^{1/2}_0(\partial\Omega)\to H^{-1/2}_0(\partial\Omega)$ is 
elliptic. In particular, 
\begin{align}
\edualW{u}{v}:=\dual{\W u}{v}_{\partial\Omega} + \dual{u}{1}_{\partial\Omega}\dual{v}{1}_{\partial\Omega}
\end{align}
is an equivalent scalar product on $H^{1/2}(\partial\Omega)=\H^{1/2}(\partial\Omega)$  with corresponding energy norm $\enorm{\cdot}$.
Additionally, there holds the mapping property $\K':H_0^{-1/2}(\partial\Omega)\to H_0^{-1/2}(\partial\Omega)$.

Let either $g=\phi$ for some $\phi\in H^{-1/2}(\Gamma)$ with $\phi\in H_0^{-1/2}(\partial\Omega)$ if $\Gamma=\partial\Omega$, or $g=(1/2-\K')\phi$ for some $\phi\in \H^{-1/2}(\Gamma)$ with $\phi\in H^{-1/2}_0(\partial\Omega)$ if $\Gamma = \partial\Omega$.
Then, the strong form~\eqref{eq:hyper strong} is equivalently stated as
\begin{align}
\label{eq:weakform2}
 \edualW{u}{v} = \dual{g}{v}_\Gamma
 \quad\text{for all }v\in \H^{1/2}(\Gamma).
\end{align}
The  Lax--Milgram lemma applies and hence \eqref{eq:hyper strong}  admits a unique solution $u\in \H^{1/2}(\Gamma)$.
If $g=\phi$, the approach is called indirect, otherwise if $g=(1/2-\K')\phi$, it is called direct.
Details and proofs are found, e.g.,  in  \cite{mclean00,hw08,steinbach08,ss11}.

\subsection{Galerkin IGABEM}
\label{sec:spaces} 
Let $p\in\N$ be a fixed polynomial degree. 
Moreover, let $\widehat\KK_\coarse$ and $\WW_\coarse$ be knots and weights as in Section~\ref{sec:splines on Gamma}.
We suppose that the initial knots $\widehat\KK_0$ additionally satisfy that $\#_0 t_{0, i}\le p$ for $i\in\{1,\dots,N_0-p\}$ and that $w_{0,1-p}=w_{0,N_0-p}$ if $\Gamma=\partial\Omega$.
Note that \eqref{eq:interpolatoric} shows that
\begin{align*}
\widehat B_{0, 1-p}(a)=\widehat B_{0, N_0-p}(b-)
=1
=\widehat B_{\coarse, 1-p}(a)=\widehat B_{\coarse, N_\coarse-p}(b-).
\end{align*}
If $\Gamma=\partial\Omega$, the assumption $w_{0,1-p}=w_{0,N_\gamma-p}$ thus ensures that $\widehat W_0(a)=\widehat W_0(b-)$, which further yields that  $w_{\coarse, 1-p}=w_{\coarse, N_\coarse-p}$.

We introduce the  ansatz space 
\begin{align}\label{eq:hypsing X0}
\YY_\coarse:=\begin{cases}\set{V_\coarse\in\SS^p({\mathcal{K}}_\coarse, \mathcal{W}_\coarse)}{ V_\coarse(\gamma(a))=V_\coarse(\gamma(b-))}\subset H^{1}(\Gamma)&\text{ if }\Gamma=\partial\Omega,\\
\set{V_\coarse\in\SS^p({\mathcal{K}}_\coarse, \mathcal{W}_\coarse)}{0= V_\coarse(\gamma(a))=V_\coarse(\gamma(b-))}\subset \H^{1}(\Gamma)&\text{ if }\Gamma\subsetneqq\partial\Omega.
\end{cases}
\end{align}
According to Section~\ref{sec:splines},
it holds that
\begin{align}\label{eq:hypsing basis}
\YY_\coarse=\begin{cases}{\rm span}\Big(\set{  R_{\coarse, i,p}}{i=2-p,\dots,N_\coarse-p-1}\cup \{  R_{\coarse, 1-p,p}+  R_{\coarse, N_\coarse-p,p}\}\Big)&\text{ if }\Gamma=\partial\Omega,\\
{\rm span}\set{  R_{\coarse, i,p}}{i=2-p,\dots,N_\coarse-p-1}&\text{ if }\Gamma\subsetneqq\partial\Omega.
\end{cases}
\end{align}
In each case, the corresponding sets even form a basis of $\YY_\coarse$. 
We abbreviate the (continuous) basis functions $R^{\rm c}_{\coarse, i,p}:=R_{\coarse, i,p}$ for $i=2-p,\dots,N_\coarse-p-1$ and  $ R^{\rm c}_{\coarse, i,p}:=R_{\coarse, 1-p,p}+  R_{\coarse, N_\coarse-p,p}$ for $i=1-p$.

We  assume the additional regularity $\phi\in L^2(\Gamma)$ and  approximate $\phi$ by $\phi_\coarse:=\Pi_\coarse \phi$, where $\Pi_\coarse$ is the $L^2$-orthogonal projection onto the space of transformed piecewise polynomials 
\begin{align}
\PP^p(\QQ_\coarse):=\set{\widehat \Psi_\coarse\circ\gamma^{-1}}{\widehat\Psi_\coarse \text{ is }\text{pcw.\ polynomial of degree } p \text{ with breakpoints }\widehat\NN_\coarse }.
\end{align}
Note that $\dual{\phi}{1}_{\partial\Omega}=0$ in case of $\Gamma=\partial\Omega$ also yields that $\dual{\phi_\coarse}{1}_{\partial\Omega}=0$.
With $g_\coarse:=\phi_\coarse$ resp.\ $g_\coarse:=(1/2-\K')\phi_\coarse$ , the corresponding Galerkin approximation $U_\coarse\in\YY_\coarse$ reads 
\begin{align}\label{eq:hypsing Galerkin}
 \edualW{U_\coarse}{V_\coarse} = \dual{g_\coarse}{V_\coarse}_\Gamma 
\,\text{ for all }V_\coarse\in \YY_\coarse. 
\end{align}
At least for the direct method, we empirically observed that working with $g$ instead of with $g_\coarse$ leads to strong implementational instabilities (since a possible singularity of $g$ meets the singularity of the integral kernels of the boundary integral operators.)

Note that \eqref{eq:hypsing Galerkin} is equivalent to solving the finite-dimensional linear system
\begin{align}\label{eq:Galerkin_hyper}
\big(\edualW{R^{\rm c}_{\coarse,i',p}}{R^{\rm c}_{\coarse,i,p}}\big)_{i,i'=k-p}^{N_\coarse-p-1} \cdot (c_{\coarse,i'})_{i'=k-p}^{N_\coarse-p} = \big(\dual{g_\coarse}{R^{\rm c}_{\coarse,i,p}}_\Gamma\big)_{i=k-p}^{N_\coarse-p-1}, 
\end{align}
where $k=1$  for $\Gamma=\partial\Omega$ and $k=2$ for $\Gamma\subsetneqq\partial\Omega$.
Then, $U_\coarse = \sum_{i'=k-p}^{N_\coarse-p} c_{\coarse,i'} R^{\rm c}_{\coarse,i',p}$.
The computation of the matrix and the right-hand side vector in~\eqref{eq:Galerkin_hyper} is realized in {\tt WMatrix.c} and {\tt RHSVectorHyp.c}, where the required projection $\Pi_\coarse$ is realized in {\tt PhiApprox.c}.
They can be called in {\sc Matlab} via {\tt buildWMatrix}, {\tt buildRHSVectorHyp}, and {\tt buildPhiApprox}; see Appendix~\ref{sec:overview}.  
For details on the implementation, we refer to Section~\ref{sec:matrices} and \ref{sec:right vectors}.

\subsection{\textsl{A~posteriori} error estimation}
\label{sec:a_posteriori_hyper}

In the hyper-singular case, we consider the following two different node-based estimators:
the $(h-h/2)$-estimator
\begin{align}\label{eq:hhW}
\eta_{\W,\rm hh2,\coarse}^2:=\sum_{x\in\NN_\coarse} \eta_{\W,\rm hh2,\coarse}(x)^2\quad\text{with}\quad  \eta_{\W,\rm hh2,\coarse}(x)^2:= \norm{h_\coarse^{1/2}\partial_\Gamma(u_\fine-u_\coarse)}{L^2(\omega_\coarse(x))}^2,
\end{align}
where  $\widehat\KK_\fine$ are the uniformly refined knots, which result from $\widehat \KK_\coarse$ by adding the midpoint of each element $\widehat Q\in\widehat\QQ_\coarse$ with multiplicity one (see also Algorithm~\ref{alg:refinement});
and the weighted-residual estimator
\begin{align}\label{eq:resW}
\eta_{\W,\rm res,\coarse}^2:=\sum_{x\in\NN_\coarse} \eta_{\W,\rm res,\coarse}(x)^2\text{ with }  \eta_{\W,\rm res,\coarse}(x)^2:= \norm{h_\coarse^{1/2}(g_\coarse-\W U_\coarse)}{L^{2}(\omega_\coarse(x))}^2,
\end{align}
where the additional regularity $\phi\in L^2(\Gamma)$ ensures that $g_\coarse -\W U_\coarse\in L^2(\Gamma)$; see Section~\ref{sec:hypsing}.
To account for the approximation of $g$ by $g_\coarse$, we additionally consider the oscillations
\begin{align}\label{eq:osc}
\osc_{\W,\coarse}^2:=\sum_{x\in\NN_\coarse} \osc_{\W,\coarse}(x)^2 
\quad\text{with}\quad
\osc_{\W,\coarse}(x):=
\norm{h_\coarse^{1/2}(1-\Pi_\coarse)\phi}{L^2(\omega_\coarse(x))}^2.
\end{align}

The computation of the estimators and oscillations is implemented in \texttt{HHEstHyp.c}, \texttt{ResEstHyp.c}, and \texttt{OscHyp.c}. 
They can be called in {\sc Matlab} via \texttt{buildHHEstHyp}, \texttt{buildResEstHyp}, and \texttt{buildOscHyp}.
The residual estimators require the evaluation of the boundary integral operators $\mathfrak{W}$ and $\mathfrak{K}'$ applied to some function, which is implemented in \texttt{ResidualHyp.c}. 
The evaluations can be used in {\sc Matlab} via \texttt{evalW.c} and \texttt{evalRHSHyp.c}; see Appendix~\ref{sec:overview}.
The implementation of the estimators and the evaluations are discussed in Section~\ref{sec:implementation}.

\subsection{Adaptive IGABEM algorithm}

To enrich a given ansatz space $\XX_\coarse$, we compute one of the error estimators $\eta_\coarse$ and determine a set of marked nodes with large indicator $\eta_\coarse(x)$.
By default, we then apply an adapted version of the refinement Algorithm~\ref{alg:refinement} with the following two obvious modifications: 
First, we only consider $p\in\N$ instead of $p\in\N_0$ as input.
Second, in {\rm (iii)}, we only increase the multiplicity if it is less or equal to $p-1$ to ensure that the basis functions are at least continuous.
As before, the required marking is realized in {\tt markNodesElements.m}. 
Together with the computation of the new weights, the $h$-refinement is realized in  {\tt refineBoundaryMesh.m}, while the multiplicity increase is realized in {\tt increaseMult.m}.
We stress that we have also implemented two further relevant strategies that rely on $h$-refinement only:
Replace {\rm (iii)} by adding all elements $Q\in\QQ_\coarse$ containing one of the other nodes in $\MM_\coarse$ to $\MM_{\coarse}'$;
or replace {\rm (iii)} by adding all elements $Q\in\QQ_\coarse$ containing one of the other nodes in $\MM_\coarse$ to $\MM_{\coarse}'$
and  insert the midpoints in {\rm (iv)} with multiplicity $p$.  
The first strategy leads to refined splines of full regularity, whereas the second one leads to lowest possible regularity, i.e., mere continuity.

We fix the considered error estimator $\eta_\coarse\in \big\{(\eta_{\W,\rm hh2,\coarse}^2+ \osc_\coarse^2)^{1/2},$ $(\eta_{\W,\rm res,\coarse}^2+\osc_\coarse^2)^{1/2}\big\}$. 
The corresponding error indicators $\eta_\coarse(x)$ are defined accordingly.

\begin{algorithm}\label{the algorithm2}
\textbf{Input:} Polynomial degree $p\in \N$, initial knot vector $\widehat\KK_0$, initial weights~$\mathcal{W}_0$, marking parameter $0<\theta\le1$.\\
\textbf{Adaptive loop:} For each $\ell=0,1,2,\dots$ iterate the following steps {\rm(i)--(iv)}:
\begin{itemize}
\item[\rm(i)] Compute  approximation $U_\ell\in\XX_\ell$ by solving~\eqref{eq:hypsing Galerkin}.
\item[\rm(ii)] Compute refinement indicators $\eta_\ell({x})$
for all nodes ${x}\in\NN_\ell$.
\item[\rm(iii)] Determine a minimal set of nodes $\MM_\ell\subseteq\NN_\ell$ such that
\begin{align}\label{eq:Doerfler}
 \theta\,\eta_\ell^2 \le \sum_{{x}\in\MM_\ell}\eta_\ell({x})^2.
\end{align}
\item[\rm(iv)] Generate refined knot vector $\widehat\KK_{\ell+1}$ via the adapted version of Algorithm~\ref{alg:refinement}.
\end{itemize}
\textbf{Output:} Approximate solutions $U_\ell$ and estimators $\eta_\ell$ for all $\ell \in \N_0$.
\end{algorithm}

The adaptive algorithm is realized in {\tt IGABEMHyp.m}. 
The considered problem as well as the used parameters can be changed there as pleased.
In Section~\ref{sec:implementation}, we discuss how the arising (singular) boundary integrals are computed.
We also mention that Algorithm~\ref{the algorithm2}~(iii) is realized in {\tt markNodesElements.m}, which also determines the corresponding set of marked elements $\MM_\ell'$ and the nodes whose multiplicity should be increased from Algorithm~\ref{alg:refinement}~(ii)--(iii).



\subsection{Numerical experiments}\label{section:numerics_hyp}

\begin{figure}
\centering
	\includegraphics[width=.475\textwidth,clip=true]{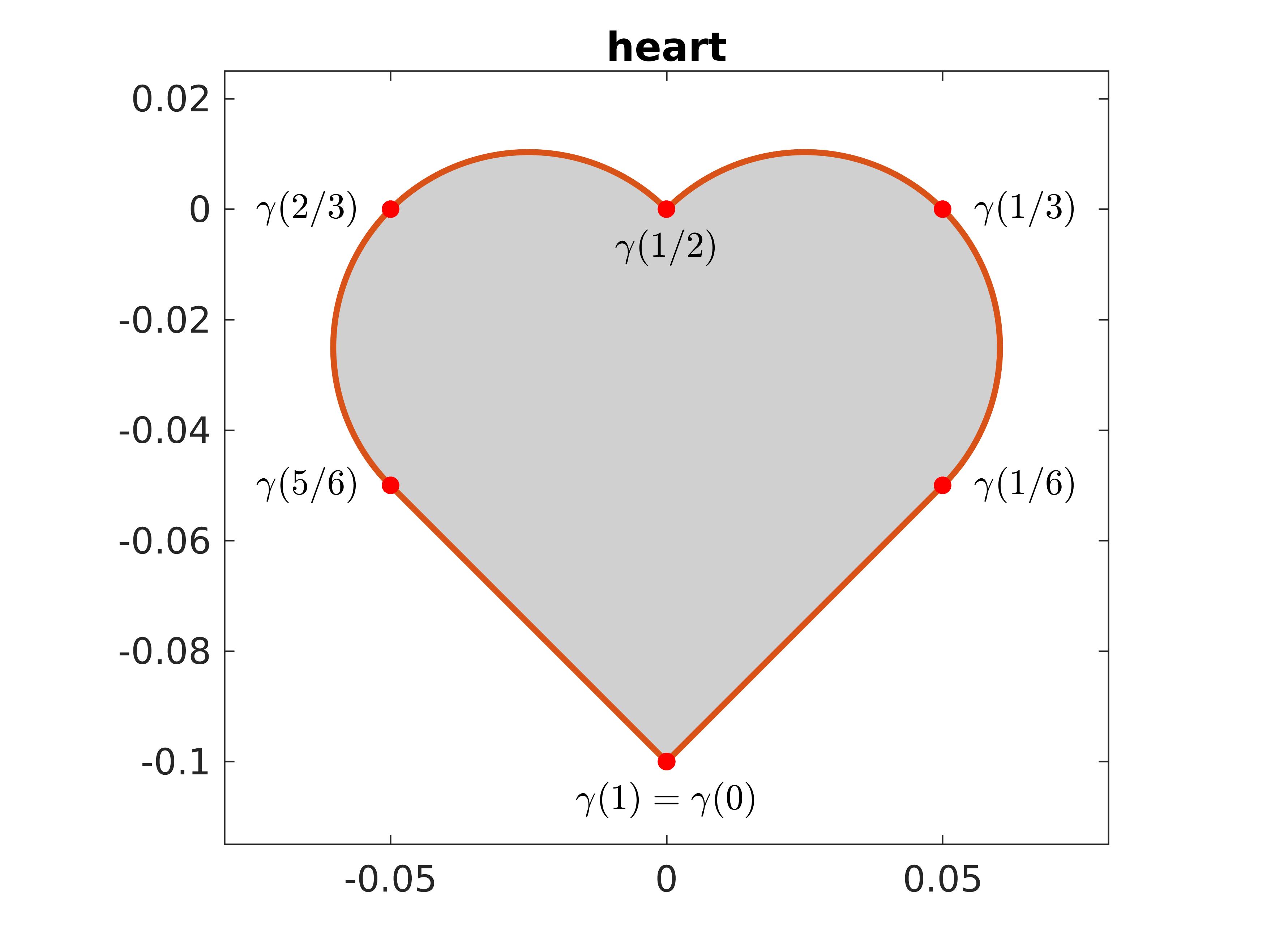}
	\caption{Geometry and initial nodes for examples from Section~\ref{section:numerics_hyp}.}
	\label{fig:geometries_hyp}
\end{figure}

In this section, we empirically investigate the performance of Algorithm~\ref{the algorithm2} for different ansatz spaces, and error estimators. 
Figure~\ref{fig:geometries_hyp} shows the different geometries, namely a circle and a heart geometry whose boundary can be parametrized via rational splines of degree~$2$. 
As initial ansatz space, we either consider the same rational splines, i.e., $\widehat\KK_0=\widehat\KK_\gamma$ and $\WW_0=\WW_\gamma$, splines of degree $p$ and smoothness $C^{p-1}$, 
i.e., 
\begin{align*}
\widehat\KK_0 =(\widehat x_{\gamma,1}, \widehat x_{\gamma,2},\dots,\widehat x_{\gamma,n_\gamma-1},\underbrace{\widehat x_{\gamma,n_\gamma},\dots,\widehat x_{\gamma,n_\gamma}}_{\# = p+1})
\end{align*} 
and $\WW_0=(1,\dots,1)$, 
or continuous piecewise polynomials of degree $p$, i.e., 
\begin{align*}
\widehat\KK_0 =(\underbrace{\widehat x_{\gamma,1},\dots, \widehat x_{\gamma,1}}_{\#=p},\dots,\underbrace{\widehat x_{\gamma,n_\gamma-1},\dots, \widehat x_{\gamma,n_\gamma-1}}_{\#=p},\underbrace{\widehat x_{\gamma,n_\gamma},\dots,\widehat x_{\gamma,n_\gamma}}_{\#=p+1})
\end{align*} 
and $\WW_0=(1,\dots,1)$.
In the latter case, we always consider $h$-refinement with new knots having multiplicity $p$ as before Algorithm~\ref{the algorithm2}.

\subsubsection{Singular problem on heart geometry}
\label{section:hyp_singular_heart}
Similarly to Section~\ref{section:weak_singular_pacman}, we prescribe an exact solution $P$ of the Laplace problem in polar coordinates $(x_1, x_2)=r(\cos\beta,\sin\beta)$ with $\beta\in(-\frac{3\pi}{2},\frac{\pi}{2})$ as
\begin{align*}
P(x_1,x_2):= r^{2/3}\cos(\tau\beta)
\end{align*}
on the heart geometry, cf.\ Figure~\ref{fig:geometries_hyp}. 
The solution $u$ of the the hyper-singular equation~\eqref{eq:hyper strong} with $g=(1/2-\K')(\partial P/\partial\nu)$ is just the restriction $P$ onto $\partial\Omega$, which has a generic singularity at the origin. 

\begin{figure}
\centering
	\small{\quad\quad Heart geometry}\\
	\vspace{2mm}
	\tiny{\hspace*{12mm} Without multiplicity increase}\\
	\vspace{2mm}
	\includegraphics[width=.475\textwidth,clip=true]{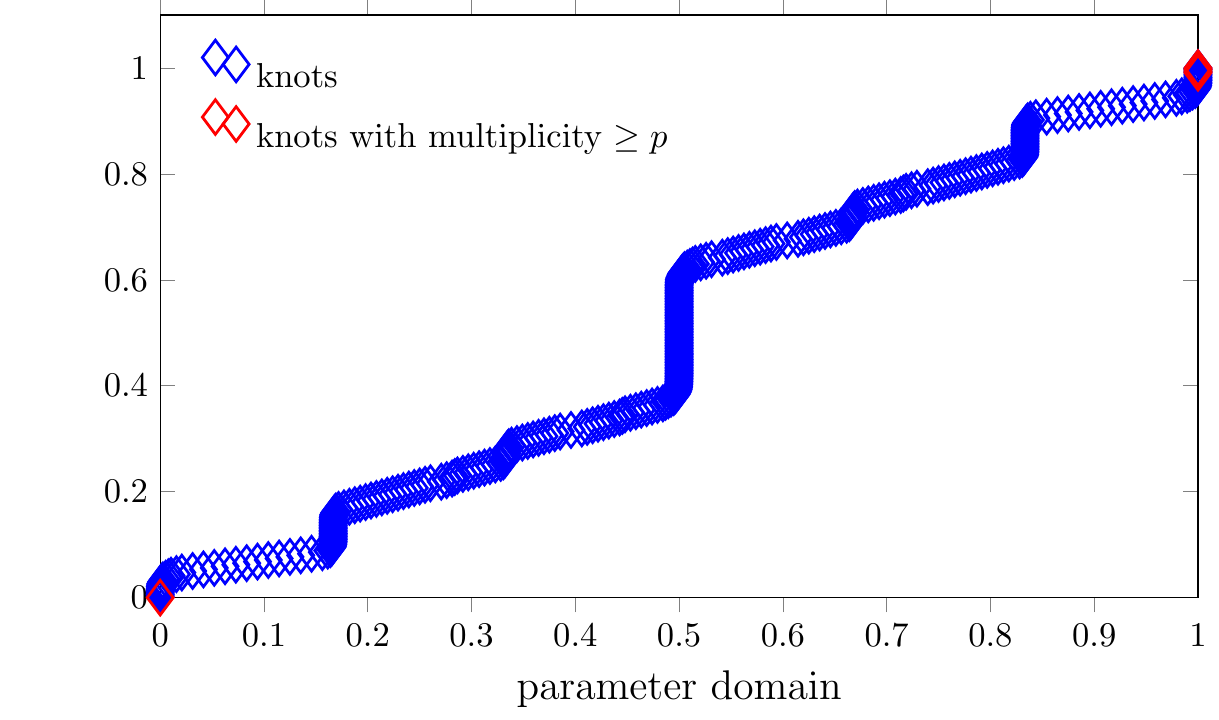}\quad
	\includegraphics[width=.475\textwidth,clip=true]{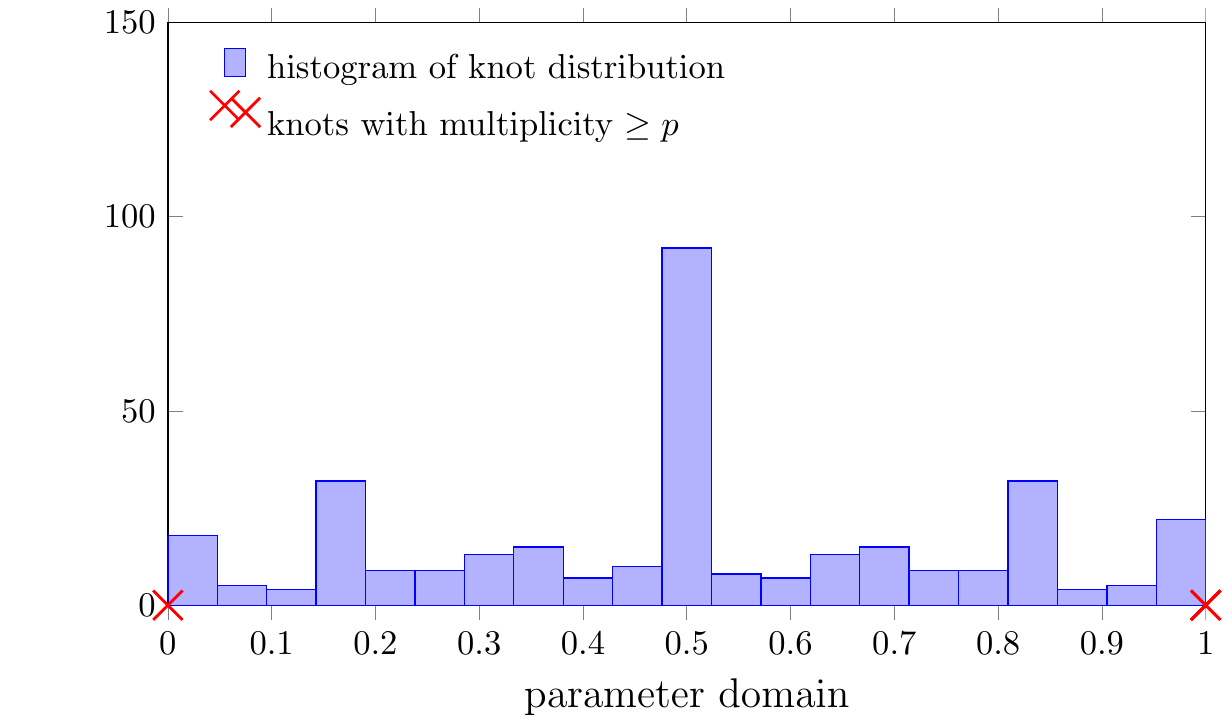}
	\tiny{\hspace*{12mm} With multiplicity increase}\\
	\vspace{2mm}
	\includegraphics[width=.475\textwidth,clip=true]{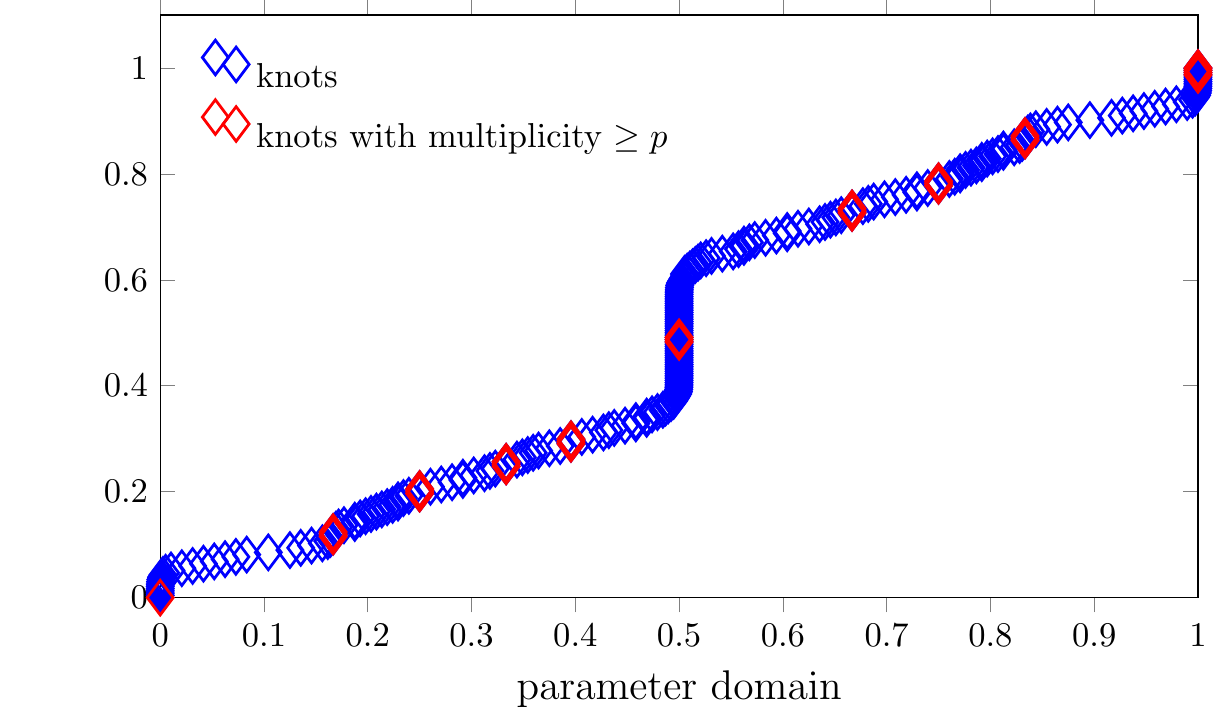}\quad
	\includegraphics[width=.475\textwidth,clip=true]{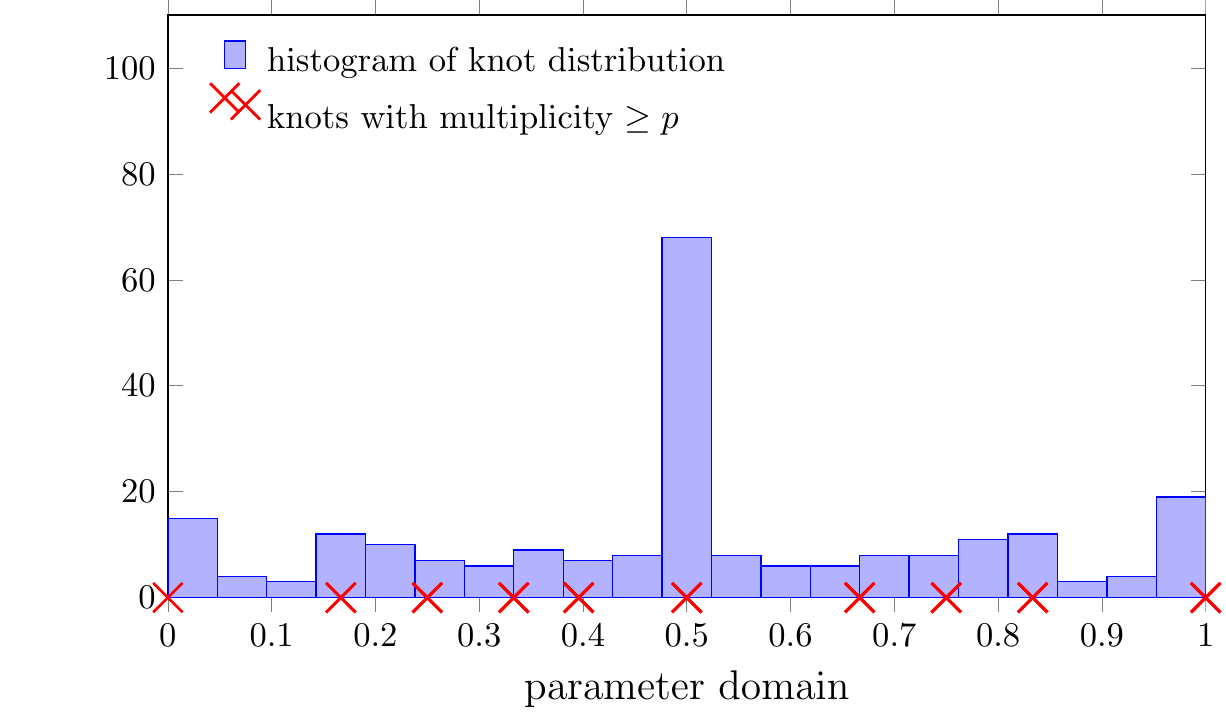}
	\caption{Example from Section~\ref{section:hyp_singular_heart}: Distribution and histogram of the knots of the first step of Algorithm~\ref{the algorithm2} where the error estimator is below $10^{-6}$ without multiplicity increase (top) as well as with mulitplicity increase (bottom) for splines of order $p=3$ and the weighted-residual estimator with $\theta=0.5$.}
	\label{fig:hyp_heart_dist}
\end{figure}

In Figure~\ref{fig:hyp_heart_dist}, the knot distribution (i.e., the relative number of knots lower or equal than the parameter points on the $x$-axis) as well as a histogram of the knots in the parameter domain $[a,b]$ of the first step of Algorithm~\ref{the algorithm2} where the residual error estimator is below $10^{-6}$ are plotted over the parameter domain. 
For splines of degree $p=3$ and the weighted-residual estimator with $\theta=0.5$, Algorithm~\ref{the algorithm2} without multiplicity increase (top) and Algorithm~\ref{the algorithm2} with multiplicity increase (bottom) heavily refine the mesh towards $\gamma(1/2)=(0,0)$, where $u$ has a generic singularity. For the latter one, we can see that multiplicity increase leads to a smoother knot distribution where less $h$-refinement takes place besides the singularity.

\begin{figure}
\centering
	\small{\hspace*{10mm} Heart geometry}\\
	\vspace{2mm}
	\begin{minipage}{.475\textwidth}
	\centering
	\tiny{\hspace{15mm} $(h-h/2)$-estimator}
	\vspace{2mm}
	\end{minipage}\quad
	\begin{minipage}{.475\textwidth}
	\centering
	\tiny{\hspace{15mm} Weighted-residual estimator}
	\vspace{2mm}
	\end{minipage}
	\tiny{\hspace*{12mm} Error}\\
	\vspace{2mm}
	\includegraphics[width=.475\textwidth,clip=true]{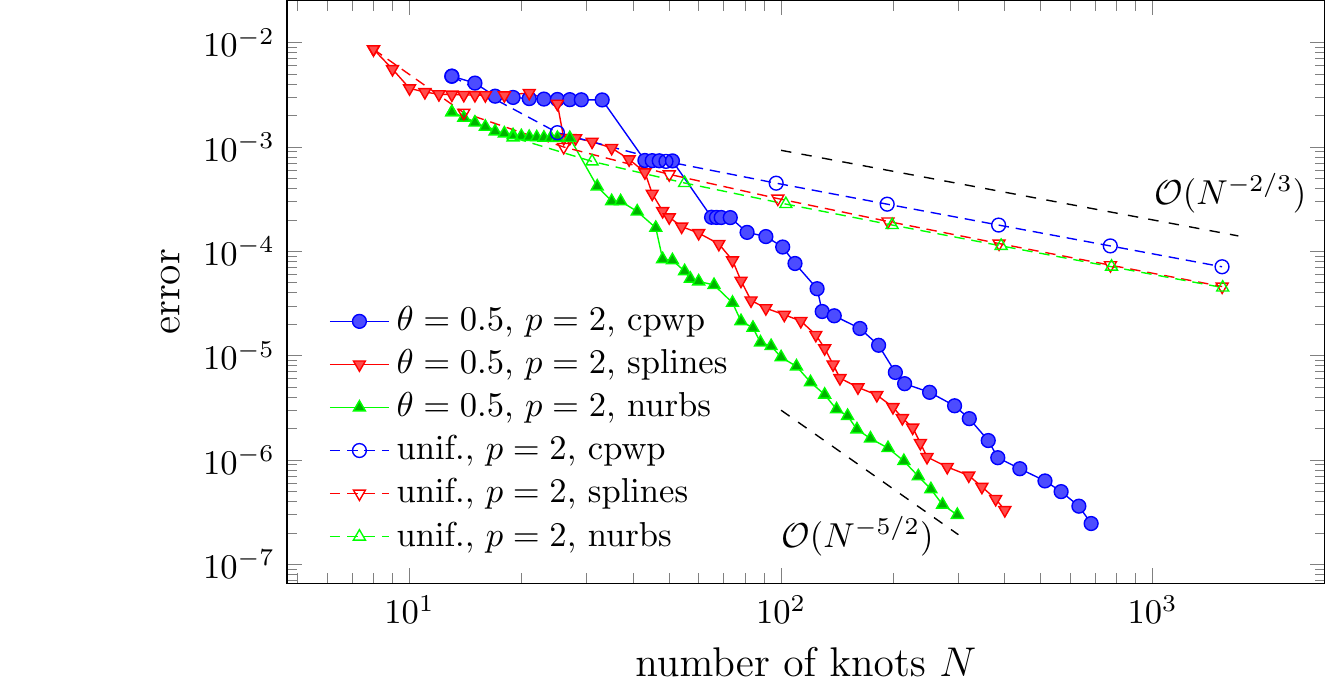}\quad
	\includegraphics[width=.475\textwidth,clip=true]{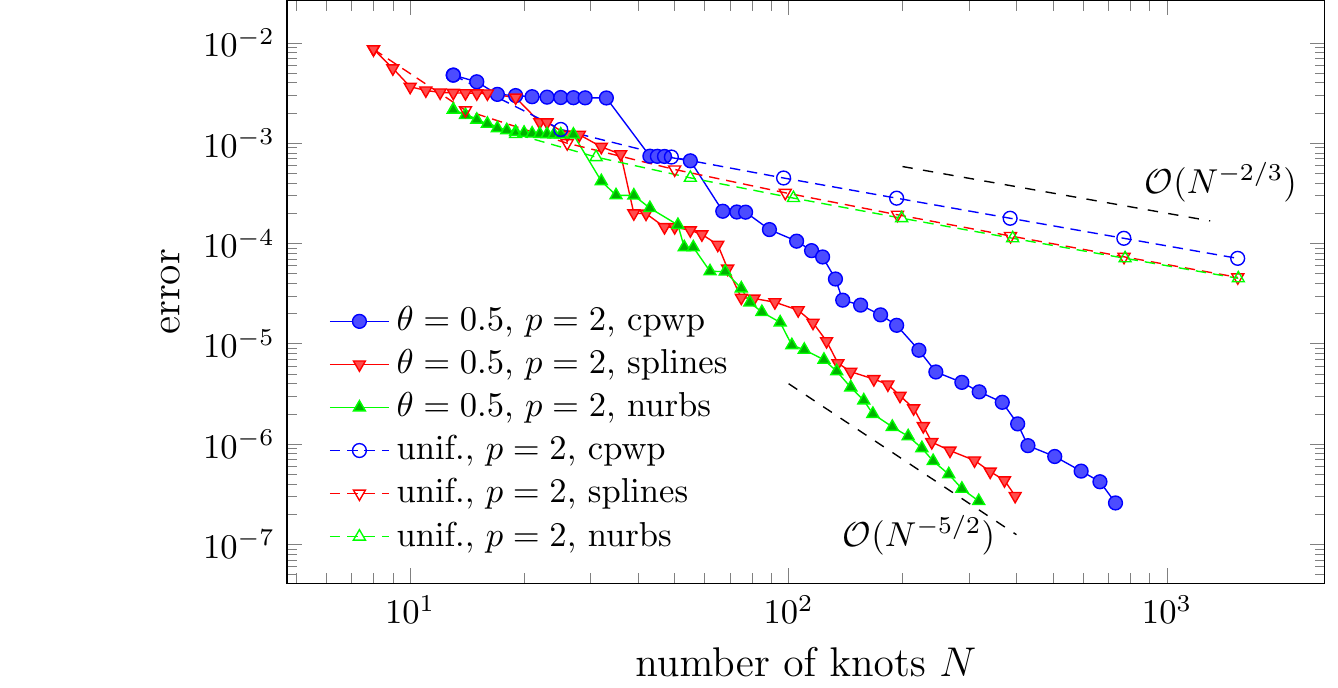}
	\vspace{2mm}
	\tiny{\hspace*{12mm} Cont.\ piecewise polynomials}\\
	\vspace{2mm}
	\includegraphics[width=.475\textwidth,clip=true]{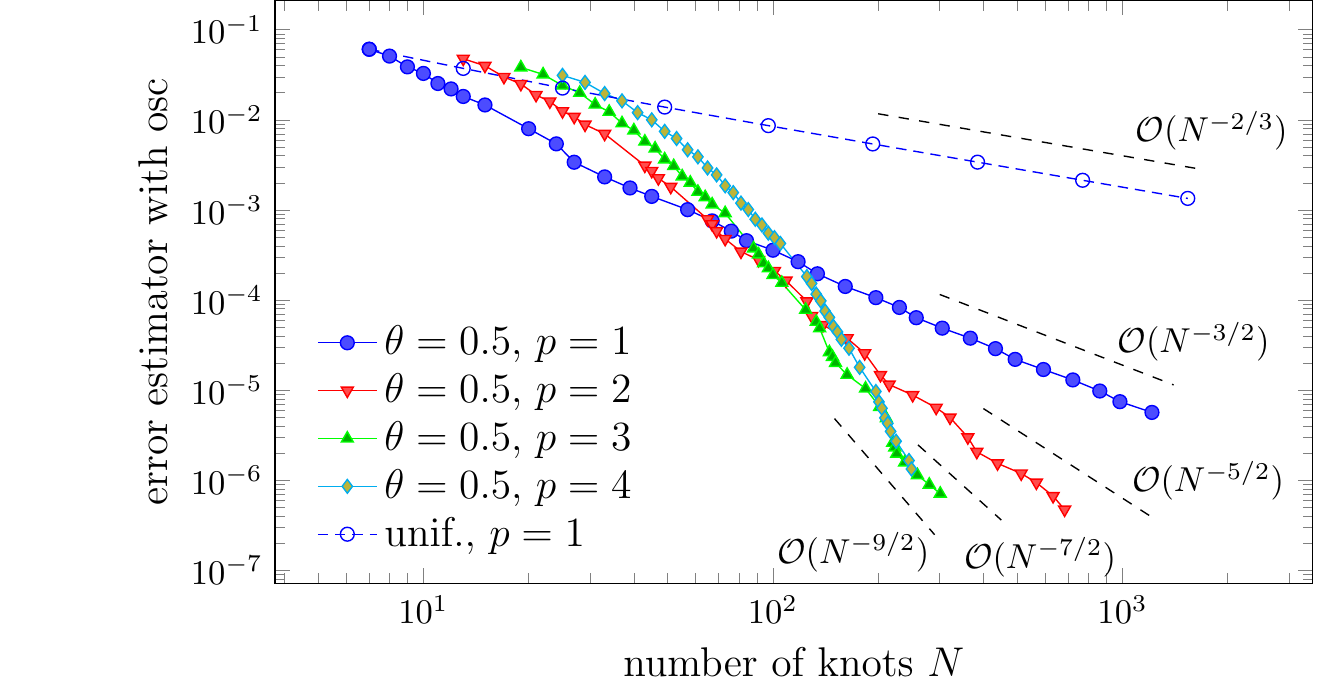}\quad
	\includegraphics[width=.475\textwidth,clip=true]{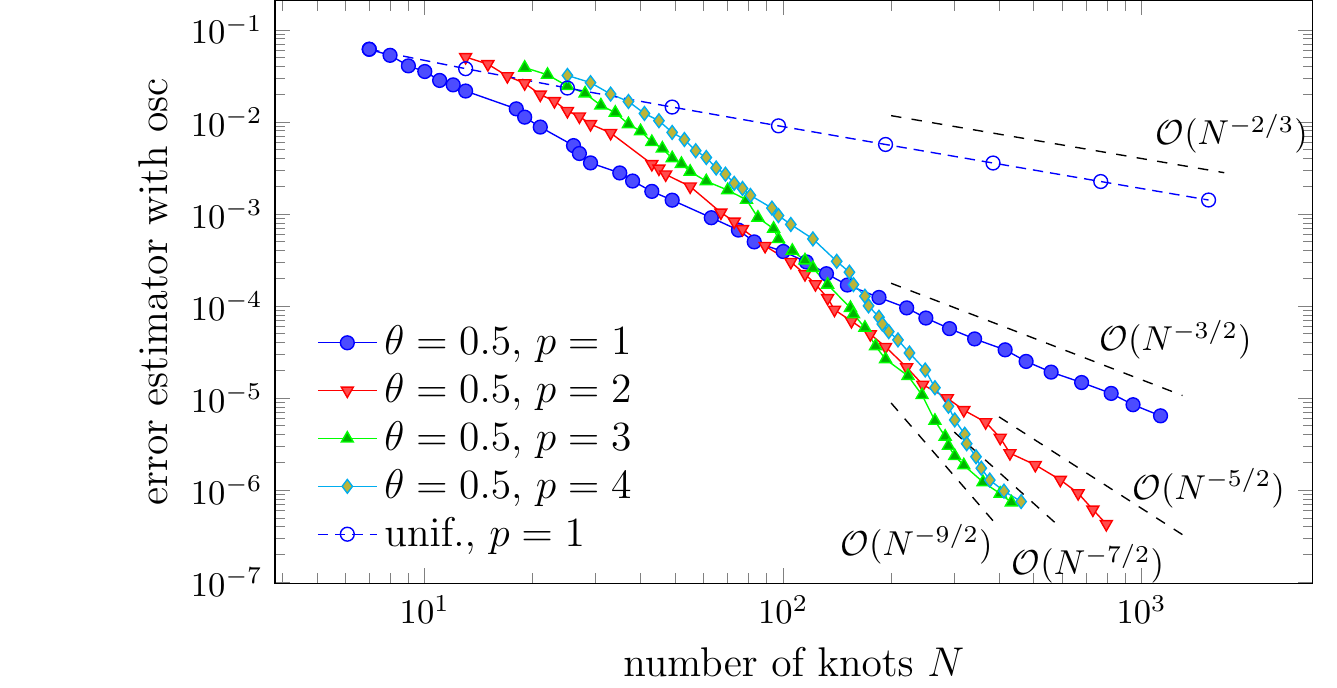}
	\vspace{2mm}
	\tiny{\hspace*{12mm} Splines}\\
	\vspace{2mm}
	\includegraphics[width=.475\textwidth,clip=true]{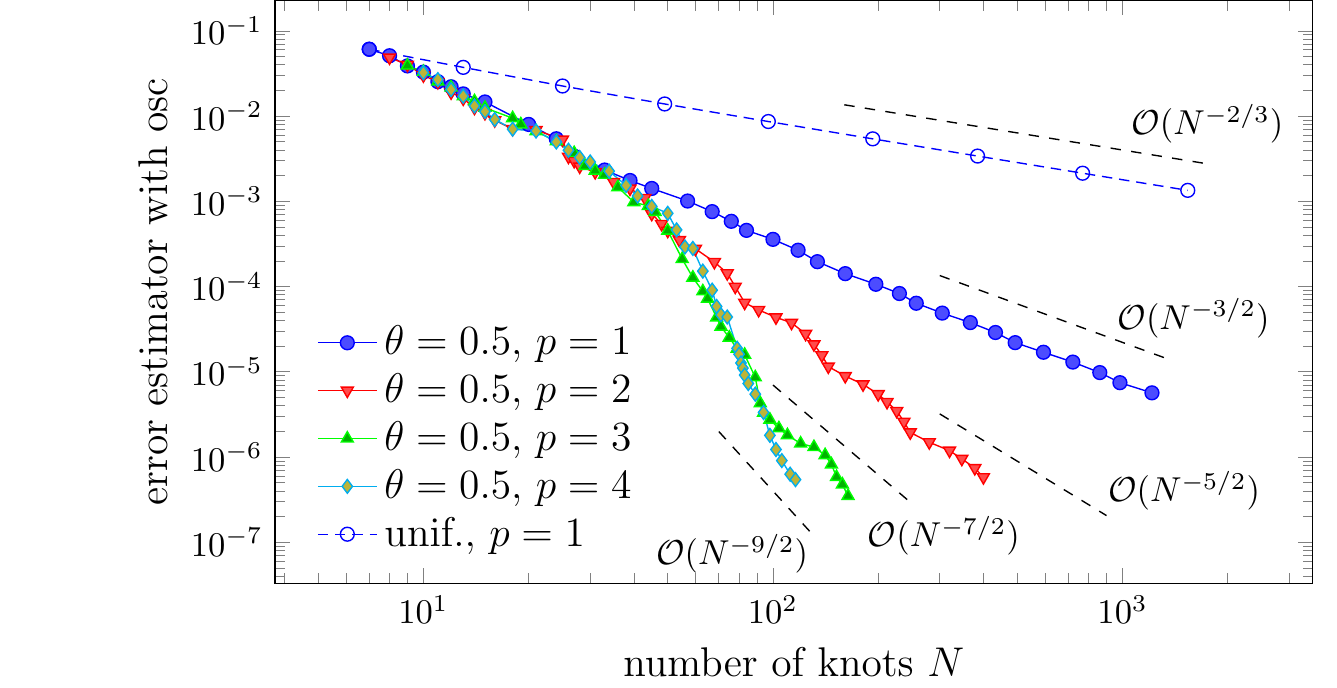}\quad
	\includegraphics[width=.475\textwidth,clip=true]{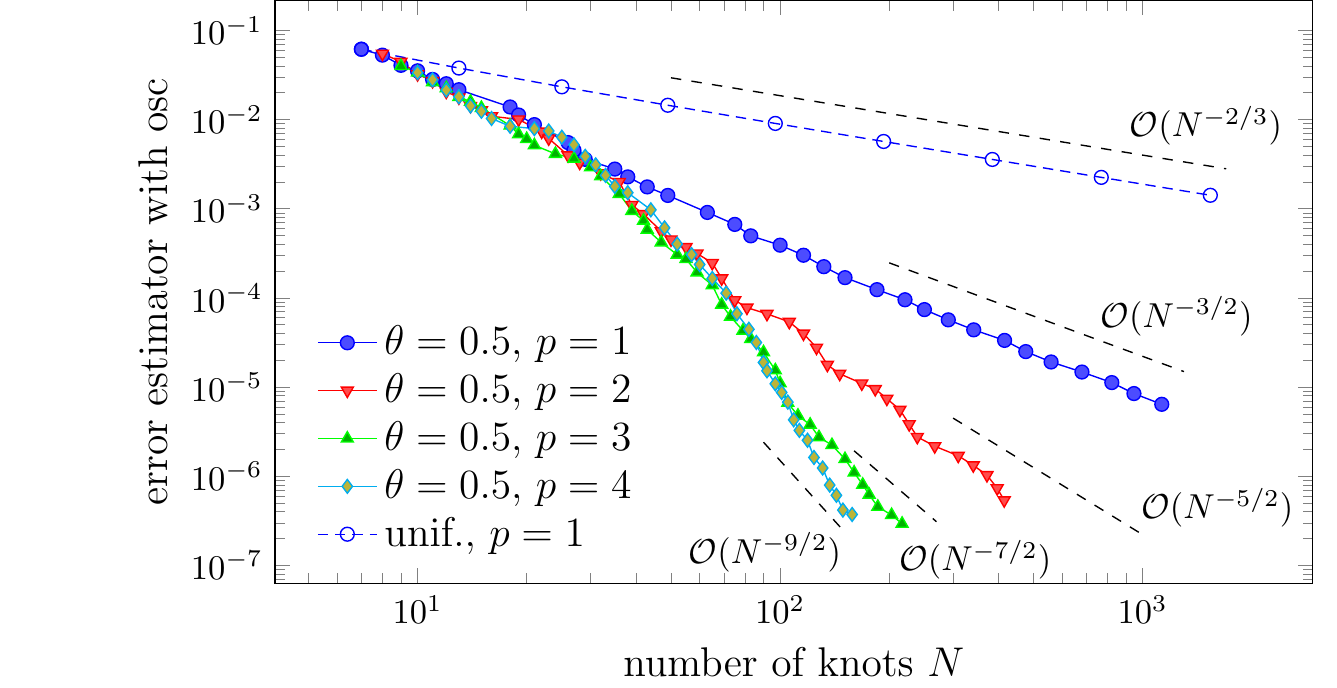}

	\caption{Example from Section~\ref{section:hyp_singular_heart}: (Approximate) error $\norm{u-U_\ell}{\mathfrak W}$ and error estimator $\eta_\ell$ with oscillations with respect to the number of knots $N$ for $\theta=0.5$, $p\in\{1,2,3,4\}$, different ansatz spaces (continuous piecewise polynomials, splines, NURBS), and different estimators ($(h-h/2)$, weighted-residual). Left: $(h-h/2)$-estimator. Right: Weighted-residual estimator. Top: Comparison of $p\in\{1,2,3,4\}$ for continuous piecewise polynomials. Middle: Comparison of $p\in\{1,2,3,4\}$ for splines. Bottom: Comparison of $p\in\{2,3,4\}$ for different ansatz spaces.}
	\label{fig:hyp_singular_comparison}
\end{figure}

In Figure~\ref{fig:hyp_singular_comparison}, we compare uniform refinement with $\theta=1$ and adaptive refinement with $\theta=0.5$, different ansatz spaces, and different estimators for Algorithm~\ref{the algorithm2}.
First (top), for $p=2$ fixed, we compare (an approximation of the) error for the different ansatz spaces and the algorithm steered either by the $(h-h/2)$-estimator (top, left) or the weighted-residual estimator (top, right). 
Instead of the exact error $\norm{u-U_\ell}{\mathfrak W}$, we compute $\norm{U_+-U_\ell}{\mathfrak W}$,  where  $\widehat\KK_\fine$ are the uniformly refined knots, which result from $\widehat \KK_\ell$ by adding the midpoint of each element $\widehat Q\in\widehat\QQ_\coarse$ with multiplicity one (see also Algorithm~\ref{alg:refinement}).
We get similar plots, where uniform mesh refinement leads to the suboptimal rate of convergence $\OO(N^{-2/3})$ because the solution lacks regularity. However, for all adaptive cases, Algorithm~\ref{the algorithm2} regains the optimal rate of convergence $\OO(N^{-1/2-p})$, 
where splines exhibit a better multiplicative constant than piecewise polynomials.
Next (middle), we consider the estimators for continuous piecewise polynomials. For both the $(h-h/2)$-estimator (middle, left) and the weighted-residual estimator (middle, right), we compare different orders $p\in\{1,2,3,4\}$ for the ansatz spaces. Again, uniform mesh refinement leads to the suboptimal rate $\OO(N^{-2/3})$ (only displayed for $p=1$), whereas the adaptive strategy regains the optimal order of convergence $\OO(N^{-1/2-p})$. 
Lastly, we get similar results for splines (bottom).


\section{Implementational aspects}
\label{sec:implementation}

\subsection{Computation of Galerkin matrices}
\label{sec:matrices}
The computation of the Galerkin matrices of~\eqref{eq:Galerkin_weak} and~\eqref{eq:Galerkin_hyper} is realized in {\tt VMatrix.c} and {\tt WMatrix.c}, and can be called in {\sc Matlab} via {\tt buildVMatrix} and {\tt buildWMatrix}. 
 Recall that $\edualW{R^{\rm c}_{\coarse,i',p}}{R^{\rm c}_{\coarse,i,p}}$ involves the term $\dual{R^{\rm c}_{\coarse,i',p}}{1}_{\partial\Omega}\dual{R^{\rm c}_{\coarse,i,p}}{1}_{\partial\Omega}$ if $\Gamma=\partial\Omega$, 
 which can  be approximated by standard tensor-Gauss quadrature.
 For the term $\dual{\W R^{\rm c}_{\coarse,i',p}}{R^{\rm c}_{\coarse,i,p}}_\Gamma$, we employ Maue's formula
\begin{align}\label{eq:maue}
\dual{\W u}{v}_\Gamma=\dual{\V\partial_\Gamma u}{\partial_\Gamma v}_\Gamma=\edualV{\partial_\Gamma u}{\partial_\Gamma v}\quad\text{for all }u,v\in \H^1(\Gamma);
\end{align}
see \cite{maue49}.
Note that $\partial_\Gamma R^{\rm c}_{\coarse,i',p}$ and $\partial_\Gamma R^{\rm c}_{\coarse,i,p}$ can be easily calculated via \eqref{eq:derivative of splines} and \eqref{eq:explicit derivative}.
Thus, it is sufficient to explain how to compute $\edualV{\phi}{\psi}$ for (at least) $\QQ_\coarse$-piecewise continuous functions $\phi,\psi\in \H^{-1/2}(\Gamma)$.

We start with
\begin{align*}
\edualV{\phi}{\psi}&=\int_\Gamma\int_\Gamma \phi(y) \psi(x) G(x,y) \d y\d x=\sum_{Q\in\QQ_\coarse}\sum_{Q'\in\QQ_\coarse} \int_Q\int_{Q'}  \phi(y) \psi(x) G(x,y) \d y\d x.
\end{align*}
Note that only elements $Q$ and $Q'$  that are contained in the supports of $\psi$ and $\phi$, respectively, contribute to the sum.
Recall from~\eqref{eq:locality} that (derivatives of) NURBS have  local support. 
Let $Q,Q'\in\QQ_\coarse$ with $Q=\gamma(\widehat Q)$, $\widehat Q=[\widehat x_{\coarse,j-1},\widehat x_{\coarse,j}]),$ $Q'=\gamma(\widehat Q')$, $\widehat Q'=[\widehat x_{\coarse,j'-1},\widehat x_{\coarse,j'}]$ for some $j,j'\in\{1,\dots,n_\coarse\}$.
Then, elementary integral transformations show that
\begin{align*}
 &\int_Q\int_{Q'}  \phi(y) \psi(x) G(x,y) \d y\d x=\int_{\widehat Q}\int_{\widehat Q'} \phi\big(\gamma(t)\big) \psi\big(\gamma(s)\big) G\big(\gamma(s),\gamma(t)\big)|\gamma'(t)| \,|\gamma'(s)| \d t \d s
 \\
 &=\int_0^1\int_0^1 \phi\big(\gamma_{j'}(\tau)\big)\psi\big(\gamma_j(\sigma)\big) G\big(\gamma_j(\sigma),\gamma_{j'}(\tau)\big) |\gamma_{j'}'(\tau)| \,|\gamma_j'(\sigma)| \d \tau \d \sigma,
\end{align*}
where $\gamma_j(\sigma):=\gamma\big(\widehat x_{\coarse,j-1}+\sigma(\widehat x_{\coarse,j}-\widehat x_{\coarse,j-1})\big)$ 
and $\gamma_{j'}(\tau):=\gamma\big(\widehat x_{\coarse,j'-1}+\tau(\widehat x_{\coarse,j'}-\widehat x_{\coarse,j'-1})\big)$.
We abbreviate $\widetilde \phi(\tau):=\phi\big(\gamma_{j'}(\tau)\big)|\gamma_{j'}'(\tau)|$ and $\widetilde \psi(\sigma):=\phi\big(\gamma_j(\sigma)\big)|\gamma_j'(\sigma)|$,
which gives 
\begin{align*}
\int_Q\int_{Q'}  \phi(y) \psi(x) G(x,y) \d y\d x=\int_0^1\int_0^1 \widetilde\phi(\tau)\widetilde\psi(\sigma) G\big(\gamma_j(\sigma),\gamma_{j'}(\tau)\big)  \d \tau \d \sigma.
\end{align*}
Now, we distinguish three cases.

\noindent\textbf{Case 1 (no intersection):} We assume that $Q\cap Q'=\emptyset$. 
In this case, the integrand is (at least) continuous and we can use standard tensor-Gauss quadrature to approximate the integral. 

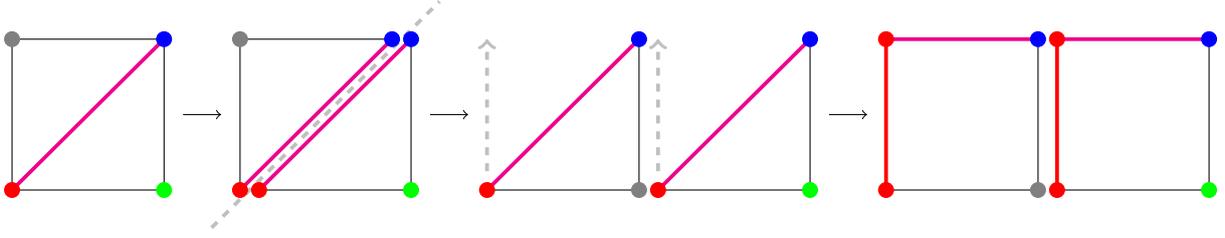
\begin{figure}
\begin{tikzpicture}
\tikzmath{\c=2;\s=0.25;\a=0.5;\e=0.25;}
\draw (0,0)
  -- (\c,0) 
  -- (\c,\c) 
  -- (0,\c)
  -- cycle;
\draw[line width = 0.5mm,magenta] (0,0) -- (\c,\c);
\fill[red] (0,0) circle (3pt);
\fill[green] (\c,0) circle (3pt);
\fill[blue] (\c,\c) circle (3pt);
\fill[gray] (0,\c) circle (3pt);
\draw[->]        (\c+\s,\c/2)   -- (\c+\s+\a,\c/2);
\draw (\c+2*\s+\a,0) 
  -- (\c+2*\s+\a,\c) 
  -- (2*\c+2*\s+\a,\c) 
  -- cycle;
\draw[line width = 0.5mm,magenta] (\c+2*\s+\a,0) -- (2*\c+2*\s+\a,\c);
\draw[dashed, line width=0.5mm, lightgray ] (\c+2.5*\s+\a-2*\e,-2*\e) -- (2*\c+2.5*\s+\a+2*\e,\c+2*\e);
\draw (\c+2*\s+\a+\e,0) 
  -- (2*\c+2*\s+\a+\e,0) 
  -- (2*\c+2*\s+\a+\e,\c) 
  -- cycle;
\draw[line width = 0.5mm,magenta] (\c+2*\s+\a+\e,0) -- (2*\c+2*\s+\a+\e,\c);
\fill[red] (\c+2*\s+\a,0) circle (3pt);
\fill[blue] (2*\c+2*\s+\a,\c) circle (3pt);
\fill[gray] (\c+2*\s+\a,\c) circle (3pt);
\fill[red] (\c+2*\s+\a+\e,0) circle (3pt);
\fill[green] (2*\c+2*\s+\a+\e,0) circle (3pt);
\fill[blue] (2*\c+2*\s+\a+\e,\c) circle (3pt);
\draw[->]        (2*\c+3*\s+\a+\e,\c/2)   -- (2*\c+3*\s+2*\a+\e,\c/2);
\draw (2*\c+4*\s+2*\a+\e,0) 
  -- (3*\c+4*\s+2*\a+\e,0) 
  -- (3*\c+4*\s+2*\a+\e,\c) 
  -- cycle;
\draw[line width = 0.5mm,magenta] (2*\c+4*\s+2*\a+\e,0) -- (3*\c+4*\s+2*\a+\e,\c);
\draw[->, dashed, line width=0.5mm, lightgray ]        (2*\c+4*\s+2*\a+\e,\e)   -- (2*\c+4*\s+2*\a+\e,\c);
\draw (3*\c+4*\s+2*\a+2*\e,0) 
  -- (4*\c+4*\s+2*\a+2*\e,0) 
  -- (4*\c+4*\s+2*\a+2*\e,\c) 
  -- cycle;
\draw[line width = 0.5mm,magenta] (3*\c+4*\s+2*\a+2*\e,0)  -- (4*\c+4*\s+2*\a+2*\e,\c);
\draw[->, dashed, line width=0.5mm, lightgray ]   (3*\c+4*\s+2*\a+2*\e,\e)   -- (3*\c+4*\s+2*\a+2*\e,\c);
\fill[red] (2*\c+4*\s+2*\a+\e,0) circle (3pt);
\fill[gray] (3*\c+4*\s+2*\a+\e,0) circle (3pt);
\fill[blue] (3*\c+4*\s+2*\a+\e,\c) circle (3pt);
\fill[red] (3*\c+4*\s+2*\a+2*\e,0) circle (3pt);
\fill[green] (4*\c+4*\s+2*\a+2*\e,0) circle (3pt);
\fill[blue] (4*\c+4*\s+2*\a+2*\e,\c) circle (3pt);
\draw[->]        (4*\c+5*\s+2*\a+2*\e,\c/2)   -- (4*\c+5*\s+3*\a+2*\e,\c/2);
\draw (4*\c+6*\s+3*\a+2*\e,0) 
  -- (5*\c+6*\s+3*\a+2*\e,0) 
  -- (5*\c+6*\s+3*\a+2*\e,\c) 
  -- (4*\c+6*\s+3*\a+2*\e,\c) 
  -- cycle;
\draw[magenta, line width=0.5mm] (5*\c+6*\s+3*\a+2*\e,\c) -- (4*\c+6*\s+3*\a+2*\e,\c);
\draw[red, line width=0.5mm] (4*\c+6*\s+3*\a+2*\e,0) -- (4*\c+6*\s+3*\a+2*\e,\c);
\draw (5*\c+6*\s+3*\a+3*\e,0) 
  -- (6*\c+6*\s+3*\a+3*\e,0) 
  -- (6*\c+6*\s+3*\a+3*\e,\c) 
  -- (5*\c+6*\s+3*\a+3*\e,\c) 
  -- cycle;
\draw[magenta, line width=0.5mm] (6*\c+6*\s+3*\a+3*\e,\c)  -- (5*\c+6*\s+3*\a+3*\e,\c);
\draw[red, line width=0.5mm] (5*\c+6*\s+3*\a+3*\e,0) -- (5*\c+6*\s+3*\a+3*\e,\c);
\fill[red] (4*\c+6*\s+3*\a+2*\e,0) circle (3pt);
\fill[gray] (5*\c+6*\s+3*\a+2*\e,0) circle (3pt);
\fill[blue] (5*\c+6*\s+3*\a+2*\e,\c) circle (3pt);
\fill[red] (4*\c+6*\s+3*\a+2*\e,\c) circle (3pt);
\fill[red] (5*\c+6*\s+3*\a+3*\e,0) circle (3pt);
\fill[green] (6*\c+6*\s+3*\a+3*\e,0) circle (3pt);
\fill[blue] (6*\c+6*\s+3*\a+3*\e,\c) circle (3pt);
\fill[red] (5*\c+6*\s+3*\a+3*\e,\c) circle (3pt);
\end{tikzpicture}
\caption{Visualization of the integral transformations from Case~2 and~3 in Section~\ref{sec:matrices}: the reflection $(\sigma,\tau)\mapsto(\tau,\sigma)$ for the upper triangle and the (inverse) Duffy transformation $(\sigma,\tau)\mapsto(\sigma,\tau/\sigma)$ for both triangles. 
The colors indicate the sets onto which each of the colored points and lines are mapped. 
}
\label{fig:duffy}
\end{figure}

\noindent\textbf{Case 2 (identical elements):} We assume that $Q=Q'$, which also implies that $j=j'$. We split the integral into two summands, use Fubini's theorem as well as the reflection $(\sigma,\tau)\mapsto(\tau,\sigma)$ for the second one, and use the Duffy transformation 
$(\sigma, \tau)\mapsto (\sigma,\sigma\tau)$ with Jacobi determinant $\sigma$ (see Figure~\ref{fig:duffy} for a visualization) 
\begin{align*}
&\int_0^1\int_0^1 \widetilde\phi(\tau)\widetilde\psi(\sigma) G\big(\gamma_j(\sigma),\gamma_{j'}(\tau)\big)  \d \tau \d \sigma
\\
&=\int_0^1\int_0^\sigma \widetilde\phi(\tau)\widetilde\psi(\sigma) G\big(\gamma_j(\sigma),\gamma_j(\tau)\big)  \d \tau \d \sigma
+\int_0^1\int_\sigma^1 \widetilde\phi(\tau)\widetilde\psi(\sigma) G\big(\gamma_j(\sigma),\gamma_j(\tau)\big)  \d \tau \d \sigma
\\
&=\int_0^1\int_0^\sigma \widetilde\phi(\tau)\widetilde\psi(\sigma) G\big(\gamma_j(\sigma),\gamma_j(\tau)\big)  \d \tau \d \sigma
+\int_0^1\int_0^\sigma \widetilde\phi(\sigma)\widetilde\psi(\tau) G\big(\gamma_j(\tau),\gamma_j(\sigma)\big)  \d \tau \d \sigma
\\
&=\int_0^1\int_0^1 \Big(\widetilde\phi(\sigma\tau)\widetilde\psi(\sigma) G\big(\gamma_j(\sigma),\gamma_j(\sigma\tau)\big)
+\widetilde\phi(\sigma)\widetilde\psi(\sigma\tau) G\big(\gamma_j(\sigma\tau),\gamma_j(\sigma)\big)\Big) \sigma \d \tau \d \sigma.
\intertext{Recall that $G(x,y)=-\frac{1}{2\pi} \log|x-y|$.
Thus, we further obtain (with the transformation $(\sigma,\tau)\mapsto(\sigma,1-\tau)$ for the last integral) that}
&=-\frac{1}{2\pi}\int_0^1\int_0^1 \big(\widetilde\phi(\sigma\tau)\widetilde\psi(\sigma) + \widetilde\phi(\sigma)\widetilde\psi(\sigma\tau)\big)
\log\left(\frac{\big|\gamma_j(\sigma)-\gamma_j(\sigma\tau)\big|}{\sigma-\sigma\tau} (\sigma-\sigma\tau)\right) \sigma \d \tau \d \sigma
\\
&=-\frac{1}{2\pi}\int_0^1\int_0^1 \big(\widetilde\phi(\sigma\tau)\widetilde\psi(\sigma)+\widetilde\phi(\sigma)\widetilde\psi(\sigma\tau)\big)
 \log\left(\frac{\big|\gamma_j(\sigma)-\gamma_j(\sigma\tau)\big|}{\sigma-\sigma\tau}\right) \sigma \d \tau \d \sigma
\\
&\qquad-\frac{1}{2\pi}\int_0^1\int_0^1 \big(\widetilde\phi(\sigma\tau)\widetilde\psi(\sigma)+\widetilde\phi(\sigma)\widetilde\psi(\sigma\tau)\big) \log(\sigma) \sigma \d \tau \d \sigma
\\
&\qquad-\frac{1}{2\pi}\int_0^1\int_0^1 \big(\widetilde\phi(\sigma(1-\tau))\widetilde\psi(\sigma)+\widetilde\phi(\sigma)\widetilde\psi(\sigma(1-\tau))\big) \log(\tau) \sigma \d \tau \d \sigma. 
\end{align*}
With the fundamental theorem of calculus and the fact that $\gamma$ is piecewise smooth with $|\gamma'|>0$, it is easy to see that $(\sigma,\tau)\mapsto{\big|\gamma_j(\sigma)-\gamma_j(\sigma\tau)\big|}/(\sigma-\sigma\tau)$ is smooth and (uniformly) larger than $0$; see \cite[Lemma~5.2]{gantner14} for details.
Thus, one can use standard tensor-Gauss quadrature to approximate the first integral.
Note that the term  $\log(\sigma)\sigma$ is continuous, but not even $C^1$.
Hence, we use tensor-Gauss quadrature with weight function $\log(\sigma)$ and $\log(\tau)$ for the second and third integral, respectively. 

\noindent\textbf{Case 3 (common node):} We assume that $Q\cap Q'$ contains only one point.
Without loss of generality, we further assume that the singularity is at $(\sigma,\tau)=(0,1)$, i.e., $\widehat x_{\coarse,j-1}=\widehat x_{\coarse,j'}$ or $\widehat x_{\coarse,j-1}=a\wedge \widehat  x_{\coarse,j'}=b$ if $\Gamma=\partial\Omega$.
We rotate the integration domain by $\pi/2$, i.e., $(\sigma,\tau)\mapsto(\tau,1-\sigma)$, which transforms the singularity to $(\sigma,\tau)=(0,0)$, and then employ the same transformations as in Case~2 (see also Figure~\ref{fig:duffy})
\begin{align*}
\int_0^1\int_0^1 \widetilde\phi(\tau)\widetilde\psi(\sigma) G\big(\gamma_j(\sigma),\gamma_{j'}(\tau)\big)  \d \tau \d \sigma
&=\int_0^1\int_0^1 \widetilde\phi(1-\sigma)\widetilde\psi(\tau) G\big(\gamma_j(\tau),\gamma_{j'}(1-\sigma)\big)  \d \tau \d \sigma
\\
&=\int_0^1\int_0^1 \Big(\widetilde\phi(1-\sigma)\widetilde\psi(\sigma\tau) G\big(\gamma_j(\sigma\tau),\gamma_{j'}(1-\sigma)\big)
\\
&\quad+\widetilde\phi(1-\sigma\tau)\widetilde\psi(\sigma) G\big(\gamma_j(\sigma),\gamma_{j'}(1-\sigma\tau)\big)\Big) \sigma \d \tau \d \sigma.
\end{align*}
Recall that $G(x,y)=-\frac{1}{2\pi} \log|x-y|$.
Thus, we further obtain that
\begin{align*}
&=-\frac{1}{2\pi}\int_0^1\int_0^1 \widetilde\phi(1-\sigma)\widetilde\psi(\sigma\tau) \log\left(\frac{\big|\gamma_j(\sigma\tau)-\gamma_{j'}(1-\sigma)\big|}{\sigma}\,\sigma\right)\sigma \d \tau \d \sigma
\\
&\qquad-\frac{1}{2\pi}\int_0^1\int_0^1 \widetilde\phi(1-\sigma\tau)\widetilde\psi(\sigma) \log\left(\frac{\big|\gamma_j(\sigma)-\gamma_{j'}(1-\sigma\tau)\big|}{\sigma}\,\sigma\right)\sigma \d \tau \d \sigma
\\
&=-\frac{1}{2\pi}\int_0^1\int_0^1 \widetilde\phi(1-\sigma)\widetilde\psi(\sigma\tau) \log\left(\frac{\big|\gamma_j(\sigma\tau)-\gamma_{j'}(1-\sigma)\big|}{\sigma}\right)\sigma \d \tau \d \sigma
\\
&\qquad-\frac{1}{2\pi}\int_0^1\int_0^1 \widetilde\phi(1-\sigma\tau)\widetilde\psi(\sigma) \log\left(\frac{\big|\gamma_j(\sigma)-\gamma_{j'}(1-\sigma\tau)\big|}{\sigma}\right)\sigma \d \tau \d \sigma
\\
&\qquad-\frac{1}{2\pi}\int_0^1\int_0^1 \big(\widetilde\phi(1-\sigma)\widetilde\psi(\sigma\tau)+ \widetilde\phi(1-\sigma\tau)\widetilde\psi(\sigma)\big) \log(\sigma) \sigma \d \tau \d \sigma
.
\end{align*}
With the fundamental theorem of calculus and the fact that $\gamma$ is piecewise smooth with $|\gamma'|>0$, it is easy to see that $(\sigma,\tau)\mapsto{\big|\gamma_j(\sigma\tau)-\gamma_{j'}(1-\sigma)\big|}/\sigma$ 
and $(\sigma,\tau)\mapsto{\big|\gamma_j(\sigma)-\gamma_{j'}(1-\sigma\tau)\big|}/\sigma$ are smooth and (uniformly) larger than $0$; see \cite[Lemma~5.3]{gantner14} for details.
Thus, one can use standard tensor-Gauss quadrature to approximate the first and second integral.
Note that the term  $\log(\sigma)\sigma$ is continuous, but not even $C^1$.
Hence, we use tensor-Gauss quadrature with weight function $\log(\sigma)$ for the third integral. 

%

\subsection{Computation of right-hand side vectors}
\label{sec:right vectors}
The computation of the right-hand side vectors of \eqref{eq:Galerkin_weak} and \eqref{eq:Galerkin_hyper}  is realized in {\tt RHSVectorWeak.c} and {\tt RHSVectorHyp.c}, where the required projection $\Pi_\coarse$ is realized in {\tt PhiApprox.c}.
They can be called in {\sc Matlab} via {\tt buildRHSVectorWeak}, {\tt buildRHSVectorHyp}, and {\tt buildPhiApprox}.
Clearly, the terms $\dual{u/2}{R_{\coarse,i,p}}_\Gamma$ and $\dual{(\Pi_\coarse\phi)/2}{R^{\rm c}_{\coarse,i,p}}_\Gamma$ can  be approximated by standard tensor-Gauss quadrature.
 For the term $\dual{\K' (\Pi_\coarse\phi)}{R^{\rm c}_{\coarse,i,p}}_\Gamma$, we use the fact that $\K'$ is the adjoint operator of $\K$
\begin{align}
\dual{\K'\Pi_\coarse\phi}{R^{\rm c}_{\coarse,i,p}}_\Gamma=\dual{\K R^{\rm c}_{\coarse,i,p}}{\Pi_\coarse\phi}_\Gamma.
\end{align}
Thus, it is sufficient to explain how to compute $\dual{\K u }{\psi}$ for some (at least) $\QQ_\coarse$-piecewise continuous functions $u\in \H^{1/2}(\Gamma)$, $\psi\in \H^{-1/2}(\Gamma)$.
For $x,y\in\Gamma$, we abbreviate 
\begin{align}\label{eq:G nu}
G_{\nu}(x,y):=  \frac{\partial_y}{\partial\nu(y)} G(x,y)=\frac{1}{2\pi}\frac{(x-y)\cdot\nu(y)}{|x-y|^2}.
\end{align} 
We start with
\begin{align*}
\dual{\K u}{\psi}&=\int_\Gamma\int_\Gamma u(y) \psi(x) G_\nu(x,y) \d y\d x=\sum_{Q\in\QQ_\coarse}\sum_{Q'\in\QQ_\coarse} \int_Q\int_{Q'}  u(y) \psi(x) G_\nu(x,y) \d y\d x.
\end{align*}
Note that only elements $Q$ and $Q'$  that are contained in the supports of $\psi$ and $u$, respectively, contribute to the sum.
Recall from~\eqref{eq:locality} that NURBS have local support. 
Let $Q,Q'\in\QQ_\coarse$ with $Q=\gamma(\widehat Q)$, $\widehat Q=[\widehat x_{\coarse,j-1},\widehat x_{\coarse,j}]),$ $Q'=\gamma(\widehat Q')$, $\widehat Q'=[\widehat x_{\coarse,j'-1},\widehat x_{\coarse,j'}])$ for some $j,j'\in\{1,\dots,n_\coarse\}$.
As in Section~\ref{sec:matrices}, elementary integral transformations show that
\begin{align*}
 &\int_Q\int_{Q'}  u(y) \psi(x) G_\nu(x,y) \d y\d x=\int_0^1\int_0^1 \widetilde u(\tau)\widetilde\psi(\sigma) G_\nu\big(\gamma_j(\sigma),\gamma_{j'}(\tau)\big)  \d \tau \d \sigma,
\end{align*}
where $\gamma_j(\sigma):=\gamma\big(\widehat x_{\coarse,j-1}+\sigma(\widehat x_{\coarse,j}-\widehat x_{\coarse,j-1})\big)$ 
and $\gamma_{j'}(\tau):=\gamma\big(\widehat x_{\coarse,j'-1}+\tau(\widehat x_{\coarse,j'}-\widehat x_{\coarse,j'-1})\big)$, 
and $\widetilde u(\tau):=u\big(\gamma_{j'}(\tau)\big)|\gamma_{j'}'(\tau)|$ and $\widetilde \psi(\sigma):=u\big(\gamma_j(\sigma)\big)|\gamma_j'(\sigma)|$.
Now, we distinguish three cases.

\noindent\textbf{Case 1 (no intersection):} We assume that $Q\cap Q'=\emptyset$. 
In this case, the integrand is (at least) continuous and we can use standard tensor-Gauss quadrature to approximate the integral. 

\noindent\textbf{Case 2 (identical elements):} We assume that $Q=Q'$, which also implies that $j=j'$. 
Note that Taylor expansion together with the fact that $\gamma$ is piecewise $C^2$ shows that 
\begin{align}\label{eq:normal product}
(x-y)\cdot\nu(y)=\OO(|x-y|^2) \quad\text{for all } x,y \in Q; 
\end{align}
see, e.g., \cite[Equation~(5.4)]{gantner14} for a detailed calculation.
Since $\gamma$ is even piecewise $C^\infty$, this argument even yields that $(\sigma,\tau)\mapsto G_\nu\big(\gamma_j(\sigma),\gamma_{j}(\tau)\big)$ is smooth.
Thus, in principle, one can use standard tensor-Gauss quadrature to approximate the integral.  
However, to avoid the computation of $G_\nu(\gamma_j(\sigma),\gamma_{j}(\tau))$ in the limit case $\sigma=\tau$, 
we make the same steps as in Case~2 of Section~\ref{sec:matrices}, i.e., we split the integral into two summands, use Fubini's theorem as well as the reflection $(\sigma,\tau)\mapsto(\tau,\sigma)$ for the second one, and use the Duffy transformation 
$(\sigma, \tau)\mapsto (\sigma,\sigma\tau)$ with Jacobi determinant $\sigma$ (see Figure~\ref{fig:duffy} for a visualization).
 Altogether, this results in
\begin{align*}
&\int_0^1\int_0^1 \widetilde u(\tau)\widetilde\psi(\sigma) G_\nu\big(\gamma_j(\sigma),\gamma_{j'}(\tau)\big)  \d \tau \d \sigma
\\
&=\int_0^1\int_0^1 \Big(\widetilde u(\sigma\tau)\widetilde\psi(\sigma) G_\nu\big(\gamma_j(\sigma),\gamma_j(\sigma\tau)\big)
+\widetilde u(\sigma)\widetilde\psi(\sigma\tau) G_\nu\big(\gamma_j(\sigma\tau),\gamma_j(\sigma)\big)\Big) \sigma \d \tau \d \sigma.
\end{align*}
In particular, this shifts the limit case to the boundary $\sigma=0\vee\tau=1$, while tensor-Gauss quadrature only evaluates the integral in the interior of the unit square.

\begin{remark}
We mention that the smoothness of $(\sigma,\tau)\mapsto G_\nu(\gamma_j(\sigma),\gamma_{j}(\tau))$ hinges on the considered Laplace equation.
It is also satisfied for the Helmholtz equation, but fails, for instance, for the Lam\'e equation from linear elasticity. 
Nevertheless, the same transformations as above yield a final integrand that is (at least) continuous if the Lam\'e equation is considered.
\end{remark}

\noindent\textbf{Case 3 (common node):} We assume that $Q\cap Q'$ contains only one point.
Without loss of generality, we further assume that the singularity is at $(\sigma,\tau)=(0,1)$, i.e., $\widehat x_{\coarse,j-1}=\widehat x_{\coarse,j'}$ or $\widehat x_{\coarse,j-1}=a\wedge \widehat  x_{\coarse,j'}=b$ if $\Gamma=\partial\Omega$.
We make the same steps as in Case~3 of Section~\ref{sec:matrices}, i.e., we rotate the integration domain by $\pi/2$, i.e., $(\sigma,\tau)\mapsto(\tau,1-\sigma)$, which transforms the singularity to $(\sigma,\tau)=(0,0)$, and then employ the same transformations as in Case~2 (see also Figure~\ref{fig:duffy}).
Altogether this results in 
\begin{align*}
\int_0^1\int_0^1 \widetilde u(\tau)\widetilde\psi(\sigma) G_\nu\big(\gamma_j(\sigma),\gamma_{j'}(\tau)\big)  \d \tau \d \sigma
=\int_0^1\int_0^1 \Big(\widetilde u(1-\sigma)\widetilde\psi(\sigma\tau) G_\nu\big(\gamma_j(\sigma\tau),\gamma_{j'}(1-\sigma)\big)
\\
+\widetilde u(1-\sigma\tau)\widetilde\psi(\sigma) G_\nu\big(\gamma_j(\sigma),\gamma_{j'}(1-\sigma\tau)\big)\Big) \sigma \d \tau \d \sigma.
\end{align*}
Note that $G_\nu\big(\gamma_j(\sigma\tau),\gamma_{j'}(1-\sigma)\big)\sigma=G_\nu\big(\gamma_j(\sigma),\gamma_{j'}(1-\sigma\tau)\big)\sigma=\OO(1)$.
In particular, one can show with the fundamental theorem of calculus and the fact that $\gamma$ is piecewise $C^1$ with $|\gamma'|>0$ that these terms are (at least) continuous.
Thus, one can use standard tensor-Gauss quadrature to approximate the integral.

\subsection{Evaluation of single-layer operator $\V$}\label{sec:eval V}

Let $\widehat \KK_\coarse$ be a knot vector.
Let $x\in\Gamma$  and $s\in[a,b]$ with $x=\gamma(s)$.
Moreover, let $Q\in\QQ_\coarse$ with $x\in Q$ and $s\in\widehat Q:=\gamma^{-1}(Q)=[\widehat x_{\coarse,j-1},\widehat x_{\coarse,j}]$.
For $\phi\in \H^{-1/2}(\Gamma)$ (at least) $\QQ_\coarse$-piecewise continuous, 
we want to evaluate $[\V\phi](x)$.
These values are required to compute the Faermann estimator in Section~\ref{sec:compute Faermann} and the weighted-residual estimator in Section~\ref{sec:compute weighted}.
We also mention that they can be used to implement a collocation method instead of a Galerkin method to approximate the solution of \eqref{eq:weak strong}. 
The evaluation is realized in {\tt ResidualWeak.c} and can be called in {\sc Matlab} via {\tt evalV}. 

We start with 
\begin{align*}
[\V\phi](x)=\int_\Gamma \phi(y)G(x,y)\d y
= \sum_{\widehat Q'\in\widehat\QQ_\coarse} \int_{\widehat Q'} \phi\big(\gamma(t)\big) G\big(\gamma(s),\gamma(t)\big)\big|\gamma'(t)\big| \d t.
\end{align*}
For all $\widehat Q'\in\widehat\QQ_\coarse$ with $\widehat Q'\neq \widehat Q$, the integrands are (at least) continuous and we can use standard tensor-Gauss quadrature 
after transforming the integration domain to the unit square.
Recall that $G(x,y)=-\frac{1}{2\pi}\log|x-y|$.
With the abbreviation $\overline \phi(t):=\phi\big(\gamma(t)\big)\big|\gamma'(t)\big|$, it remains to consider 
\begin{align*}
&\int_{\widehat Q'} \overline\phi(t) G\big(\gamma(s),\gamma(t)\big)\d t
=\int_{\widehat x_{\coarse,j-1}}^s \overline\phi(t) G\big(\gamma(s),\gamma(t)\big) \d t
+\int_s^{\widehat x_{\coarse,j}} \overline\phi(t) G\big(\gamma(s),\gamma(t)\big) \d t
\\
&=(s-\widehat x_{\coarse,j-1})\int_0^1 \overline\phi\big(\widehat x_{\coarse,j-1}+\tau(s-\widehat x_{\coarse,j-1})\big) G\Big(\gamma(s),\gamma\big(\widehat x_{\coarse,j-1}+\tau(s-\widehat x_{\coarse,j-1})\big)\Big) \d \tau
\\
&\qquad+(\widehat x_{\coarse,j}-s)\int_0^1 \overline\phi\big(s+\tau(\widehat x_{\coarse,j}-s)\big) G\Big(\gamma(s),\gamma(s+\tau\big(\widehat x_{\coarse,j}-s)\big)\Big) \d \tau
\\
&=-\frac{1}{2\pi}(s-\widehat x_{\coarse,j-1})\int_0^1 \overline\phi\big(\widehat x_{\coarse,j-1}+\tau(s-\widehat x_{\coarse,j-1})\big) \log\left(\frac{\Big|\gamma(s)-\gamma\big(\widehat x_{\coarse,j-1}+\tau(s-\widehat x_{\coarse,j-1})\big)\Big|}{\tau}\right) \d \tau
\\
&\qquad-\frac{1}{2\pi}(s-\widehat x_{\coarse,j-1})\int_0^1 \overline\phi\big(\widehat x_{\coarse,j-1}+\tau(s-\widehat x_{\coarse,j-1})\big) \log(\tau) \d \tau
\\
&\qquad-\frac{1}{2\pi}(\widehat x_{\coarse,j}-s)\int_0^1 \overline\phi\big(s+\tau(\widehat x_{\coarse,j}-s)\big) \log\left(\frac{\Big|\gamma(s)-\gamma\big(s+\tau(\widehat x_{\coarse,j}-s)\big)\Big|}{\tau}\right) \d \tau
\\
&\qquad-\frac{1}{2\pi}(\widehat x_{\coarse,j}-s)\int_0^1 \overline\phi\big(s+\tau(\widehat x_{\coarse,j}-s)\big) \log(\tau) \d \tau.
\end{align*}
With the fundamental theorem of calculus and the fact that $\gamma$ is piecewise smooth with $|\gamma'|>0$, it is easy to see that $\tau\mapsto\big|\gamma(s)-\gamma\big(\widehat x_{\coarse,j-1}+\tau(s-\widehat x_{\coarse,j-1})\big)\big|/\tau$ and $\tau\mapsto |\gamma(s)-\gamma\big(s+\tau(\widehat x_{\coarse,j}-s)\big)\big|/\tau$ are smooth and (uniformly) larger than $0$;  see \cite[Lemma~5.6]{gantner14} for details.
Thus, we can use Gauss quadrature with weight function $1$ and $\log(\tau)$ to approximate the final integrals.

\subsection{Evaluation of double-layer operator $\K$}\label{sec:eval K}

Let $\widehat \KK_\coarse$ be a knot vector.
Let $x\in\Gamma$  and $s\in[a,b]$ with $x=\gamma(s)$.
Moreover, let $Q\in\QQ_\coarse$ with $x\in Q$ and $s\in\widehat Q:=\gamma^{-1}(Q)=[\widehat x_{\coarse,j-1},\widehat x_{\coarse,j}]$.
For $u\in \H^{1/2}(\Gamma)$ (at least) $\QQ_\coarse$-piecewise $C^1$, 
we want to evaluate $[\K u](x)$.
These values are required to compute the Faermann estimator in Section~\ref{sec:compute Faermann} and the weighted-residual estimator in Section~\ref{sec:compute weighted}.
We also mention that they can be used to implement a collocation method\footnote{Attention has to be paid for the direct method. 
Indeed, the representation $f(x)=[(1/2+\mathfrak{K})u](x)$ for the right-hand side is actually only valid for $x\in\Gamma\setminus\NN_\gamma$.
For the general representation (depending additionally on the angle at corners $x\in\NN_\gamma$), we refer to, e.g., \cite[Lemma~6.8 and 6.11]{steinbach08}.} 
instead of a Galerkin method to approximate the solution of \eqref{eq:weak strong}. 
The evaluation is realized in {\tt ResidualWeak.c} and can be called in {\sc Matlab} via {\tt evalRHSWeak}.

Recall the abbreviation $G_\nu(x,y)$ of \eqref{eq:G nu}.
We start with 
\begin{align*}
[\K u](x)=\int_\Gamma \phi(y)G_\nu(x,y)\d y
= \sum_{\widehat Q'\in\widehat\QQ_\coarse} \int_{\widehat Q'} u\big(\gamma(t)\big) G_\nu\big(\gamma(s),\gamma(t)\big)\big|\gamma'(t)\big| \d t.
\end{align*}
For all $\widehat Q'\in\widehat\QQ_\coarse$ with $\widehat Q'\neq \widehat Q$, the integrands are (at least) continuous and we can use standard tensor-Gauss quadrature 
after transforming the integration domain to the unit square.
Recall that $G(x,y)=-\frac{1}{2\pi}\log|x-y|$.
With the abbreviation $\overline u(t):=u\big(\gamma(t)\big)\big|\gamma'(t)\big|$, it remains to consider 
\begin{align*}
&\int_{\widehat Q'} \overline u(t) G_\nu\big(\gamma(s),\gamma(t)\big)\d t
=\int_{\widehat x_{\coarse,j-1}}^s \overline u(t) G_\nu\big(\gamma(s),\gamma(t)\big) \d t
+\int_s^{\widehat x_{\coarse,j}} \overline u(t) G_\nu\big(\gamma(s),\gamma(t)\big) \d t
\\
&=(s-\widehat x_{\coarse,j-1})\int_0^1 \overline u\big(\widehat x_{\coarse,j-1}+\tau(s-\widehat x_{\coarse,j-1})\big) G_\nu\Big(\gamma(s),\gamma\big(\widehat x_{\coarse,j-1}+\tau(s-\widehat x_{\coarse,j-1})\big)\Big) \d \tau
\\
&\qquad+(\widehat x_{\coarse,j}-s)\int_0^1 \overline u\big(s+\tau(\widehat x_{\coarse,j}-s)\big) G_\nu\Big(\gamma(s),\gamma(s+\tau\big(\widehat x_{\coarse,j}-s)\big)\Big) \d \tau.
\end{align*}
With the fundamental theorem of calculus and the fact that $\gamma$ is piecewise smooth with $|\gamma'|>0$, it is easy to see that $\tau\mapsto\big|\gamma(s)-\gamma\big(\widehat x_{\coarse,j-1}+\tau(s-\widehat x_{\coarse,j-1})\big)\big|/\tau$ and $\tau\mapsto |\gamma(s)-\gamma\big(s+\tau(\widehat x_{\coarse,j}-s)\big)\big|/\tau$ are smooth and (uniformly) larger than $0$;  see \cite[Lemma~5.6]{gantner14} for details.
Based on \eqref{eq:normal product}, one can show that the final two integrands are (at least) continuous; see \cite[Lemma~5.7]{gantner14} for details.
Thus, we can use standard Gauss quadrature to approximate the integrals.

\subsection{Evaluation of hyper-singular integral operator $\W$}\label{sec:eval W}

Under the assumptions of Section~\ref{sec:eval K}, we want to evaluate $[\W u](x)$ for $x\in\Gamma\setminus\NN_\gamma$.
These values are required to compute the weighted-residual estimators in Section~\ref{sec:compute weighted}.
We also mention that they can be used to implement a collocation method instead of a Galerkin method to approximate the solution of \eqref{eq:hyper strong}. 
The evaluation is realized in {\tt ResidualHyp.c} and can be called in {\sc Matlab} via {\tt evalW}.
Maue's formula \eqref{eq:maue}, $\V\partial u\in \H^1(\Gamma)$ (see Section~\ref{sec:weaksing}), and integration by parts show that $\dual{\W u}{v}_\Gamma=-\dual{\partial_\Gamma\V\partial_\Gamma u}{v}_\Gamma$
for all $v\in \H^1(\Gamma)$.
By density, we see that 
\begin{align}
\W u=-\partial_\Gamma \V \partial_\Gamma u.
\end{align}
To approximate $(\W u)(x)$, we replace $(\V\partial_\Gamma u)\circ\gamma|_{\widehat Q}$ by an interpolation polynomial, where the required point evaluations have been discussed in Section~\ref{sec:eval V}.
The arclength derivatives can be computed with \eqref{eq:derivative of splines} and  \eqref{eq:explicit derivative}.

\subsection{Evaluation of adjoint double-layer operator $\K'$}\label{sec:eval K'}

Under the assumptions of Section~\ref{sec:eval V}, we want to evaluate $[\K u](x)$ for $x\in\Gamma\setminus\NN_\gamma$.
These values are required to compute the weighted-residual estimators in Section~\ref{sec:compute weighted}.
We also mention that they can be used to implement a collocation method instead of a Galerkin method to approximate the solution of \eqref{eq:hyper strong}. 
The evaluation is realized in {\tt ResidualHyp.c} and can be called in {\sc Matlab} via {\tt evalRHSHyp}.
To approximate $[\K'\phi](x)$, one can proceed along the lines of Section~\ref{sec:eval K}; see \cite[Section~6.4]{schimanko16} for details.

\subsection{Computation of \emph{h--h/2} error estimators}

Let $\widehat \KK_\coarse$ be a knot vector with uniformly refined knot vector $\widehat \KK_\fine$ (see Section~\ref{sec:a_posteriori_weak}). 
The computation of the error indicators $\eta_{\V,\rm hh2,\coarse}(x)$ and $\eta_{\W,\rm hh2,\coarse}(x)$, $x\in\NN_\coarse$, of \eqref{eq:hhV} and \eqref{eq:hhW}, is realized in {\tt HHEstWeak.c} and {\tt HHEstHyp.c}, and can be called in {\sc Matlab} via {\tt buildHHEstWeak} and {\tt buildHHEstHyp} (which expect among others the coefficients of the Galerkin approximations on the coarse and on the fine mesh.)
Clearly, $\eta_{\V,\rm hh2,\coarse}(x)$ can be approximated by standard Gauss quadrature. 
With \eqref{eq:derivative of splines} and \eqref{eq:explicit derivative}, also the indicators $\eta_{\W,\rm hh2,\coarse}(x)$ can be approximated by standard Gauss quadrature.

\subsection{Computation of Faermann estimator}\label{sec:compute Faermann}

Let $\widehat \KK_\coarse$ be a knot vector. 
The computation of the error indicators $\eta_{\V,\rm fae, \coarse}(x)$, $x\in\NN_\coarse$, of \eqref{eq:fae} is realized in {\tt FaerEstWeak.c}, and can be called in {\sc Matlab} via {\tt buildFaerEstWeak}.
Let $Q, Q'\in\QQ_\coarse$ with $Q\cup Q'=\omega_\coarse(x)$.
We assume that $Q\neq Q'$, i.e., $x$ is not at the boundary of $\Gamma\subsetneqq \partial\Omega$.
The case $Q=Q'$ works analogously.
We abbreviate the residual $r_\coarse:=f-\V\Phi_\coarse$, which gives that
\begin{align*}
\eta_{\V,\rm fae, \coarse}(x)^2=\seminorm{r_\coarse}{H^{1/2}(\omega_\coarse(x))}^2=\seminorm{r_\coarse}{H^{1/2}(Q)}^2+2\int_Q\int_{Q'} \frac{\big|r_\coarse(y)-r_\coarse(z)\big|^2}{|y-z|^2} \d z \d y+\seminorm{r_\coarse}{H^{1/2}(Q')}^2.
\end{align*}
Given the evaluation procedures of Section~\ref{sec:eval V} and \ref{sec:eval K}, the terms $\seminorm{r_\coarse}{H^{1/2}(Q)}^2$ and $\seminorm{r_\coarse}{H^{1/2}(Q')}^2$ can be approximated exactly as in Case~2 of Section~\ref{sec:right vectors}, where one can additionally exploit the symmetry of the integrand.
The remaining integral can be approximated exactly as in  in Case~3 of Section~\ref{sec:right vectors}.

\subsection{Computation of weighted-residual error estimators}\label{sec:compute weighted}

Let $\widehat \KK_\coarse$ be a knot vector. 
The computation of the error indicators $\eta_{\V,\rm res,\coarse}(x)$ and $\eta_{\W,\rm res,\coarse}(x)$, $x\in\NN_\coarse$, of \eqref{eq:resV} and \eqref{eq:resW}, is realized in {\tt ResEstWeak.c} and {\tt ResEstHyp.c},  where the required projection $\Pi_\coarse$ is realized in {\tt PhiApprox.c}.
They can be called in {\sc Matlab} via {\tt buildResEstWeak}, {\tt buildResEstHyp}, and {\tt buildPhiApprox.c}.
To approximate the error indicators $\eta_{\V,\rm res, \coarse}(x)$, $x\in\NN_\coarse$, of \eqref{eq:resV}, we replace $(f-\V\Phi_\coarse)\circ\gamma|_{\widehat Q}$, 
$\widehat Q\in\widehat\QQ_\coarse$ with $\gamma(\widehat Q)\subseteq\omega_\coarse(x)$, by an interpolation polynomial, where the required point evaluations have been discussed in Section~\ref{sec:eval V}--\ref{sec:eval K}.
The arclength derivatives can be computed with \eqref{eq:derivative of splines} and \eqref{eq:explicit derivative}.
Finally, we use standard Gauss quadrature.
To approximate the error indicators $\eta_{\W,\rm res, \coarse}(x)$ $x\in\NN_\coarse$, of \eqref{eq:resW}, we can use standard Gauss quadrature. 
The required point evaluations have been discussed in Section~\ref{sec:eval W}--\ref{sec:eval K'}.

\subsection{Computation of data oscillations}

Let $\widehat \KK_\coarse$ be a knot vector. 
Clearly, the oscillations $\osc_{\W,\coarse}(x)$, $x\in\NN_\coarse$, of \eqref{eq:osc} can be approximated by standard Gauss quadrature. 

\appendix


\section{Overview of functions provided by \texttt{IGABEM2D}}\label{sec:overview}


\subsection*{Root folder {\tt igabem2d/}}
The folder contains the three subfolders {\tt source/}, {\tt MEX-files/}, and {\tt functions/}.
With \texttt{mexIGABEM.m}, all required paths are added, the {\sc C} files of \texttt{source/} are compiled (via the {\sc Matlab}-{\sc C} interface {\sc mex}) and the resulting {\sc mex}-files are stored in \texttt{MEX-files/}. 
(To use {\sc mex}, you first have to set it up by typing {\tt mex -setup} in {\sc MATLAB}'s command window.)
Note that the contained main files \texttt{IGABEMWeak.m} and \texttt{IGABEMHyp.m}, where the adaptive Algorithm~\ref{the algorithm} and~\ref{the algorithm2} are implemented,  automatically call \texttt{mexIGABEM.m}.
In these main files, the user can choose out of a variety of different parameters, in particular different Dirichlet data $u$ and Neumann data $\phi$, which are provided in \texttt{source/Examples.h},  as well as different NURBS geometries as in Section~\ref{sec:boundary parametrization}, which are provided in \texttt{functions/Geometries.m}. 
If the exact solution of the corresponding Laplace problem is known, it is given in \texttt{functions/Examples.m} together with its derivative, which allows for additional visualization in the main files \texttt{IGABEMWeak.m} and \texttt{IGABEMHyp.m}. 
If desired, the number of Gauss quadrature points for the involved boundary integrals can be changed in \texttt{functions/getConfigWeak.m} resp.\ \texttt{functions/getConfigHyp.m}; see also Section~\ref{sec:implementation}.

\subsection*{Subfolder {\tt igabem2d/source/}}
This folder contains all required {\sc C} files with corresponding header-files.
To call the latter in {\sc Matlab}, we use the {\sc Matlab}-{\sc C} interface {\sc mex}; see the comments on \texttt{../mexIGABEM.m}.
In \texttt{Structures.h}, structures for splines and quadrature are defined.
In \texttt{Functions.h}, elementary {\sc C} functions are implemented.
In \texttt{Splines.c} and \texttt{Splines.h}, the evaluation of NURBS and their derivatives is implemented.
The evaluations can be used in {\sc Matlab} by compiling (via {\sc mex}) the  files \texttt{evalNURBS.c} (to evaluate $\widehat R_{\coarse,i,p}$), \texttt{evalNURBSComb.c} (to evaluate $\widehat S_\coarse\in \widehat\SS^p(\widehat{\mathcal{K}}_\coarse,\WW_\coarse)$), \texttt{evalNURBSCombDeriv.c} (to evaluate the right-hand derivative $\widehat S_\coarse^{\prime_r}$ of $\widehat S_\coarse\in \widehat\SS^p(\widehat{\mathcal{K}}_\coarse,\WW_\coarse)$), \texttt{evalNURBSCurve.c} (to evaluate $\gamma$), \texttt{evalNURBSCurveDeriv.c} (to evaluate the right-hand derivative $\gamma^{\prime_r}$).

The main functions to build the Galerkin systems are \texttt{VMatrix.c}, \texttt{RHSVectorWeak.c}, \texttt{WMatrix.c}, \texttt{PhiApprox.c}, \texttt{RHSVectorHyp.c}. 
They can be used in {\sc Matlab} by compiling (via {\sc mex}) the  files \texttt{buildVMatrix.c} (to build the matrix in~\eqref{eq:Galerkin_weak}), \texttt{buildRHSVectorWeak.c} (to build the right-hand side vector in~\eqref{eq:Galerkin_weak}), \texttt{buildWMatrix.c} (to build the matrix in~\eqref{eq:Galerkin_hyper}), \texttt{buildPhiApprox.c} (to build $\Pi_\coarse \phi$), \texttt{buildRHSVectorHyp.c} (to build the right-hand side vector in~\eqref{eq:Galerkin_hyper}).

The estimators and oscillations are implemented in \texttt{HHEstWeak.c}, \texttt{FaerEstWeak.c}, \texttt{ResEstWeak.c}, \texttt{HHEstHyp.c}, \texttt{ResEstHyp.c}, \texttt{OscHyp.c}.
They can be used in {\sc Matlab} by compiling (via {\sc mex}) the files \texttt{buildHHEstWeak.c} (to build~$\eta_{\V,\rm hh2,\coarse}$), \texttt{buildFaerEstWeak.c} (to build~$\eta_{\V,\rm fae,\coarse}$), \texttt{buildResEstWeak.c} (to build~$\eta_{\V,\rm res,\coarse}$), \texttt{buildHHEstHyp.c} (to build~$\eta_{\W,\rm hh2,\coarse}$), \texttt{buildResEstHyp.c} (to build~$\eta_{\W,\rm hh2,\coarse}$), \texttt{buildOscHyp.c} (to build~$\osc_{\W,\coarse}$).
The residual estimators require the evaluation of the boundary integral operators applied to some function, which is implemented in \texttt{ResidualWeak.c} and \texttt{ResidualHyp.c}. 
The evaluations can be used in {\sc Matlab} by compiling (via {\sc mex}) the source files \texttt{evalV.c} (to evaluate $\V\Phi_\coarse$ for $\Phi_\coarse\in\XX_\coarse$), \texttt{evalRHSWeak.c} (to evaluate $\K u$), \texttt{evalW.c} (to evaluate $\W U_\coarse$ for $U_\coarse\in\YY_\coarse$), \texttt{evalRHSHyp.c} (to evaluate $\K' \Pi_\coarse\phi$).

Finally, \texttt{Examples.h} provides different Dirichlet data $u$ and Neumann data $\phi$ that can be used in the main files \texttt{../IGABEMWeak.m} and \texttt{../IGABEMHyp.m}.

\subsection*{Subfolder \texttt{igabem2d/MEX-files/}}
Initially empty, all {\sc MEX}-files are stored here; see also the comments on \texttt{../mexIGABEM.m}.

\subsection*{Subfolder \texttt{igabem2d/functions}}
This folder contains all required {\sc Matlab} files.
As already mentioned, it contains the files \texttt{Geometries.m}  (providing different NURBS geometries),  \texttt{Examples.m} (providing solutions of the Laplace problem along with their derivative), and \texttt{getConfigWeak.m}, \texttt{getConfigHyp.m} (providing the number of quadrature points).
The D\"orfler marking~\eqref{eq:Doerfler} together with the marking of Algorithm~\ref{alg:refinement} is  implemented in \texttt{markNodesElements.m}. 
The actual refinement steps of Algorithm~\ref{alg:refinement} are implemented in \texttt{increaseMult.m}, which increases the multiplicity of marked nodes, and \texttt{refineBoundaryMesh.m}, which bisects at least all marked elements. 
Other auxiliary functions are \texttt{BasisTrafo.m} (to transform coefficients corresponding to a coarse basis to coefficients corresponding to a finer basis), \texttt{ChebyNodes.m} (to compute Chebyshev nodes), \texttt{GaussData.m} (to compute standard Gauss nodes and weights), \texttt{IntMean.m} (to compute the  integral mean over $\Gamma$), \texttt{LagrDerivMatrix.m} (to compute derivatives of Lagrange polynomials), \texttt{LogGaussData.m} (to compute Gauss nodes and weights for the logarithmic weight function), \texttt{ShapeRegularity.m} (to compute $\widehat\kappa_\coarse$).



\bibliographystyle{alpha}
\bibliography{literature}

\end{document}